\documentclass[11pt,amstex]{article}  % include bezier curves
\usepackage{color}              % Need the color package
\usepackage{graphicx}
\usepackage[]{amsmath}
\usepackage{amssymb}
\usepackage{hyperref}
\usepackage{enumerate, amsfonts, amsthm, bbm}

\definecolor{MyDarkBlue}{rgb}{0,0.08,0.50}
\definecolor{BrickRed}{rgb}{0.65,0.08,0}

\hypersetup{
colorlinks=true,       % false: boxed links; true: colored links
    linkcolor=MyDarkBlue,          % color of internal links
    citecolor=BrickRed,        % color of links to bibliography
    filecolor=red,      % color of file links
    urlcolor=cyan           % color of external links
}

\newtheorem{Lemma}{Lemma}[section]
\newtheorem{Proposition}[Lemma]{Proposition}

\newtheorem{Theorem}[Lemma]{Theorem}

\newtheorem{Corollary}[Lemma]{Corollary}

\newcommand{\prob}{\mathbb{P}}
\newcommand{\bfw}{\boldsymbol{w}}

\newcommand{\EE}{\mathcal{E}}

\newcommand{\DD}{\mathcal{D}}
\newcommand{\FF}{\mathcal{F}}
\newcommand{\II}{\mathcal{I}}
\newcommand{\JJ}{\mathcal{J}}
\newcommand{\GG}{\mathcal{G}}

\newcommand{\MM}{\mathcal{M}}

\newcommand{\WW}{\mathcal{W}}
\renewcommand{\SS}{\mathcal{S}}

\newcommand{\Rbold}{{\mathbb{R}}}

\newcommand{\RR}{\mathcal{R}}
\newcommand{\expec}{\mathbb{E}}

\newcommand{\e}{{\mathrm e}}
\newcommand{\im}{{\mathrm i}}
\newcommand{\shift}{\!\!\!\!}
\renewcommand{\Re}{{\mathrm{Re}}}

\newcommand{\R}{{\mathbb R}}

\newcommand{\eqn}[1]{\begin{equation} #1 \end{equation}}
\newcommand{\eqan}[1]{\begin{align} #1 \end{align}}
\newcommand{\lbeq}[1]{\label{#1}}
\newcommand{\refeq}[1]{(\ref{#1})}
\newcommand{\sss}{\scriptscriptstyle}
\newcommand{\Var}{{\rm Var}}

\newcommand{\op}{o_{\sss \prob}}

\newcommand{\cluster}{{\cal C}}

\newcommand{\Cmax}{\cluster_{\rm max}}

\newcommand {\vep}{\varepsilon}
\newcommand\1{\mathbbm{1}}
\newcommand{\indic}[1]{\1_{\{#1\}}}
\newcommand{\indicwo}[1]{\1_{#1}}
\newcommand {\convd}{\stackrel{d}{\longrightarrow}}
\newcommand {\convp}{\stackrel{\sss {\mathbb P}}{\longrightarrow}}
\newcommand {\convpv}{\stackrel{\sss \widetilde {\mathbb P}_v}{\longrightarrow}}

\newcommand{\nn}{\nonumber}

\newcommand{\thetastar}{\theta^*}

\newcommand{\Fvartheta}{\Lambda}
\newcommand{\FSS}{G}

\newcommand{\bfwit}{\boldsymbol{w}}

\numberwithin{equation}{section}
\newcommand{\betam}{\tilde{\beta}}

\newcommand{\ii}{{\rm i}}
\definecolor{darkgreen}{rgb}{0,.4,0}
\definecolor{darkagenta}{rgb}{.5,0,.5}
\definecolor{darkred}{rgb}{1,0,0}%was 0.85
\definecolor{darkblue}{rgb}{0,0,.4}

\begin{document}

\author{
%Elie Aidekon\thanks{
%Laboratoire de Probabilit\'es et Mod\`eles Al\'eatoires,
%Universit\'e Paris 6, 4, place Jussieu, 75005 Paris, France. E-mail: {\tt elie.aidekon@upmc.fr}}
%\and
Remco van der Hofstad
\thanks{Department of Mathematics and
Computer Science, Eindhoven University of Technology, P.O.\ Box 513,
5600 MB Eindhoven, The Netherlands. E-mail:
{\tt
rhofstad@win.tue.nl}, {\tt j.s.h.v.leeuwaarden@tue.nl}}
\and
Sandra Kliem
\thanks{
Fakult\"at f\"ur Mathematik, Universit\"at Duisburg-Essen, Thea-Leymann-Str. 9, D-45127 Essen, Germany.
E-mail: {\tt sandra.kliem@uni-due.de}}
\and
Johan S.H. van Leeuwaarden$^{\ *}$
}

\title{Cluster tails for critical power-law inhomogeneous random graphs}

\maketitle

\begin{abstract}
Recently, the scaling limit of cluster sizes for critical inhomogeneous random graphs
of rank-1 type having finite variance but infinite third moment degrees was obtained \cite{BhaHofLee09b}.
It was proved that when the degrees
obey a power law with exponent $\tau\in (3,4)$, the sequence of clusters
ordered in decreasing size and multiplied through by $n^{-(\tau-2)/(\tau-1)}$
converges as $n\rightarrow \infty$ to a sequence of decreasing non-degenerate random variables.

Here, we study the tails of the limit of the rescaled largest cluster, i.e., the
probability that the scaling limit of the largest cluster takes
a large value $u$, as a function of $u$. This extends a related result of
Pittel \cite{Pitt01} for the Erd\H{o}s-R\'enyi random graph
to the setting of rank-1 inhomogeneous random graphs
with infinite third moment degrees. We make use of delicate large deviations
and weak convergence arguments.

%It is well known that
%removing the largest critical cluster in the Erd\H{o}s-R\'enyi random graph shifts
%the location in the critical window by an amount proportional to the size of the largest
%connected component. This explains that, to leading order, the
%probability of having two large clusters of size $u$ and $a u$ with $u$ large
%in the Erd\H{o}s-R\'enyi random graph
%is similar to the product of the probabilities of having a largest cluster
%of sizes $u$ and $a u$ respectively. In our inhomogeneous setting,
%the effect of removing a large cluster is far more dramatic, and
%thus the probability of having two large clusters is much smaller
%than the product of the probabilities of having a largest cluster of size $u$ and $a u$
%respectively.
\end{abstract}

\vspace{0.3in}


\noindent
{\bf MSC2000 subject classification.}
60C05, 05C80, 90B15.

%\listoftodos

%\tableofcontents{}

%%%%%%%%%%%%%%%%%%%%%%%%%%%%%%%

\section{Introduction}
\label{sec-int}
The Erd\H{o}s-R\'enyi random graph $G(n, p)$ on the vertex set $[n]:=\{1,\ldots,n\}$ is constructed by including each of the ${n\choose 2}$ possible edges with probability $p$, independently of all other edges.
Erd\H{o}s and R\'enyi discovered the double-jump phenomenon: The size of the largest component was shown to be, in probability, of order $\log n$, $n^{2/3}$, or $n$, depending on whether the average vertex degree was less than, close to, or more than one. In 1984 Bollob\'as \cite{Boll84b} and subsequently {\L}uczak \cite{Lucz90a} showed for the scaling window $p=(1+\lambda n^{-1/3})/n$, that the largest component is of the order $n^{2/3}$. Since then, the critical, or near-critical behavior of random graphs has received tremendous attention (see \cite{Aldo97,AloSpe00,Boll01,Durr06,JanLucRuc00}).
Let $(\cluster_{\sss(i)})_{i\geq 1}$ denote the connected components of $G(n,p)$, ordered in size, i.e.,~$|\Cmax|=|\cluster_{\sss(1)}|\geq |\cluster_{\sss(2)}|\geq\ldots$
Aldous \cite{Aldo97} proved the following result:

\begin{Theorem}[Aldous \cite{Aldo97}]
\label{theo:main}
For $p=(1+\lambda n^{-1/3})/n$,  $\lambda\in\mathbb{R}$ fixed, and $n\rightarrow\infty$,
\eqn{
    \lbeq{dasdsad}
    \left(|\cluster_{\sss(1)}|n^{-2/3},|\cluster_{\sss(2)}|n^{-2/3}, \ldots\right) \convd \left(\gamma_1(\lambda),\gamma_2(\lambda), \ldots\right),
}
where $\gamma_1(\lambda)>\gamma_2(\lambda)>\ldots$ are the ordered excursions of the reflected version of the process $(W^\lambda_t)_{t\geq 0} \equiv (W_t + \lambda t -t^2/2)_{t\geq 0}$ with $(W_t)_{t\geq 0}$ a standard Wiener process.
\end{Theorem}

Theorem \ref{theo:main} says that the ordered connected components in the critical Erd\H{o}s-R\'enyi random graph are described by the ordered excursions of the reflected version of $(W^\lambda_t)_{t\geq 0}$.
The strict inequalities between the scaling limits of the ordered cluster follows from the local limit theorem proved in \cite{HofKagMul09}. In \cite{Pitt01} an exact formula was derived for the distribution function of the limiting variable $\gamma_1(\lambda)$ (of the largest component) and various asymptotic results were obtained, including
  	\eqn{
	\lbeq{pittel1}
    	\prob(\gamma_1(\lambda)>u)=
	\frac{1}{\sqrt{2\pi}u^{3/2}}\e^{-\frac{1}{8}u(u-2\lambda)^2}(1+o(1)),
	\quad u\rightarrow\infty.
   	}
%(a marginally sharper version was obtained in \cite{HofJanLee10}).
%and reads
%	\eqan{
%    	\lbeq{323}
%	\prob(\gamma_1(\lambda)>u)=8\frac{u^{1/2}
%	\exp\left(-\tfrac{1}{8}u(u-2\lambda)^2\right)}{\sqrt{8\pi}(u-2\lambda)
%	(3u-2\lambda)}\left(1+O(u^{-3/2})\right),\quad u\rightarrow\infty.
%	}
%\todotemp{Constants do not match between \refeq{323} and \refeq{pittel01}!}
The result in \eqref{pittel1} gives sharp asymptotics for the largest component in the critical Erd\H{o}s-R\'enyi graph.
%Further, in \cite{HofKagMul09}, a local limit theorem was proved for the largest
%connected components or clusters in $G(n,p)$ within the critical window where
%$p=(1+\lambda n^{-1/3})/n$. A consequence of this result is that, for any $a\in (0,1)$ and
%as $u\rightarrow \infty$,
%	\eqn{
%	\lbeq{324}
%	\prob(\gamma_1(\lambda)>u, \gamma_2(\lambda)>au)
%	=\prob(\gamma_1(\lambda)>u)\prob(\gamma_1(\lambda-u)>a u)(1+o(1)).
%	}
%Together with \eqref{pittel1}, this yields sharp asymptotics for the probability of having at least
%\emph{two} large clusters. In particular, the leading order asymptotics $\e^{-u^3/8-(au)^3/8}$
%of $\prob(\gamma_1(\lambda)>u, \gamma_2(\lambda)>au)$ agrees with that of $\prob(\gamma_1(\lambda)>u)\prob(\gamma_1(\lambda)>au)$, while
%$\prob(\gamma_1(\lambda)>u, \gamma_2(\lambda)>au)$ is \emph{much} larger than
%$\prob(\gamma_1(\lambda)>(1+a)u)$.
It was rederived in \cite{HofJanLee10} by studying the excursions of the scaling limit of the exploration process that is used to describe the limits in Theorem \ref{theo:main}. In this paper, we
follow a similar path, but then for a class of inhomogeneous random graphs and its scaling limit, and extend
\eqref{pittel1} to this setting.
%For certain critical inhomogeneous random graphs, the the scaling limit were recently identified in \cite{BhaHofLee09a,BhaHofLee09b,Turo09}.

Several recent works have studied inhomogeneity in random graphs and how it
changes the critical nature. In our model,
the vertices have a weight associated to them,
and the weight of a vertex moderates its degree. Therefore, by choosing
these weights appropriately, we can generate random graphs with highly
variable degrees. For our class of random graphs,
it is shown in \cite[Theorem 1.1]{Hofs13c} that when the weights do not vary
too much, the critical behavior is similar to the one in the
Erd\H{o}s-R\'enyi random graph. See in particular the recent
works \cite{BhaHofLee09a,Turo09}, where it was shown that if the degrees have
finite \emph{third} moment, then the scaling limit for the largest critical
components in the critical window are essentially the same (up to a trivial
rescaling) as for the Erd\H{o}s-R\'enyi random graph in Theorem \ref{theo:main}.

When the degrees
have \emph{infinite} third moment, instead, it was shown in \cite[Theorem 1.2]{Hofs13c} that the {sizes} of the largest
critical clusters are quite different.
In \cite{BhaHofLee09b} scaling limits were obtained for the sizes of the largest components at criticality for rank-1 inhomogeneous random graphs with power-law degrees with power-law exponent $\tau\in (3,4)$. For $\tau\in (3,4)$, the degrees have finite
variance but infinite third moment. It was shown that the sizes of the largest components, rescaled by $n^{-(\tau-2)/(\tau-1)}$, converge to hitting times of a \emph{thinned L\'evy process}. The latter is
a special case of the general multiplicative coalescents studied by Aldous and Limic in \cite{Aldo97} and \cite{AldLim98}. We next discuss these results in more detail.

\subsection{Inhomogeneous random graphs}
\label{sec-mod}

In our random graph model, vertices have weights,
and the edges are independent, with edge probabilities
being approximately equal to the rescaled product of the weights
of the two end vertices of the edge. While there are many
different versions of such random graphs (see below),
it will be convenient for us to work
with the so-called Poissonian random graph or Norros-Reittu model
\cite{NorRei06}. To define the model, we consider the vertex set $[n]:=\{1,2,\ldots, n\}$ and suppose each vertex is assigned a {weight}, vertex
$i$ having weight $w_i$. Now, attach an edge between vertices $i$ and $j$ with probability
    \eqn{
    \lbeq{pij-NR}
    p_{ij} = 1-{\rm exp}\Big(-\frac{w_i w_j}{\ell_n}\Big), \quad \mbox{where}\quad \ell_n=\sum_{i\in[n]} w_i.
    }
Different edges are independent. In this model, the average degree
of vertex $i$ is close to $w_i$, thus incorporating inhomogeneity
in the model.

There are many adaptations of this model,
for which equivalent results hold. Indeed,
the model considered here is a special case of the
so-called \emph{rank-1 inhomogeneous random graph}
introduced in great generality by Bollob\'as,
Janson and Riordan \cite{BolJanRio07}. It is
asymptotically equivalent with many related models,
such as the \emph{random graph with given prescribed degrees}
or Chung-Lu model, where instead
    \eqn{
    p_{ij}=\max(w_iw_j/\ell_n, 1),
    }
and which has been studied intensively by Chung and Lu (see
\cite{ChuLu02a,ChuLu02b,ChuLu03, ChuLu06c, ChuLu06}).
A further adaptation is the \emph{generalized random graph}
introduced by Britton, Deijfen and Martin-L\"of in
\cite{BriDeiMar-Lof05}, for which
    \eqn{
    \lbeq{pij-GRG}
    p_{ij} = \frac{w_i w_j}{\ell_n+w_iw_j}.
    }
See Janson \cite{Jans08a} for conditions under which these random graphs
are \emph{{asymptotically} equivalent}, meaning that all events
have asymptotically equal probabilities. As discussed in more detail
in \cite[Section 1.3]{Hofs13c}, these conditions apply
in the setting to be studied in this paper. Therefore, all results
proved here also hold for these related rank-1 models.

Let the weight sequence $\bfw = (w_i)_{i\in[n]}$ be defined by
    \eqn{
    \lbeq{choicewi}
    w_i = [1-F]^{-1}(i/n),
    }
where $F$ is a distribution function on $[0,\infty)$
for which we assume that there exists
a $\tau\in (3,4)$ and $0<c_{\sss F}<\infty$ such that
    \eqn{
    \lbeq{F-bound-tau(3,4)b}
    \lim_{x\rightarrow \infty}x^{\tau-1}[1-F(x)]= c_{\sss F},
    }
and where $[1-F]^{-1}$ is the generalized inverse function of
$1-F$ defined, for $u\in (0,1)$, by
    \eqn{
    \lbeq{invverd}
    [1-F]^{-1}(u)=\inf \{ s\colon [1-F](s)\leq u\}.
    }
By convention, we set $[1-F]^{-1}(1)=0$.
%We often make use of the fact that, with $U$ uniform on $[0,1]$,
%the random variable $[1-F]^{-1}(U)$ has distribution function $F$.

For the setting in \refeq{pij-NR} and \refeq{choicewi}, by
\cite[Theorem 3.13]{BolJanRio07}, the
number of vertices with degree $k$, which we denote by $N_k$,
satisfies
    \eqn{
    \lbeq{N-k-conv}
    N_k/n\convp \expec\Big[\mathrm{e}^{-W} \frac{W^k}{k!}\Big], \qquad k\geq 0,
    }
where $\convp$ denotes convergence in probability, and where
$W$ has distribution function $F$ appearing in \refeq{choicewi}.
We recognize the limiting distribution as a
so-called \emph{mixed Poisson distribution with mixing distribution $F$},
i.e., conditionally on $W=w$, the distribution is
Poisson with mean $w$. As discussed in more detail in \cite{Hofs13c}, since a Poisson
random variable with large parameter $w$ is closely concentrated
around its mean $w$, the tail behavior of the degrees in our
random graph is close to that of the distribution $F$. As a result,
when \refeq{F-bound-tau(3,4)b} holds, and with $D_n$ the degree of a uniformly
chosen vertex in $[n]$, $\limsup_{n\rightarrow \infty} \expec[D^a_n]<\infty$
when $a<\tau-1$ and $\limsup_{n\rightarrow \infty} \expec[D^a_n]=\infty$
when $a\geq \tau-1$. In particular, the degree of a uniformly
chosen vertex in $[n]$ has finite second, but infinite third moment
when \refeq{F-bound-tau(3,4)b} holds with $\tau\in (3,4)$.

Under the key assumption in \refeq{F-bound-tau(3,4)b},
    \eqn{
    \lbeq{[1-F]inv-bd}
    [1-F]^{-1}(u)=\big(c_{\sss F}/u\big)^{1/(\tau-1)}(1+o(1)), \quad u\downarrow 0,
    }
and
the third moment of the degrees tends to infinity, i.e.,
with $W\sim F$, we have $\expec[W^3]=\infty$. Define
    \eqn{
    \lbeq{nu-def}
    \nu = \expec[W^2]/\expec[W],
    }
so that, again by \refeq{F-bound-tau(3,4)b}, $\nu<\infty$.
{Then, by \cite[Theorem 3.1]{BolJanRio07} (see also
\cite[Section 16.4]{BolJanRio07} for a detailed discussion
on rank-1 inhomogeneous random graphs,
of which our random graph is an example),}
when $\nu> 1$, there is one giant component of size proportional to $n$,
while all other components are of smaller size $o(n)$, and when
$\nu \leq 1$, the largest connected component contains a proportion of
vertices that converges to zero in probability. Thus, the critical value
of the model is $\nu=1$. The main goal of this paper is to investigate what
happens close to the critical point, i.e., when $\nu=1$.

With the definition of the weights {in \refeq{choicewi},}
we shall write $\GG_n^0(\bfw)$ {for} the
graph constructed with the probabilities in \eqref{pij-NR}, while, for any fixed $\lambda\in \Rbold$, we
shall write $\GG_n^\lambda(\bfw)$ when we use the weight sequence
    \eqn{
    \lbeq{w-lambda-def}
    \bfwit(\lambda)=(1+\lambda n^{-(\tau-3)/(\tau-1)})\bfw.
    }

We shall assume that $n$ is so large that $1+\lambda n^{-(\tau-3)/(\tau-1)}\geq 0$,
so that ${w_i(\lambda)\geq 0}$ for all $i\in [n]$.
When $\tau>4$, so that
$\expec[W^3]<\infty$, it was shown in \cite{BhaHofLee09a, Hofs13c, Turo09}  that the scaling limit of the random graphs studied here
are (apart from a trivial scaling constant) \emph{equal} to the scaling limit
of the ordered connected components in the Erd\H{o}s-R\'enyi random graph in Theorem \ref{theo:main}.
When $\tau \in (3,4)$ the situation is entirely different, as discussed next.

Throughout this paper, we make use of the following standard notation.
We let $\convd$ denote convergence in distribution, and
$\convp$ convergence in probability. For a sequence of random variables
$(X_n)_{n\geq 1}$, {we write $X_n=\op(b_n)$ when $|X_n|/b_n\convp0$ as $n\rightarrow\infty$.}
For a non-negative function $n\mapsto g(n)$,
we write $f(n)=O(g(n))$ when $|f(n)|/g(n)$ is uniformly bounded, and
$f(n)=o(g(n))$ when $\lim_{n\rightarrow \infty} f(n)/g(n)=0$.
Furthermore, we write {$f(n)=\Theta(g(n))$ if $f(n)=O(g(n))$ and $g(n)=O(f(n))$}.
%We write that a sequence of events $(E_n)_{n\geq 1}$
%occurs \emph{with high probability} ({\bf whp}) when $\prob(E_n)\rightarrow 1$.
Finally, we abbreviate
    \eqn{
    \lbeq{defs-al-rho-eta}
    \alpha=1/(\tau-1),\qquad
    \rho=(\tau-2)/(\tau-1),
    \qquad
    \eta=(\tau-3)/(\tau-1).
    }

\subsection{The scaling limit for $\tau\in (3,4)$}
\label{sec-tau(3,4)}

We next recall two key results that we recently established in \cite{BhaHofLee09b}.

\begin{Theorem}[Weak convergence of the ordered critical clusters for $\tau\in (3,4)$ \cite{BhaHofLee09b}]
\label{thm-WC-(3,4)}
Fix the Norros-Reittu random graph with weights
$\bfw(\lambda)$ defined in {\rm {\refeq{choicewi}}} and {\rm \refeq{w-lambda-def}}.
Assume that $\nu=1$ and that {\rm \refeq{F-bound-tau(3,4)b}} holds. Then,
for all $\lambda\in \Rbold$,
    \eqn{
    \lbeq{largest-cluster-size-nota}
    \big(|\cluster_{\sss(1)}|n^{-\rho},|\cluster_{\sss(2)}|n^{-\rho},\ldots\big)
    \convd (\gamma_1(\lambda),\gamma_2(\lambda),\ldots),
    }
in the product topology, for some non-degenerate limit $(\gamma_i(\lambda))_{i\geq 1}$.
\end{Theorem}

In order to further specify the scaling limit $(\gamma_i(\lambda))_{i\geq 1}$, we need to introduce a non-negative
continuous-time process $(\SS_t)_{t\geq 0}$, referred to as  a {\it thinned L\'{e}vy process}, and defined as
    \eqn{
    \lbeq{SS-def}
    \SS_t=b -abt+ct+\sum_{i=2}^{\infty} \frac{b}{i^{\alpha}} \Big[\II_i(t)-\frac{at}{i^{\alpha}}\Big],
    }
where $a,b,c$ have been identified in \cite[Theorem 2.4]{BhaHofLee09b} as
$a=c_{\sss F}^{\alpha}/\expec[W]$,
$b=c_{\sss F}^{\alpha}$ and $c=\theta=\lambda+ \zeta$ with $\zeta \in (-\infty,0)$ the constant given in \cite[(2.18)]{BhaHofLee09b}\footnote{There is a typo in \cite[Theorem 2.4]{BhaHofLee09b}, in which $c=\theta-ab$ should read $c=\theta=\lambda+ \zeta$.}. 
%\todo[inline]{S-? If we follow this ref, it should reed $c=\theta-ab$ with $\theta=\lambda+\zeta$. BUT, I thought this was a typo in the ref, compare (2.39) and (3.25) in this ref.
%Johan: Yes, typo indeed} 
Further, here we use the notation
    \eqn{
    \lbeq{II-proc}
    \II_i(t)=\indic{T_i\leq t},
    }
where $(T_i)_{i\geq 2}$ are independent exponential random variables with mean
    \eqn{
    \lbeq{Ti-def}
    \expec[T_i]=i^{\alpha}/a.
    }
Let $H_1(0)$ denote the first hitting time of $0$ of the process $(\SS_t)_{t\geq 0}$, i.e.,
    \eqn{
    \lbeq{H1-def}
    H_1(0)=\inf\{t\geq 0\colon \SS_t=0\},
    }
and $\cluster(1)$ the connected component to which vertex 1 (with the largest weight) belongs.
We recall from \cite[Theorem 2.1 and Proposition 3.7]{BhaHofLee09b} %respectively \cite[(2.14)]{BhaHofLee09b}
%\todo[inline]{S-? check (2.14), I think you meant s.th. else. I would add Prop. 3.7 or s.th. else as we need to know that $H_1(0)$ is indeed as in (1.18) above. RvdH: I have just removed (2.14), Theorem 2.1 is what we need. SK: really? Discuss 3.}
that also $|\cluster(1)|n^{-\rho}$ converges in distribution:

\begin{Theorem}[Weak convergence of the cluster of vertex 1 for $\tau\in (3,4)$]
\label{thm-WC-C(1)}
Fix the Norros-Reittu random graph with weights
$\bfw(\lambda)$ defined in {\rm\refeq{choicewi} and \refeq{w-lambda-def}}.
Assume that $\nu=1$ and that {\rm \refeq{F-bound-tau(3,4)b}} holds. Then,
for all $\lambda\in \Rbold$,
    \eqn{
    n^{-\rho}|\cluster(1)|\convd H_1^{a}(0),
    }
with $H_1^{a}(0)$ the hitting time of $0$ of $(\SS_t)_{t\geq 0}$ with $a=c_{\sss F}^{\alpha}/\expec[W]$,
$b=c_{\sss F}^{\alpha}$, $c=\theta$.
%\todo[inline]{S-? See above, $c=\theta-ab$?! Remco: To do Johan 1 (once again)}
\end{Theorem}

By scaling, %\todo[inline]{Did we accomodate the scale-factor in our result? Sandra: I don't think so. Remco: This also plays in the LD paper... Remco: To do Johan 1 (once again)}
$H_1^{a}(0)/a$ for some $a,b,c$ has the same distribution as the hitting time 
$H_1(0)$ obtained by taking $b'=a'=1$, and $c'=c/(ab)=(\lambda+\zeta)/(ab)$.
%\todo[inline]{S-? See above, $c=\theta-ab$?! Remco: To do Johan 1 (once again)} Therefore, in what follows,
We shall reparametrize $a'=b'=1$ and let
    \eqn{
    \lbeq{SS-def-ref}
    \SS_t=1+\betam t+\sum_{i=2}^{\infty} c_i [\II_i(t)-c_i t],
    }
%\todo[inline]{S after discussion: This is the notation I chose here and in \cite[(1.3)]{AidHofKliLee13b}. This is in accordance with \cite[(2.41)]{BhaHofLee09b}. Related to To do Johan 1. Please look at this.}
where we set
    \eqn{
    \lbeq{beta-def}
    \betam=\beta-1\qquad \mbox{ with } \qquad \beta=c'=\theta/(ab)=(\lambda+ \zeta)/(ab),
    }
used the notation
		\eqn{
		\lbeq{def-ci}
		c_i = i^{-\alpha},
		}
and where $\II_i(t)$ is defined in \refeq{II-proc}--\refeq{Ti-def}
with $a'=1$.

%
%Theorem \label{thm-WC-C(1)} crucial ingredient in  the proof is the investigation of the cluster tail
%of a \emph{single} cluster. Indeed, we shall show that, for $u$ large,
%the probability that $1\in \Cmax$ is overwhelmingly large.
%We next study the probability that $H_1(0)>u$ for some $u>0$ very large,
%where $H_1(0)$ is the weak limit of $n^{-(\tau-2)/(\tau-1)} |\cluster(1)|$. \todotemp{add definition of $|\cluster(1)|$}

\subsection{Main results}
\label{sec-prop-limit(3,4)}

In this section we state our three main theorems.
The first theorem concerns the probability that $H_1^a(0)>u$ for some $u>0$ large,
where $H_1^a(0)$ is the weak limit of $n^{-\rho} |\cluster(1)|$ identified in Theorem \ref{thm-WC-C(1)}.
This is achieved by investigating the hitting time $H_1(0)$ of 0 of the process $(\SS_t)_{t\geq 0}$ in 
\refeq{SS-def-ref}.

\begin{Theorem}[Tail behavior scaling limit cluster vertex 1 for $\tau\in (3,4)$]
\label{thm-tail-C1}
When $u\rightarrow \infty$, there exists $I>0$ independent of $\betam$
and $A=A(\betam)$ and $\kappa_{ij}(\betam)\in {\mathbb R}$
such that
    	\eqn{
    	\lbeq{eqn-tail-C11}
	\prob(H_1(0)>u)= \prob(H_1^a(0)>au)=\frac{A}{u^{(\tau-1)/2}}
	\e^{-Iu^{\tau-1}+u^{\tau-1}\sum_{i+j \geq 1} \kappa_{ij}u^{-i(\tau-2)-j(\tau-3)}}(1+o(1)).
	}
\end{Theorem}
The constants $I$, $A$ and $\kappa_{ij}$ are specified in Section \ref{sec-overview}.
By scaling, these constants only depend on $a,b$ through $c'=c/(ab)=(\lambda+\zeta)/(ab)$,
any other dependence disappears since the law of $H_1(0)$ only depends on $c'$.
Since $\tau\in (3,4)$, the sum over $i,j$ such that $i+j\geq 1$ is in fact \emph{finite},
as we can ignore all terms for which $\tau-1-i(\tau-2)-j(\tau-3)\leq 0$.

We can even go one step further and study the optimal trajectory the process
$t\mapsto \SS_t$ takes in order to achieve the unlikely event
that $H_1(0)>u$ when $u$ is large. In order to describe this trajectory, we need to
introduce some further notation. In the proof, it will be crucial to \emph{tilt} the
distribution, i.e., to investigate the measure $\tilde{\prob}$
with Radon-Nikodym derivative $\e^{\theta u \SS_{u}}/\expec[\e^{\theta u \SS_{u}}]$,
for some appropriately chosen $\theta$. The selection of an appropriate $\theta$ for the thinned
L\'evy process$(\SS_t)_{t\geq 0}$ is quite subtle, and has been the main topic of our paper \cite{AidHofKliLee13b}.  The main results from paper \cite{AidHofKliLee13b} are reported in Section \ref{sec-overview}, and will play an important role in the present analysis. We refer to below \eqref{thetastaru-def} for the definition of $\theta^*$ that appears in the description of the
optimal trajectory that is identified in the following
theorem:

\begin{Theorem}[Optimal trajectory]
\label{thm-path-prop}
For $p\in [0,1]$, define
	\eqn{
	\lbeq{def_I_E}
    	I_{\sss E}(p)=(\tau-1)\int_0^{\infty}
	\Big( \frac{ \e^{\theta^* v} ( 1 - \e^{-p v} ) }{ \e^{\theta^* v}(1- \e^{-v}) + \e^{-v} } - p v\Big)
	\frac{ dv }{ v^{\tau - 1} }
    	}
with $\theta^*$ as defined below \eqref{thetastaru-def}.
Then, for $u\rightarrow \infty$, for any $\vep>0$,
    	\eqn{
    	\prob\Big(\sup_{p\in[0,1]} |\SS_{pu}-u^{\tau-2}I_{\sss E}(p)|\leq u^{\tau-2}\vep \mid H_1(0)>u\Big)
	=1-o(1).
    	}
\end{Theorem}
%\todo[inline]{Remco: Do we wish to add a simulation of the process $\SS_t$ for $t\in[0,u]$
%conditionally on $H_1(0)>u$? To do Johan 2, or drop}

%\todo[inline]{Remco: Do we wish to add a plot of the function $\tau\mapsto I=-\sup_{\vartheta} \Fvartheta(\vartheta)$! To do Johan 2, or drop}

Combining Theorem \ref{thm-tail-C1} and Theorem \ref{thm-path-prop}, and showing that, for $u$ large,
the probability that $1\in \cluster_{\sss(1)}$ is overwhelmingly large, eventually leads to the following result for the largest cluster size (recall $\gamma_1(\lambda)$ from \eqref{largest-cluster-size-nota}):

\begin{Theorem}[Tail behavior scaling limit for $\tau\in (3,4)$]
\label{thm-tail-Cmax}
When $u\rightarrow \infty$, there exists $I>0$ independent of $\lambda$
and $A=A(\lambda), \kappa_{i,j}=\kappa_{i,j}(\lambda)$ such that
	\eqn{
	\label{tail333}
    	\prob(\gamma_1(\lambda)>a u)= \frac{A}{u^{(\tau-1)/2}}\e^{-Iu^{\tau-1}
	+u^{\tau-1}\sum_{i+j \geq 1} \kappa_{ij}u^{-i(\tau-2)-j(\tau-3)}}(1+o(1)).
	}
\end{Theorem}
\medskip

The constants $I$, $A$ and $\kappa_{ij}$ are equal to those in Theorem \ref{thm-tail-C1}.

\paragraph{Brownian motion on a parabola.}
Note that substituting $\tau=4$ into \eqref{tail333} yields \linebreak $\frac{A}{u^{3/2}}\e^{-Iu^{3}+\kappa_{01} u^2+(\kappa_{10}+\kappa_{02})u}(1+o(1))$,
which agrees with the result of
Pittel in \eqref{pittel1}. This suggests a smooth transition from the
case $\tau\in (3,4)$ to the case $\tau>4$. We next further explore this relation.

%\todotemp{Expand and reposition this remark}
Consider the process $(W^\lambda_t)_{t\geq 0}=(W_t + \lambda t -t^2/2)_{t\geq 0}$ with $(W_t)_{t\geq 0}$ a standard Wiener process as mentioned in Theorem \ref{theo:main}. We now apply the technique of exponential change of measure to this process. First note that the moment generating function of
$W_u^{\lambda}$ can be computed as
%\todo[inline]{S-add at this point, $\phi(u; \vartheta)$, $\phi(u)$ not defined yet (maybe just change $=$ to $:=$ below and mention that later similar quantities for our process. Also, the text ``the main term is then'' below is not clear in this context but only in the context of Section 2.1. Remco: I have updated this argument, and have added some more details. In fact, this helps to understand what we aim for next. Do we wish to give even more insight? Discuss 5}
    	\eqn{
    	\log\phi(u; \vartheta)
	\equiv\log\expec[\e^{\vartheta u W_u^\lambda}]
	=\vartheta u(\lambda u-\tfrac12 u^2+\tfrac12 \vartheta u^2)
    	}
and let $\thetastar_u$ be the solution of
	$\thetastar_u=\arg\min_{\vartheta} \log\phi(u; \vartheta)$,
which is given by
	\eqn{
	\lbeq{thetastar-BM}
	\thetastar_u=\frac{1}{2}-\frac{\lambda}{u}.
	}
The main term is
  	\eqn{
	\lbeq{}
	 \phi(u)=\phi(u; \thetastar_u)=\expec[\e^{\thetastar_u u W_u^\lambda}]
	=\e^{-\frac18 u^3+\frac12 \lambda u^2-\frac12 \lambda^2 u}=\e^{-\frac{1}{8}u(u-2\lambda)^2}.
   	}
Noting that
	\eqn{
	\prob(\gamma_1(\lambda)>u)\leq \prob(W_u^{\lambda}>0)
	\leq \prob(\e^{\thetastar_u u W_u^{\lambda}}>1)\leq
	\expec[\e^{\thetastar_u u W_u^\lambda}],
	}
we see that this upper bound agrees to leading order with the result of Pittel in \eqref{pittel1}. In order to derive the full asymtptotics in \eqref{pittel1}, one can define the measure 	
	\eqn{
	\widetilde \prob(E)=\phi(u)^{-1}\expec[\e^{\thetastar_u u W_u^\lambda} \indicwo{E}],
	}
rewrite
	\eqn{
	\prob(\gamma_1(\lambda)>u)
	= \phi(u) \widetilde \expec[\e^{\thetastar_u u W_u^\lambda} \indic{\gamma_1(\lambda)>u}],
	}
and then deduce the asymptotics of the latter expectation in full detail.
Our analysis will be based on this intuition, now applied to a more involved, so-called thinned L\'evy, stochastic process.

%We next study the probabiliy of having \emph{two} large clusters:
%
%\begin{Theorem}[Two large clusters]
%\label{thm-tail-C2-maxs}
%When $u\rightarrow \infty$, there exists $I(a)>0$ independent of $\lambda$ such that for every
%$a\in [0,1)$,
%	\eqn{
%	\label{tail334}
%    	\prob(\gamma_1(\lambda)>u, \gamma_2(\lambda)> au)\leq
%	\e^{-I(a)u^{\tau-1}(1+o(1))}.
%	}
%Here $I(a)>I(1+a^{\tau-1})$ and $I(a)>I\times (1+a)^{\tau-1}$ for every $a\in (0,1)$.
%\end{Theorem}
%\todotemp{Can we actually prove these bounds on $I(a)$?}
%
%The fact that $I(a)>I(1+a^{\tau-1})$ implies that, unlike in \eqref{324},
%the leading order asymptotics of $\prob(\gamma_1(\lambda)>u, \gamma_2(\lambda)> au)$
%is much smaller than that of $\prob(\gamma_1(\lambda)>u)\prob(\gamma_1(\lambda)> au)$.
%Further, since $I(a)>I\times (1+a)^{\tau-1}$, we have that
%$\prob(\gamma_1(\lambda)>u, \gamma_2(\lambda)> au)$ is much \emph{smaller}
%than $\prob(\gamma_1(\lambda)>(1+a)u)$, contrary to $G(n,p)$ as described in
%\eqref{324}.
%
%
%In the next section we give a detailed overview of the proofs of these above theorems.

%%%%%%%%%%%%%%%%%%%%%%%%%%%%%%%

\section{Overview of the proofs}
\label{sec-overview}
In this section, we give the overview of the proofs of Theorems \ref{thm-tail-C1}-\ref{thm-tail-Cmax}.
The point of departure for our proofs is the conjecture that $\prob(H_1(0)>u)\approx \prob(\SS_u>0)$ for large $u$. The event $\{H_1(0)>u\}$ obviously implies $\{\SS_u>0\}$, but because of the
strong downward drift of the process $(\SS_t)_{t\geq 0}$, it seems plausible that both events are roughly equivalent.

In \cite{AidHofKliLee13b} a detailed study was presented on the large deviations behavior of the process $(\SS_t)_{t\geq 0}$. Using exponential tilting of measure the following two theorems were proved.
\begin{Theorem}[Exact asymptotics tail $\SS_u$ {\cite[Theorem 1.1]{AidHofKliLee13b}}]
\label{thm-tail-Su}
There exists $I, D>0$ and $\kappa_{ij}\in {\mathbb R}$
such that, as $u\rightarrow \infty$,
    	\eqn{
    	\lbeq{eqn-tail-C1}
	\prob(\SS_u>0)= \frac{D}{u^{(\tau-1)/2}}
	\e^{-Iu^{\tau-1}+u^{\tau-1}\sum_{i+j \geq 1} \kappa_{ij}u^{-i(\tau-2)-j(\tau-3)}}(1+o(1)).
	}
\end{Theorem}	
\begin{Theorem}[Sample path large deviations {\cite[Theorem 1.2]{AidHofKliLee13b}}]
\label{thm-sample-path-LD}
There exists a function $p\mapsto I_{\sss E}(p)$ on $[0,1]$ such that, for any $\vep>0$ and $p\in [0,1)$,
    	\eqn{
    	\lbeq{eqn-sample-paths}
	\lim_{u\rightarrow \infty}
	\prob\big(\big|\SS_{pu}-u^{\tau-2}I_{\sss E}(p)|\leq \vep u^{\tau-2}\mid\SS_u>0)=1.
	}
\end{Theorem}	

In \cite{AidHofKliLee13b} it is explained that specific challenges arise in the identification of a tilted measure due to the power-law nature of $(\SS_t)_{t\geq 0}$. General principles prescribe that the tilt should follow from a variational problem, but in the case of $(\SS_t)_{t\geq 0}$ this involves a Riemann sum that is hard to control. In \cite{AidHofKliLee13b} this Riemann sum is approximated by its limiting integral, and it is proved that the tilt that follows from the corresponding approximate variational problem is sufficient to establish the large deviations results in Theorems \ref{thm-tail-Su} and \ref{thm-sample-path-LD}. Details about this tilted measure are presented in Subsection \ref{sec-tilting}.

It is clear that Theorems \ref{thm-tail-Su} and \ref{thm-sample-path-LD} for the event $\{\SS_u>0\}$ are the counterparts of Theorems  \ref{thm-tail-C1} and \ref{thm-path-prop} for $\{H_1(0)>u\}$. Let us now sketch how we make formal the conjecture that  $\prob(H_1(0)>u)\approx\prob(\SS_u>0)$ for large $u$. We show that $\prob(H_1(0)>u)$ has the same asymptotic behavior as $\prob(\SS_u>0)$ in \eqref{eqn-tail-C1}, with the same constants except for the constant $D$. Despite the similarity of this result, the proof method we shall use is entirely different from the exponential tilting in \cite{AidHofKliLee13b}. In order to establish the asymptotics for $\prob(H_1(0)>u)$, we establish
sample path large deviations, not conditioned on the event $\{\SS_u>0\}$, but on the event $\{H_1(0)>u\}$.
This is much harder, since
we have to investigate the probability that $\SS_t>0$ for \emph{all} $t\in [0,u]$.
In order to prove these strong sample-path properties, we first prove that $\SS_t$ is close to its expected value for a finite, but large, number of $t$'s, followed by a proof that the path cannot deviate much in the small time intervals between these times. Now here is our strategy for the proofs. We extend the conjecture $\prob(H_1(0)>u)\approx \prob(\SS_u>0)$ by a conjectured sample path behavior that says that, under the tilted measure, the typical sample path of  $(\SS_t)_{t\geq 0}$ that leads to the event $\{\SS_u>0\}$ remains positive and hence
implies $\{H_1(0)>u\}$. To be more specific, we divide up this likely sample path into three parts: the early part, the middle part, and the end part. Our proof consists of treating each of these parts separately. We shall prove consecutively that with high probability the process:\\
(i) does not cross zero in the initial part of the trajectory (`no early hits');\\
(ii) is high up in the state space in the middle part of the trajectory, while experiencing small fluctuations, and therefore does not hit zero (`no middle ground');\\
(iii) is forced to remain positive until the very end.\\
In the last step, we have to be very careful, and it is in this step that it will turn out that the constant $D$ arising in the asymptotics of $\prob(\SS_u>0)$ in \refeq{eqn-tail-C1} is {\em different} from the constant $A$ arising in the asymptotics of $\prob(H_1(0)>u)$ in \refeq{eqn-tail-C11}.

We next summarize the technique of exponential tilting developed in \cite{AidHofKliLee13b}  for the
thinned  L\'evy process  $(\SS_t)_{t\geq 0}$ with $\tau\in(3,4)$, which allows us to give more details about how we shall establish the conjectured sample path behavior for each of the three parts described above.

\subsection{Tilting and properties of the tilted process}
\label{sec-tilting}
All results presented in this subsection are proved in \cite{AidHofKliLee13b}.

\paragraph{Exponential tilting.}
We use the notion of exponential tilting of measure in order to
give a convenient description of the probability of interest as follows.
    \eqn{
    \lbeq{tilt-form}
    \prob(\SS_{u}>0)
    =\phi(u;\vartheta) \widetilde\expec_{\vartheta}[\e^{-\vartheta u \SS_u}\indic{\SS_{u}>0}],
    }
where $\vartheta$ is chosen later on.
%Here, for an interval $I\subseteq [0,\infty)$, we write
%    \eqn{
%    \{\SS_{I}>0\}=\{\SS_t>0 \ \forall t\in I\},
%    }
%and
We define the measure $\widetilde \prob_{\vartheta}$ with corresponding expectation
$\widetilde \expec_{\vartheta}$ by the equality, for every event $E$,
    \eqn{
    \lbeq{tildeP-def}
    \widetilde \prob_{\vartheta}(E)=\frac{1}{\phi(u;\vartheta)}\expec_{\vartheta}[\e^{\vartheta u \SS_u}\indicwo{E}],
    }
where the normalizing constant $\phi(u;\vartheta)$ is defined as
    \eqn{
    \lbeq{phi-def}
    \phi(u; \vartheta)=\expec[\e^{\vartheta u \SS_u}].
    }
Choosing a good $\vartheta$ is rather delicate, and we explain this
in more detail now. By the independence of the indicators $(\II_i(u))_{i\geq 2}$, we obtain
that
    \eqan{
    \lbeq{expectation-exp}
    \phi(u; \vartheta)&=\expec[\e^{\vartheta u \SS_u}]=\e^{\vartheta u(1+\betam u)}\prod_{i=2}^{\infty}\e^{-\vartheta u^2 c_i^2}
    \Big(\e^{- c_iu}
    +\e^{\vartheta c_iu}
    (1-\e^{-c_iu})\Big)\\
    &=\e^{\vartheta u(1+\betam u)} \e^{\sum_{i=2}^{\infty}f(i/u^{\tau-1};\vartheta)}\nn
    %\e^{\sum_{i=2}^{\infty}
    %\Big(
    %\log\big(1+\e^{-c_iu}(\e^{-\vartheta c_iu}-1)\big)+\vartheta c_iu-\vartheta u^2 c_i^2 \Big)}.\nn
    }
with
    \eqn{
    \lbeq{def-f-tail}
    f(x;\vartheta)=\log\big(1+\e^{-x^{-\alpha}}
    (\e^{-\vartheta x^{-\alpha}}-1)\big)+\vartheta x^{-\alpha}-\vartheta x^{-2\alpha}.
    }
The function $x\mapsto f(x;\vartheta)$ is integrable at $x=0$ and at $x=\infty$,
so  the above sum can be approximated by the integral
    \eqn{
    \lbeq{sum-to-int}
    \sum_{i=2}^{\infty}f(i/u^{\tau-1};\vartheta)=u^{\tau-1}
    \int_0^{\infty}f(x;\vartheta) dx+e_\vartheta(u)
    \equiv u^{\tau-1}\Fvartheta(\vartheta)+ e_\vartheta(u),
    }
for some error term $u\mapsto e_\vartheta(u)$ given by
	\eqn{\label{sdfasdf}
	e_\vartheta(u)=\vartheta \left\{ u [\zeta(\alpha)-1] - u^2[\zeta(2\alpha)-1] \right\} + o_\vartheta(1),
	}
where $\alpha=1/(\tau-1)$ and the Riemann zeta functions $\zeta(\cdot)$ defined as
	\eqn{
 	\lbeq{rz}
 	\zeta(s)=\lim_{N\to\infty}\Big\{\sum_{n=1}^N n^{-s}-\frac{N^{1-s}}{1-s}-
	\frac12 N^{-s}\Big\}, \quad \Re (s)>-1, s\neq 1,
 	}
%\todo[inline]{S-add check against LD-article; there was a todo as I suggested writing $\Re$ differently. RvdH: Need to check this! Johan, are you willing to do this?}
where $\Re(s)$ denotes the real part of $s\in {\mathbb C}$. Equation \refeq{rz} follows from Euler-Maclaurin summation \cite[p.~333]{Har49}.
The error term in \eqref{sdfasdf} converges to $0$ uniformly for $\vartheta$ in compact sets bounded away from zero.
This implies that
    \eqn{
    \phi(u; \vartheta)
    =\e^{u^{\tau-1}\Fvartheta(\vartheta)+\vartheta u(\zeta(\alpha)+(\betam-\zeta(2\alpha)+1) u)+o_\vartheta(1)}.
    }
Let $\thetastar_u$ be the solution of
	\eqn{
	\lbeq{thetastaru-def}
	\thetastar_u=\arg\min_{\vartheta} \Big[\Fvartheta(\vartheta)
	+\vartheta u^{2-\tau}(\zeta(\alpha)+(\betam-\zeta(2\alpha)+1) u)\Big]
	}
and let $\thetastar$ be the value of $\vartheta$ where $\vartheta\mapsto \Fvartheta(\vartheta)$ is minimal.
It is not hard to see that $I\equiv -\Fvartheta(\thetastar)>0$ and that $\thetastar$ is
unique. In \cite[Lemma 3.6]{AidHofKliLee13b}, we have seen that $\thetastar_u=\thetastar+o(1)$. Further, $\thetastar>0$ by \cite[Lemma 3.5]{AidHofKliLee13b}. Define $\phi(u)=\phi(u;\thetastar_u)$. The next result investigates the main term $\phi(u)$:

\begin{Proposition}[Asymptotics of main term {\cite[Proposition 2.1]{AidHofKliLee13b}}]
\label{prop-asy-phi}
As $u\rightarrow \infty$, and with \linebreak $I=-\min_{\vartheta\geq 0} \Fvartheta(\vartheta)>0$,
there exist $\kappa_{ij}\in {\mathbb R}$ such that
    \eqn{
    \lbeq{phi-asy}
    \phi(u)=\expec[\e^{\thetastar_u u \SS_u}]=\e^{-Iu^{\tau-1}+u^{\tau-1}\sum_{i+j \geq 1} \kappa_{ij}u^{-i(\tau-2)-j(\tau-3)}}(1+o(1)).
    }
\end{Proposition}

\paragraph{Properties of the process under the tilted measure.}
From now on, we will take $\vartheta=\thetastar_u$, and we define
$\widetilde \prob=\widetilde \prob_{\thetastar_u}$ with corresponding expectation
$\widetilde \expec=\widetilde \expec_{\thetastar_u}$. In what follows, we abbreviate
$\theta=\thetastar_u$. Under this new measure, the rare event of $\SS_u$ being positive becomes quite likely. To describe these results, let us introduce some notation. Recall from \refeq{def_I_E} that, for $p\in [0,1]$,
	\eqn{
    	\lbeq{def_I_E-rep}
    	I_{\sss E}(p)
	=(\tau-1)\int_0^{\infty}\Big( \frac{ \e^{\theta^* v} ( 1 - \e^{-p v} ) }
	{ \e^{\theta^* v}(1- \e^{-v}) + \e^{-v} } - p v\Big) \frac{ dv }{ v^{\tau - 1} },
    	}
where we take $\vartheta=\thetastar$, which turns out to be the limit of $\thetastar_u$ as $u\rightarrow \infty$ (see, e.g., \cite[Lemma 3.6]{AidHofKliLee13b}).
As we see in Theorem \ref{thm-sample-path-LD}, the function $p\mapsto I_{\sss E}(p)$ will serve to describe the asymptotic mean of the process $p\mapsto \SS_{pu}$ conditionally on $\SS_u > 0$. It is not hard to check that
	\eqn{
    	I_{\sss E}(0)=0,\qquad \text{ and }
    	I_{\sss E}(1)=0,
    	}
the latter by definition of $\thetastar$, since $0 = \Fvartheta'(\thetastar) = I_{\sss E}(1)$. Finally,
	\eqn{
    	I_{\sss E}(p)>0 \mbox{ for every } p\in(0,1),
	}
and
	\eqn{
    	I_{\sss E}'(0)>0 \mbox{ and } I_{\sss E}'(1)<0.
	}
%\todo{check against more detailed version in our other article.}

\begin{Lemma}[Expectation of $\SS_t$ {\cite[Lemma 2.2]{AidHofKliLee13b}}]
\label{lem-expec-S}
As $u\to\infty$,\\
{\rm(a)}  $\widetilde \expec[\SS_{t}] = u^{\tau-2}I_{\sss E}(t/u)+ O(1+t + t |\theta^*- \theta_u^*| u^{\tau-3})$ uniformly in $t\in [0,u]$. \\
{\rm(b)}  $\widetilde \expec[\SS_{t}-\SS_u] = u^{\tau-2}I_{\sss E}(t/u)+ O(u-t + u^{-1} + |\theta^*- \theta_u^*| u^{\tau-2})$ uniformly in $t\in [u/2,u]$.\\
{\rm(c)}  $\widetilde \expec[\SS_{t}-\SS_u] = u^{\tau-3}I_{\sss E}'(1)(t-u)(1+o(1))+ O(u^{-1})$ when $u-t=o(u)$.\\
{\rm(d)}  $u\widetilde \expec[\SS_u]=o(1)$ when $u\rightarrow \infty$.
\end{Lemma}
\medskip

We will also need some consequences of the asymptotic properties of $\widetilde \expec[\SS_{t}]$. This is stated in the following corollary:

\begin{Corollary}
\label{cor-asympt-mean-S}
As $u\rightarrow \infty$,
\begin{itemize}
\item[{\rm (a)}] $\widetilde{\expec}[\SS_{t}]\geq \underline{c} t u^{\tau-3}$ and $\widetilde{\expec}[\SS_{t}]\leq \overline{c} t u^{\tau-3}$ uniformly for $t\in [\vep,u/2]$, where $0<\underline{c}<\overline{c}<\infty$;
\item[{\rm (b)}] $\widetilde{\expec}[\SS_{u-t}-\SS_u]\geq \underline{c} t u^{\tau-3}$ and
$\widetilde{\expec}[\SS_{u-t}-\SS_u]\leq \overline{c} t u^{\tau-3}$ uniformly for $t\in [T u^{-(\tau-2)},u/2]$, where $0<\underline{c}<\overline{c}<\infty$;
\item[{\rm (c)}] $\widetilde{\expec}[\SS_{t}]=\widetilde{\expec}[\SS_{t_1}](1+o(1))$ for $t\in [t_1,t_2]$ and $t_1\in [\vep, u/2]$ and $t_2-t_1=O(u^{-(\tau-2)})$;
\item[{\rm (d)}] $\widetilde{\expec}[\SS_{t}]=\widetilde{\expec}[\SS_{t_1}](1+o_{\sss T}(1))$ for $t\in [t_1,t_{2}]$ and $t_1\in [u/2, u-Tu^{-(\tau-2)}]$, where $o_{\sss T}(1)$ denotes a quantity
$c(T,u)$ such that
$\lim_{T\rightarrow \infty} \limsup_{u\rightarrow \infty}c(T,u)=0$ and $t_2-t_1=O(u^{-(\tau-2)})$.
\end{itemize}
\end{Corollary}
%\todo[inline]{Remco: added the strict inequalities in part (a), need to check whether proof still works...}

\proof Part (a) for $t\in [\vep,\vep u]$ for $\vep>0$ sufficiently small follows from Lemma \ref{lem-expec-S}(a) together with the facts that $I_{\sss E}(0)=0$, $I_{\sss E}'(0)>0$, and that $1+t + t |\theta^*- \theta_u^*| u^{\tau-3}=o(t u^{\tau-3})$. The fact that $I_{\sss E}'(0)>0$ also implies that $\underline{c}$ can be taken to be strictly positive. For $t \in [\vep u, u/2]$, Part (a) follows from the fact that $I_{\sss E}(p)>0$
for all $p\in[\vep,1/2]$ and that $1+t + t |\theta^*- \theta_u^*| u^{\tau-3}=o(u^{\tau-2})$.

Part (b) follows as Part (a), now using Lemma \ref{lem-expec-S}(b) together with the fact that $I_{\sss E}(1)=0$, $I_{\sss E}'(1)<0$.

Part (c) follows from Lemma \ref{lem-expec-S}(a), by subtracting the two terms.
Note that the error term $O(1+t_1 + t_1 |\theta^*- \theta_u^*| u^{\tau-3})$ is $o(t_1 u^{\tau-3})$ since $t_1\geq \vep$,
while $\widetilde{\expec}[\SS_{t_1}]=\Theta(t_1 u^{\tau-3})$ by Part (a) of this corollary. Further, note that
	\eqn{
	u^{\tau-2}[I_{\sss E}(t/u)-I_{\sss E}(t_1/u)]=O(u^{\tau-2} \max_{p\in [0,1]} |I_{\sss E}'(p)| (t_2-t_1))
	=O(1),
	}
which is $o(1)\widetilde{\expec}[\SS_{t_1}]$.\\
Part (d) follows again from Lemma \ref{lem-expec-S}(a) by subtracting the two terms. Note again that the error term $O(1+t_1 + t_1 |\theta^*- \theta_u^*| u^{\tau-3})$ is $o(t_1 u^{\tau-3})$,
while $\widetilde{\expec}[\SS_{t_1}]=\Theta(t_1 u^{\tau-3})$ by part (b) of this corollary and Lemma \ref{lem-expec-S}(a). Further, note that
	\eqn{
	u^{\tau-2}[I_{\sss E}(t/u)-I_{\sss E}(t_1/u)]=O(u^{\tau-2} \max_{p\in [0,1]} |I_{\sss E}'(p)| (t_2-t_1))
	=O_T(1),
	}
which is $o_{\sss T}(1)\widetilde{\expec}[\SS_{t_1}]$.
\qed

\medskip
The next lemma concerns the variance of the process.
Define, for $p\in [0,1]$,
\eqn{
    \lbeq{IV-def}
     I_{\sss V}(p)= (\tau-1)\int_0^{\infty} \frac{ \e^{\thetastar v}(1 - \e^{-p v} ) }{ \e^{\thetastar v}(1- \e^{-v}) + \e^{-v} } \Big(1- \frac{ \e^{\thetastar v} ( 1 - \e^{-p v} ) }{ \e^{\thetastar v}(1- \e^{-v}) + \e^{-v} } \Big) \frac{ dv }{ v^{\tau - 2} }
    }
    and
	\eqan{
	J_{\sss V}(p)&=  (\tau-1)\int_0^\infty \frac{ \e^{\thetastar v}(  \e^{-p v} - \e^{-v}) }
	{ \e^{\thetastar v}(1- \e^{-v}) + \e^{-v} } \Big(1- \frac{ \e^{\thetastar v} ( \e^{-p v} - \e^{-v}) }
	{ \e^{\thetastar v}(1- \e^{-v}) + \e^{-v} } \Big) \frac{ dv }{ v^{\tau - 2} },\\
	G_{\sss V}(p) &= (\tau-1)\int_0^\infty \frac{ \e^{2\thetastar v}( 1-\e^{-p v})(\e^{-pv} - \e^{-v}) }
	{ (\e^{\thetastar v}(1- \e^{-v}) + \e^{-v})^2 }{dv \over v^{\tau-2}}.
	}
Again, it is not hard to see that
	\eqn{
	\lbeq{IV-domain}
    	0<I_{\sss V}(p)<\infty \mbox{ for every } p\in(0,1], \mbox{ while } I_{\sss V}(0)=0,
	}
and
	\eqn{
	\lbeq{JV-domain}
    	0<J_{\sss V}(p)<\infty \mbox{ for every } p\in[0,1), \mbox{ while } J_{\sss V}(1)=0.
	}
%\todo{check (here, I only changed editing)}

\begin{Lemma}[Covariance structure of $\SS_t$ {\cite[Lemma 2.3]{AidHofKliLee13b}}]
\label{lem-var-S-rep}
As $u\to\infty$,\\
{\rm(a)}  $\widetilde{{\rm Var}}[\SS_{t}] = u^{\tau -3}I_{\sss V}(t/u) + O(1 + t |\theta^* - \theta_u^*| u^{\tau-4})$ uniformly in $t\in[0,u]$. \\
{\rm(b)}  $\widetilde{{\rm Var}}[\SS_{t}-\SS_u] = u^{\tau -3}J_{\sss V}(t/u) + O((u-t)u^{-1} + (u-t) |\theta^* - \theta_u^*| u^{\tau-4})$ uniformly in $t\in[0,u]$.\\
{\rm(c)}  $\widetilde {\rm Cov}[\SS_t,\SS_u-\SS_t] = -u^{\tau-3}G_{\sss V}(t/u)+ O((u-t)u^{-1} + (u-t) |\theta^* - \theta_u^*| u^{\tau-4})$ uniformly in $t\in[0,u]$.
\end{Lemma}	
%We prove Lemma \ref{lem-expec-S} and Lemma \ref{lem-var-S-rep} in Section \ref{sec-prop-tilted}.
%\todo[inline]{S after discussion: Remco, you asked if the minus sign in (c) is correct. Please see my comment in \cite{AidHofKliLee13b}. Remco: Still to check. To do Johan 3.}

The next result bounds the Laplace transform of the couple $(\SS_t,\SS_u)$:
\begin{Proposition}[Joint moment generating function of $(\SS_t,\SS_u)$  {\cite[Proposition 2.4]{AidHofKliLee13b}}]
\label{prop-laplace}
(a) As  $u\to \infty$,
	\eqn{
	\widetilde \expec\Bigg[\e^{ \lambda \frac{\SS_{t}-\widetilde \expec[\SS_t]}{\sqrt{I_{\sss V}(t/u) u^{\tau-3}}}}\Bigg]
	= \e^{{1\over 2}\lambda^2+ \Theta},
	}
where $|\Theta|\le o_u(1)$ as $u\rightarrow \infty$ uniformly in $t\in[u/2,u]$ and $\lambda$ in a compact set.\\
(b) Fix $\vep>0$ small.  As  $u\to \infty$, for any $\lambda_1, \lambda_2 \in \mathbb{R}$,
	\eqn{
	\widetilde \expec\Bigg[\e^{ \lambda_1 {\SS_{t} -\widetilde \expec[\SS_t] \over
	\sqrt{I_{\sss V}(t/u) u^{\tau-3}}} + \lambda_2  {\SS_{u} - \SS_t - \widetilde \expec[\SS_u - \SS_t]
	\over \sqrt{J_{\sss V}(t/u)u^{\tau-3}}}}\Bigg]
	= \e^{{1\over 2}\lambda_1^2 + {1\over 2}\lambda_2^2 -
	\lambda_1\lambda_2{G_{\sss V}(t/u)\over I_{\sss V}(t/u)J_{\sss V}(t/u)} + \Theta},
	}
where $|\Theta|\le o_u(1) + O(t^{3(3-\tau)/2})$ uniformly in $t\in[\vep,u-u^{-(\tau-5/2)}]$ and $\lambda_1, \lambda_2$ in a compact set.
\end{Proposition}
\medskip

By Proposition \ref{prop-laplace} and the fact that $u\widetilde \expec[\SS_u]=o(1)$
(see \cite[Lemma 4.1]{AidHofKliLee13b}), $u^{-(\tau-3)/2}\SS_u$ converges to a normal distribution with mean $0$ and variance $I_{\sss V}(1)$. We next extend this intuition by proving that the density of
$\SS_{u}$ close to zero behaves like $\left(2\pi I_{\sss V}(1)\right)^{-1/2}u^{-(\tau-3)/2}$:

\begin{Proposition}[Density of $\SS_u$ near zero {\cite[Proposition 2.5]{AidHofKliLee13b}}]
\label{prop-dens-f}
Uniformly in  $s=o(u^{(\tau-3)/2})$, the density $\widetilde f_{\SS_u}$ of $\SS_u$ satisfies
    \eqn{
    \widetilde f_{\SS_u}(s)=B u^{-(\tau-3)/2}(1+o(1)),
    }
with $B= \left(2\pi I_{\sss V}(1)\right)^{-1/2}$ and  $I_{\sss V}(p)$ defined in  \eqref{IV-def}.
Moreover,
$\widetilde f_{\SS_t}(s)$ is uniformly bounded by a constant times $u^{-(\tau-3)/2}$ for all
$s, u$ and $t \in [u/2,u]$.
\end{Proposition}
\medskip

There are three more results from \cite{AidHofKliLee13b} that will be used in this paper. The first is a description of the distribution
of the indicator processes $(\II_i(t))_{t\geq 0}$ under the measure $\widetilde \prob$.
Since our indicator processes $(\II_i(t))_{t\geq 0}$ are \emph{independent},
this property also holds under the measure $\widetilde \prob$:

\begin{Lemma}[Indicator processes under the tilted measure {\cite[Lemma 4.2]{AidHofKliLee13b}}]
\label{lem-ind-tilt}
Under the measure $\widetilde \prob$, the distribution of the indicator processes
$(\II_i(t))_{t\geq 0}$ is that of independent indicator processes. More precisely,
    \eqn{
    \lbeq{II-proc-tilt}
    \II_i(t)=\indic{T_i\leq t},
    }
where $(T_i)_{i\geq 2}$ are independent random variables with distribution
    \eqn{
    \lbeq{Ti-def-tilt}
    \widetilde \prob(T_i\leq t) =
    \begin{cases}
    \frac{\e^{\theta c_iu}(1-\e^{-c_i t})}{\e^{\theta c_iu}(1-\e^{-c_iu})
    +\e^{-c_iu}} &\text{for } t \leq u;\\
    \frac{\e^{\theta c_iu}(1-\e^{-c_iu})
    +(\e^{-c_iu}-\e^{-c_i t})}{\e^{\theta c_iu}(1-\e^{-c_iu})
    +\e^{-c_iu}} &\text{for } t>u.
    \end{cases}
    }
\end{Lemma}
\medskip

The second lemma describes what happens to the variances for small $p$ or for $p$ close to 1:
\begin{Lemma}[Asymptotic variance near extremities {\cite[Lemma 4.3(b)]{AidHofKliLee13b}}]
\label{lem-mean-var-extreme}
%{\rm(a)} As $a\to 0$, $I_{\sss E}(a)=a I_{\sss E}'(0)(1+o(1))$, while
%$I_{\sss E}(a)= -(1-a)I_{\sss E}'(1)(1+o(1))$ when $a\to 1$.\\
%{\rm(b)}
As $p\to 1$,  $J_{\sss V}(p) = -(1-p)J_{\sss V}'(1)(1+o(1))$ with $J_{\sss V}'(1)<0$, while, as $p\to 0$,
	\eqn{
	\lbeq{I_V-a-near-zero}
	I_{\sss V}(p) = p^{\tau-3} \mathcal I_{\sss V} (1+o(1)),
	\qquad
	\mbox{with}
	\qquad
	\mathcal I_{\sss V}= (\tau-1) \int_0^\infty (1 - \e^{-y}  ) \e^{-y} \frac{ dy }{ y^{\tau - 2} }.
	}
Consequently, there exist $0<\underline{c}<\overline{c}<\infty$ such that, for every $p\in[0,\vep]$ with $\vep>0$ sufficiently small,
	\eqn{
	\lbeq{I_V-a-near-zero-rep}
	\underline{c} p^{\tau-3} \leq I_{\sss V}(p) \leq \overline{c} p^{\tau-3}.
	}
\end{Lemma}
%\todo[inline]{S to Remco: I added equation-numbers \eqref{IV-domain} and \eqref{JV-domain}. That may make your extended formulation unnecessary? See three refs to this Lemma resp. its equation below. Remco: I have added \refeq{I_V-a-near-zero-rep} explicitly, as this is what we need.}
%\todotemp{RvdH: Do we really need Lemma \ref{lem-mean-var-extreme}?}
%
%\medskip

We finally rely on the following corollary that allows us to compute sums that we will encounter frequently:
\begin{Corollary}[Replacing sums by integrals in general {\cite[Corollary 3.3]{AidHofKliLee13b}}]
\label{cor-exp-mom-cs}
For every $a\in \R, a > \tau-1$ and $b>0$, there exists a constant $c(a,b)$ such that
	\eqn{
	\sum_{i=2}^{\infty} c_i^a \e^{-b c_i u}=c(a,b) u^{\tau-a-1}(1+o(1)).
	}
\end{Corollary}

\subsection{No early hits and middle ground}
In this section, we prove that the tilted process is unlikely to hit 0 until a time that is very close to
$u$. We start by investigating the early hits.

\paragraph{No early hits.}
In this step, we prove that it is unlikely that the process hits zero
early on, i.e., in the first time interval $[0,\vep]$ for some
$\vep>0$ sufficiently small. In its statement, we write
$0\in\SS_{[0,t]}$ for the event that $\{\SS_s=0\}$ for some $s\in [0,t]$, so that
$\prob(H_1(0)>u)=\prob(0\not\in \SS_{[0,u]})$.
\begin{Lemma}[No early hits]
\label{lem-no-early-hits}
For every $u\in [0,\infty)$, as $\vep\downarrow 0$,
    \eqn{
    \prob(0\in\SS_{[0,\vep]},\SS_{u}>0)
    =o_\vep(1)\prob(\SS_{u}>0),
    }
where $o_\vep(1)$ denotes a function that converges to zero as $\vep\downarrow 0$,
uniformly in $u$.
\end{Lemma}

The proof of Lemma \ref{lem-no-early-hits} follows from a straightforward
application of the FKG-inequality for independent random variables
(see \cite{ForKasGin71}, or \cite[Theorem 2.4, p.\ 34]{Grim99}).
The standard versions of the FKG-inequality hold for independent
indicator random variables, and in our case we need it for independent exponentials.
It is not hard to prove that the FKG-inequality we need holds by an
approximation argument.

\proof We note that the process $(\SS_t)_{t\geq 0}$ is a
deterministic function of the exponential random variables $(T_i)_{i\geq 2}$
(recall \refeq{SS-def}, \refeq{II-proc} and \refeq{Ti-def}). Now,
the event $\{0\in\SS_{[0,\vep]}\}$ is \emph{increasing} in
terms of the random variables $(T_i)_{i\geq 2}$ (use that $\SS_t$ only has positive jumps). Here we say that an
event $A$ is increasing when, if $A$ occurs for a realization
$(t_i)_{i\geq 2}$ of $(T_i)_{i\geq 2}$, and if $(t_i')_{i\geq 2}$
is coordinatewise larger than $(t_i)_{i\geq 2}$, then $A$
also occurs for $(t_i')_{i\geq 2}$. Clearly, the event
$\{\SS_{u}>0\}$ is \emph{decreasing} (for a definition, change the role of $t_i$ and $t'_i$ in the definition of an increasing event), so that the FKG-inequality
implies that these events are negatively correlated:
    \eqn{
    \prob(0\in\SS_{[0,\vep]},\SS_{u}>0)
    \leq \prob(0\in\SS_{[0,\vep]})\prob(\SS_{u}>0).
    }
We conclude the proof by noting that $\prob(0\in\SS_{[0,\vep]})=o_\vep(1)$
independently of $u$.
\qed
\medskip

\noindent
The key to our proof of Theorem \ref{thm-tail-C1} will be to show that $\prob(H_1(0)>u)=\Theta(\prob(\SS_{u}>0))$,
so that Lemma \ref{lem-no-early-hits} and the known asymptotics of
$\prob(\SS_{u}>0)$ imply that
it is unlikely to have an early hit of zero.

\paragraph{No middle ground.}
By \refeq{tildeP-def} (recall that $\phi(u)=\phi(u;\theta)$ with $\theta=\thetastar_u$), Lemma \ref{lem-no-early-hits} and Theorem \ref{thm-tail-Su},
    \eqn{
    \lbeq{tilt-form-rep-a}
    \prob(H_1(0)>u)
    = \phi(u) \widetilde\expec[\e^{-\theta u \SS_u} \indic{\SS_{[0,u]}>0}]=\phi(u) \widetilde\expec[\e^{-\theta u \SS_u} \indic{\SS_{[\vep,u]}>0}] +\phi(u) u^{-(\tau-1)/2} o_{\vep}(1).
    }
For $M>0$ arbitrarily fixed, we split
    \eqn{
    \prob(H_1(0)>u)
    = \phi(u) \widetilde\expec[\e^{-\theta u \SS_u} \indic{\SS_{[\vep,u]}>0, \SS_u\in[0,M/u]}]
	+\phi(u) \widetilde\expec[\e^{-\theta u \SS_u} \indic{\SS_{[\vep,u]}>0, \SS_u>M/u}]+\phi(u)u^{-(\tau-1)/2}o_{\vep}(1).
    }
By Proposition \ref{prop-dens-f}, we can bound
	\eqan{
	\widetilde\expec[\e^{-\theta u \SS_u} \indic{\SS_{[\vep,u]}>0, \SS_u>M/u}]
	&\leq \widetilde\expec[\e^{-\theta u \SS_u} \indic{\SS_u>M/u}]\leq  \int_{M/u}^{\infty} \e^{-\theta u v}f_{\SS_u}(v)dv\\
	&\leq O(u^{-(\tau-3)/2})
	\int_{M/u}^{\infty} \e^{-\theta u v}dv=O(u^{-(\tau-1)/2})\e^{-\theta M}.\nn
	}
As a result, we arrive at
	\eqn{
    	\lbeq{tilt-form-rep-b}
    	\prob(H_1(0)>u)
    	= \phi(u) \widetilde\expec[\e^{-\theta u \SS_u} \indic{\SS_{[\vep,u]}>0, \SS_u\in[0,M/u]}]
	+\phi(u)u^{-(\tau-1)/2}o_{\sss M}(1)+\phi(u)u^{-(\tau-1)/2}o_{\vep}(1),
    }
where $o_{\sss M}(1)$ denotes a quantity $c(M,u)$ such that $\limsup_{M\rightarrow \infty} \limsup_{u\rightarrow \infty} c(M,u)=0$.

We continue to prove that the dominant contribution to the
expectation of the right-hand side of \refeq{tilt-form} originates from
paths that remain positive until time $u-t$ for $t=Tu^{-(\tau-2)}$.

\begin{Proposition}[No middle ground]
\label{prop-no-middle-ground}
Fix $\varepsilon>0$. For every $u\in [0,\infty)$ and $\vep, M>0$ fixed,
    \eqn{
    \widetilde\prob (0\in\SS_{[\vep,u-Tu^{-(\tau-2)}]}, \SS_u \in[0,M/u])
    \leq o_{\sss T}(1)u^{-(\tau-1)/2},
    }
where we recall that $o_{\sss T}(1)$ denotes a quantity
$c(T,u)$ such that
$\lim_{T\rightarrow \infty} \limsup_{u\rightarrow \infty}c(T,u)=0$.
\end{Proposition}
\medskip

We prove Proposition \ref{prop-no-middle-ground} in Section \ref{sec-no-middle}.
%In the proof, we also prove the following strong concentration properties of the process $(\SS_t)_{t\in [0,u]}$,
%needed in the proof of Theorem \ref{thm-path-prop}.
%
%\begin{Lemma}[Concentration of $\SS_t$]
%\label{lem-conc-S}
%Fix $\varepsilon>0$. As $u\to\infty$,
%    	\eqn{
%    	\widetilde\prob\Big(\sup_{a\in[0,1]} |\SS_{au}-u^{\tau-2}I_{\sss E}(a)|
%	>\vep u^{\tau-2}, \SS_u \in[0,M/u] \Big)=o(u^{-(\tau-1)/2}).
%   	}
%\end{Lemma}

%\todo[inline]{S-? how do we get $a \in \{0,1\}$? $a=0$ is trivial as $S_0=1$, ok. For $a=1$ we could use $S_u \in [0,M/u]$ with high enough probability; see todo S-xxx (*). RvdH: Indeed. Since we are assuming that $\SS_u\geq 0$, for $a=1$ the event is thus $\SS_u\geq \vep u^{\tau-2}$. This is exponentially unlikely by the Chernoff bound. We need to check whether this is explained in sufficient detail where the proof is given.
%Discuss 13.}
\medskip

By \refeq{tilt-form-rep-b} and Proposition \ref{prop-no-middle-ground},
    \eqn{
    \lbeq{tilt-form-rep}
    \prob(H_1(0)>u)
    = \phi(u) \widetilde\expec[\e^{-\theta u \SS_u} \indic{\SS_{[u-T u^{-(\tau-2)},u]}>0}]
    + \phi(u)u^{-(\tau-1)/2}\big[o_{\vep}(1)+o_{\sss M}(1)+o_{\sss T}(1)\big].
    }
Since $\vep$, $M$ and $T$ are arbitrary, it now suffices to identify the
asymptotics of the expectation appearing on the right-hand side of \refeq{tilt-form-rep}.

\subsection{Remaining positive near the end}
\label{sec-pos-end}
To prove Theorem \ref{thm-tail-C1}, by Proposition \ref{prop-asy-phi} and equation (\ref{tilt-form-rep}),
it suffices to prove that, with $\gamma=(\tau-1)/2$,
    \eqn{
    \lbeq{near-end}
    \widetilde\expec[\e^{-\theta u \SS_u}\indic{\SS_{[u-t,u]}>0, \SS_u \in[0,M/u]}]
    =(A+o_{\sss T}(1))u^{-\gamma}(1+o(1)),
    }
where $t=Tu^{-(\tau-2)}$.
In the above expectation, we see two terms. The term $\e^{-\theta u \SS_u}$
forces $\SS_u$ to be small, more precisely, $\SS_u=\Theta(1/u)$ for $u$ large,
while the term $\indic{\SS_{[u-t,u]}>0}$ forces the path to remain positive until
time $u$. We now study these two effects.

%%%%%%%%%%%%%%%%%%%%%%%%%%%%%%%%%%%%%%%%%%%%%%%

In order to investigate the probability that $\SS_{[u-t,u]}>0$,
we proceed as follows. Let
    \eqn{
    \lbeq{index-def}
    \JJ(u)=\{j\colon \II_j(u)=1\}
    }
denote the set of indices for which $T_j\leq u$. We condition on the
set $\JJ(u)$. Note that $\SS_u$ is measurable with respect to $\JJ(u)$. We now rewrite $\SS_{u-t}$ in a convenient form. For this, recall \eqref{SS-def-ref} and write
%that
%    \eqn{
%    \SS_t=1+\betam t+\sum_{i=2}^{\infty} c_i[\II_i(t)-c_it].
%    }
    \eqan{
    \SS_{u-t}&=\frac{t}{u}+\frac{u-t}{u}\SS_u +\sum_{i=2}^{\infty} c_i[\II_i(u-t)-\frac{u-t}{u}\II_i(u)]\\
    &=\frac{t}{u}+\frac{u-t}{u}\SS_u -\sum_{i=2}^{\infty} c_i[\indic{T_i\in (u-t,u]}-\frac{t}{u}\II_i(u)].\nn
    }
Thus, with
    \eqan{
    \lbeq{Qu-def}
    Q_u(t)
    &\equiv u \SS_{u-t} - t - (u-t) \SS_u
    =-\sum_{i=2}^{\infty} c_i[u\indic{T_i\in (u-t,u]}-t\II_i(u)],
    }
we have that $\SS_{u-t}>0$ precisely when $Q_u(t) > -t-(u-t)\SS_u.$
We rewrite
    \eqn{
    Q_u(t)=-\sum_{i\in \JJ(u)} c_i[u\indic{T_i\in (u-t,u]}-t].
    }
Note that, for any $t=o(u)$,
    \eqan{
    \lbeq{near_end-as-int}
    \widetilde\expec[\e^{-\theta u \SS_u}\indic{\SS_{[u-t,u]}>0}]
    &=
    \frac{1}{u}\int_0^{\infty} \e^{-\theta v}
    \widetilde\prob\big(\SS_{[u-t,u]}>0\mid u\SS_u=v\big) \widetilde f_{\SS_u}(v/u)dv\\
    &=
    \frac{1}{u}\int_0^{\infty} \e^{-\theta v}
    \widetilde\prob\big(Q_u(s)> -v-s+sv/u \,\forall s\in [0,t]\mid u\SS_u=v\big) \widetilde f_{\SS_u}(v/u)dv.
    \nn
    }
We aim to use dominated convergence on the above integral, and we start by proving pointwise convergence.
By Proposition \ref{prop-dens-f}, $\widetilde f_{\SS_u}(v/u)=B u^{-(\tau-3)/2}(1+o(1))$ pointwise in $v$
(in fact, even when $v=o(u^{(\tau-1)/2})$).
This leads us to study, for all $v>0$,
    \eqn{
    \lbeq{gu(v)-def}
    g_{u,t}(v)\equiv\widetilde\prob\big(Q_u(s)> -v-s+sv/u \,\forall s\in [0,t]\mid u\SS_u=v\big).
    }
%\todo[inline]{S-??? (*1) Shouldn't this be some $g_{u,t}(v)$? Then \refeq{gu(v)-conv} becomes $g_{u,t}(v) \rightarrow g_t(v) \equiv \prob\big(\sup_{0 \leq s \leq t} (L_s-\kappa s) \leq v\big)$. Actually, I get $\prob\big(\sup_{0 \leq s \leq t} (L_s-(\kappa+1) s) \leq v\big)$? Now fix \refeq{int-rewrite}-\eqref{A-def} ... I admit, I am confused. If you change this (or not), also have a look at todo (*2). Remco: I agree. Check the consequences of this change. Discuss 14.}
We split
    \eqn{
    \lbeq{split-Q}
    Q_u(t)=A_u(t)-B_u(t),
    }
where
    \eqn{
    \lbeq{AB-defs}
    B_u(t)\equiv u\sum_{i\in \JJ(u)} c_i[\indic{T_i\in (u-t,u]}-\widetilde \prob(T_i>u-t \mid T_i\leq u)],
    \qquad
    A_u(t)\equiv -\sum_{i\in \JJ(u)} c_i[u\widetilde \prob(T_i>u-t \mid T_i\leq u)-t].
    }
Thus, $(A_u(t))_{t\in [0,u]}$ is \emph{deterministic} given $\JJ(u)$,
while $(B_u(t))_{t\in [0,u]}$ is \emph{random} given $\JJ(u)$.
The main result for the near-end regime is the following proposition, which proves
that $g_u(v)$ converges pointwise. % (the proof is deferred to Section \ref{section-pf-A-B}).

\begin{Proposition}[Weak conditional convergence of time-reversed process]
\label{prop-weak-cond-conv}
~\\
(a) As $u\rightarrow \infty$, conditionally on $u\SS_{u}=v$,
    \eqn{
    \lbeq{aim-A}
    (A_u(t u^{-(\tau-2)}))_{t\geq 0}\convd (\kappa t)_{t\geq 0},
    }
where $\kappa\in (0,\infty)$ is given by
    \eqn{
    \lbeq{kappa-def}
    \kappa=\int_0^{\infty} x^{-\alpha} \frac{\e^{\theta x^{-\alpha}}\e^{-x^{-\alpha}}}
    {\e^{\theta x^{-\alpha}}(1-\e^{-x^{-\alpha}})+\e^{-x^{-\alpha}}}
    \big[\e^{x^{-\alpha}}-1-x^{-\alpha} \big]dx.
    }
(b) As $u\rightarrow \infty$, conditionally on $u\SS_{u}=v$,
    \eqn{
    \lbeq{aim-M}
    (B_u(tu^{-(\tau-2)}))_{t\geq 0}\convd (L_t)_{t\geq 0},
    }
where $(-L_t)_{t\geq 0}$ is a L\'evy process with no positive jumps and with Laplace transform
    \eqn{
    \lbeq{Levy-def}
    \expec[\e^{a (-L_s)}]=\e^{s \int_{-\infty}^0 (\e^{a z}-1-a z)\Pi(dz)}, \quad a \geq 0,
    }
and characteristic measure
    \eqn{
    \lbeq{charact-meas-def}
    \Pi(dz)=(\tau-1)\frac{(-z)^{-(\tau-1)}\e^{-\theta z}}{\e^{-\theta z}(1-\e^z)
    +\e^z} \e^z dz.
    }
\end{Proposition}
\medskip

Proposition \ref{prop-weak-cond-conv} is proved in Section \ref{section-pf-A-B},
and determines the precise constant $A$ from \eqref{eqn-tail-C11}, as we now explain
in more detail.

We proceed by investigating some properties of the supremum of the L\'{e}vy process
from \eqref{aim-M} that we need later on. Note in particular that the distribution of $L_s$ in \eqref{Levy-def}
does not depend on $v$. With a slight abuse of notation, also the distribution of the
limiting process $(L_s)_{s\geq 0}$ shall
be denoted by $\prob$.

\begin{Lemma}[Supremum of the L\'{e}vy process]
\label{lem-sup-levy}

Let $I_\infty \equiv \inf_{t \geq 0} (-L_t+\kappa t)$. Then
    \eqn{
    \lbeq{I-infty-W}
    \prob(I_\infty \geq -v) = \WW(v)/\WW(\infty),
    }
where $\WW\colon[0,\infty) \rightarrow [0,\infty)$ is the unique continuous increasing function that has Laplace transform
    \eqn{
    \lbeq{W-eq-def}
    \int_0^\infty \e^{-a x} \WW(x) dx = \frac{1}{\psi(a)}, \quad a>\Psi(0),
    }
where the Laplace exponent $\psi$ is given by $\expec[\e^{a (\kappa t-L_t)}]=\e^{t\psi(a)}$ and is
computed in \eqref{laplace-exp} below, while $\Psi(0)$ is the largest solution of the
equation $\psi(a)=0$, and $\WW(\infty)=1/\psi'(0)=1/\kappa$ is a constant.
\end{Lemma}

\proof We rewrite \eqref{Levy-def} to see that $X_s \equiv -L_s+\kappa s$ is a
L\'evy process with no positive jumps and Laplace exponent
    \eqan{
    \lbeq{laplace-exp}
    \psi(a) &= \kappa a + \int_{(-\infty,0)} \left(\e^{az}-1-az \right) \Pi(dz) \\
    &= \beta' a + \int_{(-\infty,0)} \left(\e^{az}-1-az \indic{z>-1} \right) \Pi(dz) \nn
    }
    with
        \eqan{
     \beta' &= \kappa + \int_{(-\infty,0)} (-z) \indic{z \leq -1} \Pi(dz) > 0
    }
as defined in \cite[Section VII.1]{Bert96}. Indeed, recall from \cite[Section VII.1]{Bert96}
that $\expec[\e^{aX_s}] = \e^{s \psi(a)}$ and note that our $\beta'$ corresponds to $a$ in \cite{Bert96}. Also note from \eqref{kappa-def} that $\kappa>0$. Thus $\psi'(0+) = \kappa>0$ and \cite[Corollary 2(ii) in Section VII.1]{Bert96} yields that $X_s$ drifts to $\infty$ (for a definition, see \cite[Theorem 12(ii) in Section VI.3]{Bert96}). This in turn implies (see \cite[Proof of Theorem 8, in Section VII.2]{Bert96})
    \eqn{
    \prob(I_\infty \geq -v) = \WW(v)/\WW(\infty),
    }
where $\WW$ is given in the statement of \cite[Theorem 8, in Section VII.2]{Bert96}. For the definition of $\Psi$ see before \cite[Theorem 1 of Section VII.1]{Bert96}. Also note from the second equation of the proof of \cite[Proof of Theorem 8, in Section VII.2]{Bert96} that $\WW(\infty)>0$. To see that $\WW(\infty)=1/\psi(0)$, note that if $a\downarrow 0$,
    \eqn{
    \int_0^\infty \e^{-a x} \WW(x) dx= \frac{\WW(\infty)}{a (1+o(1))}.
    }
Now, $\psi(0)=0$, so that $1/\psi(a)=1/(a\psi'(0))(1+o(1))$ as $a\downarrow 0$, which identifies $\WW(\infty)=1/\psi'(0)=1/\kappa$.
\qed
\medskip

By Proposition \ref{prop-weak-cond-conv} and the continuity of $\WW$ in Lemma \ref{lem-sup-levy}, with $\MM_T=\sup_{0\leq s\leq T} (L_s-\kappa s)$
for each $v\geq 0$ and for $t=Tu^{-(\tau-2)}$, for $u \rightarrow \infty$,
%\todo[inline]{S: check the two points I added/changed}
    	\eqn{
    	\lbeq{gu(v)-conv}
    	g_{u,t}(v)\rightarrow g_T(v)\equiv \prob\big(\MM_T \leq v\big).
	}
Further, as $T\rightarrow \infty$,
	\eqn{
	\lbeq{gu(v)-conv-rep}
	g_T(v)\downarrow g(v)\equiv \prob(\sup_{0\leq s < \infty} (L_s-\kappa s)\leq v)=\frac{\WW(v)}{\WW(\infty)}.
    	}
Now we are ready to complete the proofs of our main results.

\subsection{Completion of the proofs}
\label{sec-compl-proofs}

\paragraph{Completion of the proof of Theorem \ref{thm-tail-C1}.}
We start by completing the proof of Theorem \ref{thm-tail-C1}. Recall that it remains to prove \eqref{near-end} with $\gamma=(\tau-1)/2$. By \eqref{near_end-as-int} and \eqref{gu(v)-def}, we need to compute
    \eqn{
	\label{int-rewrite}
    \widetilde\expec[\e^{-\theta u \SS_u}\indic{\SS_{[u-t,u]}>0, \SS_u\in[0,M/u]}]
    =
    \frac{1}{u^{(\tau-1)/2}}\int_0^{M} \e^{-\theta v}
    g_{u,t}(v)\big[u^{(\tau-3)/2}\widetilde f_{\SS_u}(v/u)\big]dv,
    }
where $t=Tu^{-(\tau-2)}$. A similar problem was encountered in \cite[Proof of Theorem 1.1]{AidHofKliLee13b}, which is restated here as Theorem \ref{thm-tail-Su}, apart from the
fact that there the function
$g_{u,t}(v)$ was absent.

We wish to use bounded convergence. For this, we note that
$u^{(\tau-3)/2}\widetilde f_{\SS_u}(v/u)\rightarrow B$ by Proposition \ref{prop-dens-f} for each $v$
(in fact, for all $v=o(u)$), while, by \refeq{gu(v)-conv}--\refeq{gu(v)-conv-rep}, $g_{u,t}(v)\rightarrow g_T(v)$, which, in turn, converges to $g(v)$ as $T\rightarrow \infty$.
Further, since $g_{u,t}(v)\leq 1$ and $u^{(\tau-3)/2}\widetilde f_{\SS_u}(v/u)$ is uniformly bounded (see Proposition \ref{prop-dens-f}), the integrand $\e^{-\theta v}g_{u,t}(v)\big[u^{(\tau-3)/2}\widetilde f_{\SS_u}(v/u)\big]$ is uniformly bounded by a constant.
Thus, by the Bounded Convergence Theorem,
    \eqan{
    \lbeq{compl-conv}
    \widetilde\expec[\e^{-\theta u \SS_u}\indic{\SS_{[u-t,u]}>0, \SS_u\in[0,M/u]}]
    &=\frac{B}{u^{(\tau-1)/2}}\int_0^M \e^{-\theta v}
    g_T(v)dv(1+o(1))\\
	&=\frac{B}{u^{(\tau-1)/2}}\int_0^M \e^{-\theta v}
    g(v)dv(1+o(1)+o_{\sss T}(1)).\nn
    }
This identifies (recall \eqref{near-end}, \eqref{I-infty-W}, \eqref{W-eq-def} and \eqref{gu(v)-conv})
    \eqn{
    \lbeq{A-def}
    A=B\int_0^{\infty} \e^{-\theta v}
    g(v)dv
    =B\int_0^{\infty} \e^{-\theta v}
    \prob\big(\MM\leq v\big)dv=\frac{B}{\theta} \expec[\e^{-\theta \MM}]=\frac{B\psi'(0)}{\psi(\theta)}.
    }
Since $D=B/\theta$ by \cite[(7.4)]{AidHofKliLee13b} and $\prob\big(\MM\leq v\big)<1$ for every $v$, we also immediately obtain that $A\in(0,D)$.
%where $M=\sup_{s\geq 0} (L_s-\kappa s)$. To see that this is well-defined, see Lemma \ref{lem-sup-levy} above.
%We further compute, using \refeq{W-eq-def},
%    \eqn{
%    A=B\int_0^{\infty} \e^{-\theta v}
%    g(v)dv=B \int_0^{\infty} \e^{-\theta v} W(v)/W(\infty)
%    =B/(\psi(\theta)W(\infty))=B\psi'(0)/\psi(\theta).
%    }
%\todotemp{Is $\psi(\theta)$ a nice constant?}
%Alternatively, by reversing the integral, we
%can rewrite, using $f_{\sss M}$ to denote the density of $M$,
%    \eqn{
%    \int_0^{\infty} \e^{-\theta v}
%    \prob\big(M\leq v\big)dv
%    =\int_0^{\infty} \e^{-\theta v}
%    \int_{0}^v f_{\sss M}(m) dm dv
%    =\frac{1}{\theta} \int_0^{\infty} \e^{-\theta m}
%    f_{\sss M}(m) dm=\frac{1}{\theta} \expec[\e^{-\theta M}].
%    }
%This identifies
%    \eqn{
%    \lbeq{A-comp}
%    A=\frac{B}{\theta} \expec[\e^{-\theta M}],
%    }
and completes the proof of Theorem \ref{thm-tail-C1}.
\qed

\paragraph{Path properties: Proof of Theorem \ref{thm-path-prop}.}
%Theorem \ref{thm-path-prop} follows directly from Theorem \ref{thm-tail-C1},
%which can be reformulated as
%$\widetilde \prob(H_1(0)>u)=\Theta(u^{-(\tau-1)/2})$ and from Lemma \ref{lem-conc-S}.
We bound, using that $\{H_1(0)>u\}\subseteq \{\SS_u>0\}$,
	\eqan{
	&\prob\Big(\sup_{p\in[0,1]} |\SS_{pu}-u^{\tau-2}I_{\sss E}(p)| >\vep u^{\tau-2}\mid H_1(0)>u\Big)\\
	%&\qquad=\frac{\widetilde \expec\Big[\indic{\sup_{a\in[0,1]}
	%|\SS_{au}-u^{\tau-2}I_{\sss E}(a)| > \vep u^{\tau-2}, H_1(0)>u}
	%\e^{-\theta u\SS_u}\Big]}
	%{\widetilde \expec\Big[\indic{H_1(0)>u}\e^{-\theta u\SS_u}\Big]}\nn\\
	&\qquad\leq \prob\Big(\sup_{p\in[0,1]} |\SS_{pu}-u^{\tau-2}I_{\sss E}(p)| >\vep u^{\tau-2}\mid\SS_u>0\Big)\frac{\prob(\SS_u>0)}{\prob(H_1(0)>u)}.\nn
	}
By Theorems \ref{thm-tail-Su} and \ref{thm-tail-C1}, the ratio of probabilities converges to $D/A\in(0,\infty)$, while, by Theorem \ref{thm-sample-path-LD}, the conditional probability converges to 0.
This completes the proof of Theorem \ref{thm-path-prop}.
\qed

\paragraph{Completion of the proof of Theorem \ref{thm-tail-Cmax}.}
%\todo[inline]{Johan: this is really the last part. Last time we went through it, but I find it too risky to make the changes. Your notes are in you mail box.} 

We finally complete the proof of Theorem \ref{thm-tail-Cmax} using
Theorem \ref{thm-tail-C1} and recalling \eqref{beta-def}. Denote
    \eqn{
    \lbeq{SS(i)-def-ref}
    \SS_t^{\sss (i)}=c_i+\betam_i t+\sum_{j=1\colon j\neq i}^{\infty} c_j [\II_j(t)-c_j t],
    }
%\todo[inline]{S-??? I gave the ref. Formula (3.76) in the ref reads differently. I think there are 3 typos inside:
%	$$
%	\SS_t^{\sss (i)} = b \ch{c_i} - abt c_i^\ch{2} + ct
%	+ \sum_{j=\ch{2}: j \neq i}^\infty b c_j [\II_j(t) - atc_j]
%	$$
%Remco: Still to be done. To do Johan 4.
%}
where $\betam_i=(\lambda+\zeta)/(ab)-c_i^2$ (see \cite[Remark 3.9]{BhaHofLee09b} and recall $a',b',c'$ from above \eqref{SS-def-ref}). The intuition for the above formula is that
    \eqn{
    \lbeq{SS(i)-def-ref-rep}
    \SS_t^{\sss (i)}=\sum_{j\geq 1}^{\infty} c_j [\II_j(t)-c_j t],
    }
where we slightly abuse notation to set $\II_i(0)=1$ for the process $(\SS_t^{\sss (i)})_{t\geq 0}$. Since $(\SS_t^{\sss (i)})_{t\geq 0}$ describes the scaling limit of the exploration process of the cluster of vertex $i\geq 1$, while $\II_j(t)$ has the interpretation as the indicator that vertex $j$ is found in the exploration before time $t$, it is reasonable to set $\II_i(0)=1$ for $(\SS_t^{\sss (i)})_{t\geq 0}$. \footnote{We take this opportunity to correct some typos in \cite{BhaHofLee09b}.
In \cite[(3.76)]{BhaHofLee09b}, the term $-abti^{-\alpha}=-abtc_i$ should be replaced by
$-abti^{-2\alpha}=-abtc_i^2.$ This corresponds to the choice of $\betam_i=(\lambda+\zeta)/(ab)-c_i^2$ here. Further, in \cite[(3.79)]{BhaHofLee09b}, the product over $j\in[i-1]$ should be over $q\in[i-1]$.
}

Define
    \eqn{
    \lbeq{H(i)-def}
    H^{\sss (i)}(0)=\inf\{t\geq 0\colon \SS_t^{\sss (i)}=0\}.
    }
Then, $H_1(0)=H^{\sss (1)}(0)$. Let $\cluster(i)$  be the connected component to which vertex $i$ belongs, and let $\cluster_{\sss \leq}(i)$ be the set  $\cluster(i)$ if none of the vertices $j\in[i-1]=\{1,\ldots,i-1\}$ belongs to $\cluster(i)$, and the empty set $\varnothing$ otherwise. We know from \cite[(3.78)]{BhaHofLee09b} and the scaling explained around \eqref{SS-def-ref} that $n^{-\rho}|\cluster(i)|\convd a \cdot H^{\sss (i)}(0)$ for each $i\geq 1$ with $\rho=(\tau-2)/(\tau-1)$ (cf.\ \eqref{defs-al-rho-eta}).
Finally, denote
    \eqn{
    H_i(0)=\begin{cases}
    0 &\text{if }\exists j<i \text{ such that } \II_j(H^{\sss (i)}(0))=1;\\
    H^{\sss (i)}(0) &\text{otherwise.}
    \end{cases}
    }
Then, by \cite[(3.79)]{BhaHofLee09b}, $n^{-\rho}|\cluster_{\sss \leq}(i)| \convd a \cdot H_i(0)$. This provides us with the appropriate background to complete the proof of Theorem \ref{thm-tail-Cmax}.

We start with the lower bound. By construction,
$\gamma_1(\lambda)\geq a \cdot H_1(0)$ (see \cite[Theorems 1.1 and 2.1]{BhaHofLee09b} and recall that $\cluster_{\sss(i)}$ denotes the $i^{\mathrm{th}}$-largest connected component).
Therefore,
    \eqn{
    \prob(\gamma_1(\lambda)>au)\geq \prob(H_1(0)>u),
    }
and thus the lower bound follows from Theorem \ref{thm-tail-C1}.

For the upper bound, we use that (cf.\ \cite[Theorems 1.1]{BhaHofLee09b})
    \eqn{
    \prob(\gamma_1(\lambda)>au)
	=
	\lim_{n\rightarrow\infty}
	\prob(\exists i\colon n^{-\rho}|\cluster_{\sss \leq}(i)|\geq au)\leq
	\lim_{n\rightarrow\infty} \sum_{i\geq 1}
	\prob(n^{-\rho}|\cluster_{\sss \leq}(i)|\geq au).
	}
By the weak convergence of $ n^{-\rho}|\cluster_{\sss \leq}(i)|$ and the fact that
there are with high probability only finitely many clusters that are larger than $\vep n^{\rho}$
(as proved in \cite[Theorem 1.6]{BhaHofLee09b}),
%\todotemp{I erased the middle part in \eqref{thm1-6-a}:}
	\eqn{
	\lbeq{thm1-6-a}
	\prob(\gamma_1(\lambda)>au)\leq
%	\sum_{i\geq 1} \prob(H^{\sss (i)}(0)>\ch{au})=
    	\prob(H_1(0)>u) +\sum_{i\geq 2} \prob(H_i(0)>u).
   	}
The first term is the main term, and we prove that
$\sum_{i\geq 2} \prob(H_i(0)>u)=o(\prob(H_1(0)>u))$ now.

For this, we note that
    \eqan{
    \prob(H_i(0)>u)&=\prob(\II_j(u)=0 \, \forall j\in [i-1], \SS_{[0,u]}^{\sss (i)}>0)\\
	&=\prob(\II_j(u)=0 \, \forall j\in [i-1])\prob(\SS_{[0,u]}^{\sss (i)}>0\mid \II_j(u)=0 \, \forall j\in [i-1]).\nn
    }
We can rewrite, on the event $\{\II_j(u)=0 \, \forall j\in [i-1]\}$,
	\eqn{
	\SS_{t}^{\sss (i)}=\tfrac{(\lambda+\zeta)}{ab}t+\sum_{j\geq i+1} c_j(\II_j(t)-c_jt)
	+c_i-\sum_{j=1}^i c_j^2t\leq 
	\tfrac{(\lambda+\zeta)}{ab}t+\sum_{j\geq i+1} c_j(\II_j(t)-c_jt)
	+c_1-\sum_{j=1}^i c_j^2t.
	}
Therefore,
	\eqan{
	\prob(\SS_{[0,u]}^{\sss (i)}>0\mid \II_j(u)=0 \, \forall j\in [i-1]) 
	&\leq \prob\Big(\tfrac{(\lambda+\zeta)}{ab}t+\sum_{j\geq i+1} c_j(\II_j(t)-c_jt)
	+c_1-\sum_{j=1}^i c_j^2t>0~\forall t\in[0,u]\Big) \\
	&=\prob\big(\SS_{[0,u]}^{\sss (1)}>0\mid \II_j(u)=0 \, \forall j\in [i]\setminus \{1\}\big).\nn
	}
	
The event $\big\{\II_j(u)=0\forall j\in [i]\setminus \{1\}\big\}$ is decreasing (recall the notions used in the proof of Lemma \ref{lem-no-early-hits}) in the
random variables $(T_i)_{i\geq 2}$, while the event
$\{\SS_{[0,u]}^{\sss (1)}>0\}$ is increasing. Thus, by the FKG-inequality,
    	\eqn{
    	\lbeq{thm1-6-aa}
    	\prob\big(\SS_{[0,u]}^{\sss (1)}>0\mid \II_j(u)=0 \, \forall j\in [i]\setminus \{1\}\big)
	\leq 	\prob(\SS_{[0,u]}^{\sss (1)}>0)=\prob(H_1(0)>u).
    	}
We can identify
    \eqn{
    \lbeq{thm1-6-b}
    \prob(\II_j(u)=0\forall j\in [i-1])
    =\e^{-\sum_{j=1}^{i-1} c_j u}.
    }
%\ch{Further, we next show that $\SS_t^{\sss (i)}$ is stochastically dominated by $\SS_t^{\sss (1)}$, which we abbreviate as $\SS_t^{\sss (i)}\preceq\SS_t^{\sss (1)}$. For this, we write
%	\eqn{
%	\SS_t^{\sss (i)}=\SS_t^{\sss (1,i)}+c_1[\II_1(t)-1],
%	\qquad\qquad
%	\SS_t^{\sss (1)}=\SS_t^{\sss (1,i)}+c_i[\II_i(t)-1],
%	}
%where $\SS_t^{\sss (1,i)}$ denotes the common parts of $\SS_t^{\sss (1)}$ and $\SS_t^{\sss (i)}$, which is independent from $(\II_1(t),\II_i(t))$. Therefore, $\SS_t^{\sss (i)}\preceq\SS_t^{\sss (1)}$ follows if 
%$c_1[\II_1(t)-1]\preceq c_i[\II_i(t)-1]$. For this, we compute, for $x<0$ and using that $c_1\geq c_i$,
%	\eqn{
%	\prob(c_1[\II_1(t)-1]\leq x)=\prob(\II_1(t)=0)=1-\e^{-c_1t}
%	\geq 1-\e^{-c_i t}=\prob(c_i[\II_i(t)-1]\leq x).
%	}
%For $x\geq 0$, both sides are equal to 1, so that $\prob(c_1[\II_1(t)-1]\leq x)\geq \prob(c_i[\II_i(t)-1]\leq x)$ holds for every $x\in {\mathbb R}$, which proves that $c_1[\II_1(t)-1]\preceq c_i[\II_i(t)-1]$. 
%%This stochastic domination holds for every fixed $t\geq 0$, but since $\II_i(t)=\indic{T_i\leq t}$, it follows that $T_i\preceq . As a result, we can couple  $\SS_t^{\sss (i)}$ and $\SS_t^{\sss (1)}$ such that $\SS_t^{\sss (i)}\leq \SS_t^{\sss (1)}$ occurs a.s. 
%Therefore,
%    \eqn{
%    \lbeq{thm1-6-c}
%    \prob(\SS_{[0,u]}^{\sss (i)}>0)\leq \prob(\SS_u^{\sss (1)}>0) \leq \prob(\SS_u^{\sss (1)}>0).
%    }
Combining \eqref{thm1-6-a}, \eqref{thm1-6-aa}--\eqref{thm1-6-b} we arrive at
    \eqn{
    \prob(\gamma_1(\lambda)>au)
    \leq \prob(H_1(0)>u)\Big[1+\sum_{i\geq 2}\e^{-\sum_{j=1}^{i-1} c_j u}\Big].
    }
Since $c_j=j^{-\alpha}$ with $\alpha\in (1/3,1/2)$,
$\sum_{j=1}^{i-1} c_j\geq (i-1) c_{i-1}=(i-1)^{1-\alpha}$.
Therefore,
    \eqn{
    \sum_{i\geq 2}\e^{-\sum_{j=1}^{i-1} c_j u}=o(1).
    }
%Further, by Theorems \ref{thm-tail-C1} and \ref{thm-tail-Su}, $\prob(\SS_u^{\sss (1)}>0)/\prob(H_1(0)>u)$ remains uniformly bounded. 
This completes the proof of Theorem \ref{thm-tail-Cmax}.
\qed

\section{No middle ground: Proof of Proposition \ref{prop-no-middle-ground}}
\label{sec-no-middle}

In this section, we show that the probability to hit zero in the time interval $[\varepsilon,u-Tu^{-(\tau-2)}]$, where $T$ is a constant, becomes negligible as $T\rightarrow \infty$.

The strategy of proof is as follows. We start in Proposition \ref{prop-integer-times} by investigating
the value of $\SS_t$ at some discrete times $(t_k)_{k\geq 1}$ in $[0,u]$
and show that with high probability $\SS_t$ does not deviate far from its mean.
Next, in Proposition \ref{prop-no-hit-small-interval}, we show
that it is unlikely for the process $(\SS_t)_{t\geq 0}$
to make a substantial deviation in the interval $[t_k,t_{k+1}]$
from its value in $t_k$.

We start with a preparatory lemma that will allow us to give bounds on the asymptotic
parameters appearing in the upcoming proofs:

\begin{Lemma}[Asymptotics of parameters]
\label{lem-LD-par}
There exists $K\geq 1$ such that
	\eqn{
	\lbeq{first-bd-sec-mom}
	\widetilde\prob\big(\sum_{i=2}^{\infty} c_i^2 \II_i(u)\geq K u^{\tau-3}\big)
	\leq Cu^{-(\tau-1)},
	}
and, for all $|\lambda| \leq \delta u$ with $\delta>0$ sufficiently small, there exists $K>0$
such that
	\eqn{
	\lbeq{sec-bd-sec-mom}
	\widetilde\prob\Big(\sum_{i=2}^{\infty} c_i [1-\II_i(u/2)]\big(\e^{\lambda c_i}-1-\lambda c_i\big)
	\geq K\lambda^2u^{\tau-4}\Big)
	\leq C u^{-(\tau-1)}.
	}
\end{Lemma}

\proof
We use the second moment method. With Lemma \ref{lem-ind-tilt} we compute that
	\eqn{
	\widetilde\expec\big[\sum_{i=2}^{\infty} c_i^2 \II_i(u)\big]
	\leq \sum_{i=2}^{\infty} c_i^2 C(\theta) \big( 1-\e^{-c_i u} \big).
	}
Split the sum into $i$ with  $c_i u \leq 1$ and $c_i u>1$. For the first, we bound
$1-\e^{-c_i u}\leq O(1) c_iu$, for the latter, we bound $1-\e^{-c_i u}\leq 1$, to obtain
	\eqn{
	\widetilde\expec\big[\sum_{i=2}^{\infty} c_i^2 \II_i(u)\big]
	\leq O(1) \sum_{i\colon c_iu\leq 1} c_i^3u
	+O(1) \sum_{i\colon c_iu>1} c_i^2= O(1) u^{\tau-3}(1+o(1)),
	}
the latter by an explicit computation using that $c_i=i^{1/(\tau-1)}$.

Further, with Corollary \ref{cor-exp-mom-cs}
	\eqn{
	\widetilde\Var\big(\sum_{i=2}^{\infty} c_i^2 \II_i(u)\big)
	\leq \sum_{i=2}^{\infty} c_i^4 (1-\widetilde \prob(T_i\leq u))
	\leq C(\theta) \sum_{i=2}^{\infty} c_i^4\e^{-c_i u} = O(1)u^{\tau-5}(1+o(1)).
	}
The Chebychev inequality now proves \refeq{first-bd-sec-mom}.

For \refeq{sec-bd-sec-mom}, we again compute
	\eqn{
	\widetilde\expec\Big[\sum_{i=2}^{\infty} c_i [1-\II_i(u/2)]\big(\e^{\lambda c_i}-1-\lambda c_i\big)\Big]
	\leq C(\theta) \sum_{i=2}^{\infty} c_i \e^{-c_iu/2}[\e^{\lambda c_i}-1-\lambda c_i]
	= C(\theta) \sum_{i=2}^{\infty} c_i \e^{-c_iu/2}\e^{|\lambda| c_i}(\lambda c_i)^2/2.
	}
Thus, for $|\lambda|\leq \delta u$ and again using
Corollary \ref{cor-exp-mom-cs}, we obtain
	\eqn{
	\widetilde\expec\Big[\sum_{i=2}^{\infty} c_i [1-\II_i(u/2)]\big(\e^{\lambda c_i}-1-\lambda c_i\big)\Big]
	=O(\lambda^2 u^{\tau-4}).
	}
Further,
	\eqan{
	\widetilde\Var\Big(\sum_{i=2}^{\infty} c_i [1-\II_i(u/2)]\big(\e^{\lambda c_i}-1-\lambda c_i\big)\Big)
	&\leq C(\theta) \sum_{i=2}^{\infty} c_i^2\e^{-c_i u/2}\big(\e^{\lambda c_i}-1-\lambda c_i\big)^2\\
	&\leq
	C(\theta) |\lambda|^4  \sum_{i=2}^{\infty} c_i^6\e^{-c_i u/2}\e^{2|\lambda| c_i}
	=O(|\lambda|^4 u^{\tau-7}).\nn
	}
Again the claim follows from the Chebychev inequality.	
\qed
\bigskip

\medskip

We continue to show that the probability for $\SS_t$ to deviate far from its mean at some \emph{discrete} times in the time interval $[\vep,u-Tu^{-(\tau-2)}]$ is small when $T$ is large enough:

\begin{Proposition}[Probability to deviate far from mean at discrete times]
\label{prop-integer-times}
Let $\eta>0$ and $\delta_u=u^{-(\tau-2)}$. For any $\vep>0$ and $M>0$,
	\eqn{
	\limsup_{u\to\infty} u^{(\tau-1)/2}\widetilde \prob \left(\exists \, k\in\mathbb{N} \mbox{ s.t. }
	k\delta_u \in [\vep,u-T\delta_u]\colon \,
	\big|\SS_{k\delta_u} - \widetilde \expec [\SS_{k\delta_u}]\big| >\eta
	\widetilde \expec [\SS_{k\delta_u}], \SS_u\in[0,M/u]\right)
	= o_{\sss T}(1),
	}
where we recall the definition of $o_{\sss T}(1)$ from {\rm Proposition \ref{prop-no-middle-ground}}.
\end{Proposition}

%\todotemp{Remco: Replace $m$ in $t\in [\vep,u/2]$ by $\vep$??}

\proof The proof is split between the cases $t\in [\vep,u/2]$,   $t\in [u/2, u-\vep]$ and $t\in [u-\vep, u-u^{-(\tau-2)}]$, where $\vep>0$ is some arbitrary constant.

\paragraph{{\bf Proof for $t\in [\vep,u/2]$.}}We start by proving the proposition for $t\in [\vep,u/2]$, for which we use Proposition \ref{prop-laplace} with $\lambda_1=\pm 1$ and $\lambda_2=0$ to see that, for any $x>0$,
	\eqn{
	\widetilde \prob\Big( \big|{\SS_{t} - \widetilde \expec[\SS_{t}]
	\over \sqrt{I_{\sss V}(t/u) u^{\tau-3}}}\big| >x\Big)
	\le \centering c\e^{-x},
	}
where we note that the $\e^{\Theta}$ error term can be put inside the constant $c$ since $|\Theta|\leq o_u(1)+O(t^{3(\tau-3)})$  and $t\geq \vep$ is strictly positive. By \eqref{I_V-a-near-zero-rep} in Lemma \ref{lem-mean-var-extreme}, $I_{\sss V}(p) \le c p^{\tau-3}$ for all $p\in [0,1/2]$. Applying this to $p=t/u$ yields
	\eqn{
	\widetilde \prob\left( |\SS_t  -\widetilde \expec[\SS_{t}] | > c x t^{(\tau-3)/2} \right) \le c\e^{-x}.
	}
%Since $t\leq u/2$, in turn, this shows that
%	\eqn{
%	\widetilde \prob\left( |\SS_t  -\widetilde \expec[\SS_{t}] | >
%	c x t^{\frac{\tau-3}{2}}\right) \le c(\vep)\e^{-x}.
%	}
By Corollary \ref{cor-asympt-mean-S}(a), we have $\widetilde \expec[\SS_t]/(tu^{\tau-3}) \in [\underline{c},\overline{c}]$ for $t\in [\vep,u/2]$ and some constants $\underline{c},\overline{c}>0$. Therefore, taking $x=a \eta t^{\frac{1}{2} (5-\tau)}u^{\tau-3}$
for some $a>0$ chosen appropriately,
	\eqn{
	\lbeq{eq:integer1}
	\widetilde \prob\left( |\SS_t-   \widetilde\expec[\SS_t]|> \eta \widetilde \expec[\SS_{t}]  \right)
	\le c\e^{-a \eta t^{\frac{1}{2}(5-\tau)}u^{\tau-3}}.
	}
We take $t=k\delta_u$ for $k\delta_u \in [\vep,u/2]$, so that there are at most $u/\delta_u=u^{\tau-1}$possible values of $k$. Thus,
 	\eqn{
	\lbeq{eq:integer3a}
	\widetilde \prob\left( \exists t_k\in [\vep,u/2]\colon |\SS_{t_k}-\widetilde\expec[\SS_{t_k}]|
	>\eta \widetilde \expec[\SS_{t_k}]  \right)
	\le c(\epsilon)u^{\tau-1}\e^{-a \eta u^{(\tau-1)/2}}.
	}
This proves the proposition for $k\delta_u \in [\vep,u/2]$.

\paragraph{{\bf Proof for $t\in [u/2,u-\vep]$.}} We continue by proving the proposition for $t\in [u/2,u-\vep]$, for which we again use Proposition \ref{prop-laplace} with $\lambda_1=\pm1$ and $\lambda_2=0$ to see that, for any $x>0$,
	\eqn{
	\widetilde \prob\Big( \big|{\SS_{t}- \widetilde \expec[\SS_{t}]
	\over \sqrt{I_{\sss V}(t/u) u^{\tau-3}}}\big| >x\Big)
	\le \centering c\e^{-x}.
	}
By Lemma \ref{lem-mean-var-extreme} and the fact that $I_{\sss V}(p)>0$ for every $p\in(0,1)$, we obtain that there exists a constant $c>0$ such that $I_{\sss V}(p) \ge c$ for all $p\in [1/2, 1-\vep]$.
%\todotemp{Remco: Check references!}
Applying this to $p=t/u$ yields
	\eqn{
	\widetilde \prob\left( |\SS_t -\widetilde \expec[\SS_{t}] | > c x u^{(\tau-3)/2}\right) \le c\e^{-x}.
	}
By Lemma \ref{lem-expec-S}(d) and Corollary \ref{cor-asympt-mean-S}(b), we have $\widetilde \expec[\SS_t] u^{-(\tau-2)} \in [\underline{c},\overline{c}]$ for all $t\in [u/2,u-\vep]$ and some constants $\underline{c}=\underline{c}(\vep),\overline{c}=\overline{c}(\vep)>0$. Therefore, taking $x=a \eta u^{(\tau-1)/2}$
for some $a>0$ chosen appropriately,
	\eqn{
	\lbeq{eq:integer2b}
	\widetilde \prob\left( |\SS_t- \widetilde\expec[\SS_t]|> \eta \widetilde \expec[\SS_{t}]  \right)
	\le c\e^{-a \eta u^{(\tau-1)/2}}.
	}
We take $t=t_k=k\delta_u$ for $k\delta_u \in [u/2, u-\vep]$, so that there are at most $u/\delta_u=u^{\tau-1}$
possible values of $k$. Thus,
	\eqn{
	\lbeq{eq:integer3b}
	\widetilde \prob\left( \exists t_k\in [u/2, u-\vep]\colon |\SS_{t_k}-   \widetilde\expec[\SS_{t_k}]|> \eta \widetilde \expec[\SS_{t_k}]  \right)
	\le c(\epsilon)u^{\tau-1}\e^{-a \eta u^{(\tau-1)/2}}.
	}
This proves the proposition for $k\delta_u \in [u/2,u-\vep]$.

\paragraph{{\bf Proof for $t\in[u-\vep,u-Tu^{-(\tau-2)}]$: Rewrite.}}
The proof for $t\in[u-\vep,u-Tu^{-(\tau-2)}]$ is the hardest, and is split into three steps. We start by rewriting the event of interest. We define $s=u-t$ and investigate $\SS_{u-s}$ in what follows, so that now $s\in [Tu^{-(\tau-2)},\vep]$.

Recall the definition of $Q_u(s)$ in \refeq{Qu-def},
    \eqan{
    \lbeq{Qu-def-repeat}
    Q_u(s)
    &= u \SS_{u-s} - s - (u-s) \SS_u
    =-\sum_{i=2}^{\infty} c_i[u\indic{T_i\in (u-s,u]}-s\II_i(u)],
    }
so that $|\SS_{u-s}-\widetilde \expec[\SS_{u-s}]|>\eta \widetilde \expec[\SS_{u-s}]$ precisely when
	\eqn{
	\lbeq{rewrite-step-1}
	|Q_u(s) -\widetilde \expec[Q_u(s)]+(u-s) \big(\SS_u-\widetilde \expec[\SS_u]\big)|
	> \eta u \widetilde \expec[\SS_{u-s}].
	}
When $\SS_u\in [0,M/u]$ and using that $u\widetilde \expec[\SS_u]=o(1)$ by Lemma \ref{lem-expec-S}(d), we therefore obtain that if \refeq{rewrite-step-1} holds, then
	\eqn{
	|Q_u(s) -\widetilde \expec[Q_u(s)]|
	> \eta u \widetilde \expec[\SS_{u-s}]-M+o(1).
	}
By Lemma \ref{lem-expec-S}(d) and Corollary \ref{cor-asympt-mean-S}(b), we have that $\widetilde \expec[\SS_{u-s}]\geq c su^{\tau-3}$ for some $c>0$. Therefore, $\eta u \widetilde \expec[\SS_{u-s}]\geq c \eta T$, so that, by taking $T=T(M)$ sufficiently large, we obtain that%, for any $s\in [Tu^{-(\tau-2)},\vep]$,
	\eqan{
	\lbeq{rewrite-step-2}
	&\widetilde \prob\Big(\exists s_k\in [Tu^{-(\tau-2)},\vep]\colon\big|\SS_{u-s_k} - \widetilde \expec [\SS_{u-_k}]\big| >\eta
	\widetilde \expec [\SS_{u-s_k}], \SS_u\in[0,M/u]\Big)\\
	&\qquad \leq \widetilde \prob\Big(\exists s_k\in [Tu^{-(\tau-2)},\vep]\colon|Q_u(s_k) -\widetilde \expec[Q_u(s_k)]|
	> \eta c s_ku^{\tau-2}, \SS_u\in[0,M/u]\Big).\nn
	}
We condition on $\JJ(u)$ from (\ref{index-def}), and note that $\SS_u$ is measurable w.r.t\ $\JJ(u)$ to obtain
	\eqan{
	&\widetilde \prob\Big(\exists s_k\in [Tu^{-(\tau-2)},\vep]\colon|Q_u(s_k) -\widetilde \expec[Q_u(s_k)]|
	> \eta c s_ku^{\tau-2}, \SS_u\in[0,M/u]\Big)\\
	&\qquad =\widetilde \expec\Big[\indic{\SS_u\in[0,M/u]} \widetilde \prob\big(\exists s_k\in [Tu^{-(\tau-2)},\vep]\colon|Q_u(s_k) -\widetilde \expec[Q_u(s_k)]|
	> \eta c s_ku^{\tau-2}\mid \JJ(u)\big)\Big].\nn
	}
This is the starting point of our analysis.
%We use Markov's inequality to bound
%	\eqan{
%	&\widetilde \prob\Big(|Q_u(t) -\widetilde \expec[Q_u(t)]|
%	> \eta c tu^{\tau-2}/2, \SS_u\in[0,M/u]\Big)\\
%	&\qquad \leq (\eta c tu^{\tau-2})^{-4}
%	\widetilde \expec\Big[\indic{\SS_u\in[0,M/u]}
%	\big(Q_u(t) -\widetilde \expec[Q_u(t)]\big)^4\Big].\nn
%	}
%The idea is that $\widetilde \expec\Big[\big(Q_u(t) -\widetilde \expec[Q_u(t)]\big)^4\Big]=O((tu^{\tau-2})^2)$. This suggests that
%	\eqn{
%	\widetilde \expec\Big[\indic{\SS_u\in[0,M/u]}
%	\big(Q_u(t) -\widetilde \expec[Q_u(t)]\big)^4\Big]=O((tu^{\tau-2})^2)	
%	\widetilde \prob\big(\SS_u\in[0,M/u]\big).
%	}
%In total, this would give
%	\eqan{
%	&\widetilde \prob\Big(|Q_u(t) -\widetilde \expec[Q_u(t)]|
%	> \eta c tu^{\tau-2}/2, \SS_u\in[0,M/u]\Big)=O(\eta^{-4} (tu^{\tau-2})^{-2}).\nn
%	}
%This bound is true for any $t\in[Tu^{-(\tau-2)},u/2]$. Taking $t=t_k=k u^{-(\tau-2)}$ and summing out over $k\geq T$ leads to
%	\eqn{
%	u^{(\tau-1)/2}\widetilde \prob\Big(\exists k\geq T\colon |Q_u(t) -\widetilde \expec[Q_u(t)]|
%	> \eta c tu^{\tau-2}/2, \SS_u\in[0,M/u]\Big)
%	\leq O(\eta^{-4}) \sum_{k\geq T}  k^{-2}=O(\eta^{-4}/T)=o_{\sss T}(1),
%	}
%when we take $T=T(\eta)$ sufficiently large, as required.
%\todotemp{Remco: Add in details!}
%
%\noindent
%{\bf Old proof:}
We split, writing $\eta'=\eta/2$,
	\eqan{
	&\widetilde \prob\big(\exists s_k\in [Tu^{-(\tau-2)},\vep]\colon|Q_u(s_k) -\widetilde \expec[Q_u(s_k)]|
	> \eta c s_ku^{\tau-2}\mid \JJ(u)\big)\\
	&\qquad \leq \widetilde \prob\big(\exists s_k\in [Tu^{-(\tau-2)},\vep]\colon|Q_u(s_k) -\widetilde \expec[Q_u(s_k)\mid \JJ(u)]|
	> \eta' c s_ku^{\tau-2}\mid \JJ(u)\big)\nn\\
	&\qquad \qquad+
	\indic{\exists s_k\in [Tu^{-(\tau-2)},\vep]\colon |\widetilde \expec[Q_u(s_k)\mid \JJ(u)]-\widetilde \expec[Q_u(s_k)]|>\eta' c s_ku^{\tau-2}}.\nn
	}
%\todo[inline, color=magenta]{Remco: Update from here!}
We conclude using the union bound that
%\todo[inline]{Sandra: need to note that we can add $\exists t_k$ ... in \refeq{rewrite-step-2} on both sides or better directly work with $t_k$'s from \refeq{rewrite-step-2} onwards. Remco:  Discuss 19.}
	\eqan{
	\lbeq{cond-expec-split}
	&\widetilde \prob\Big(\exists s_k\in [Tu^{-(\tau-2)},\vep]\colon|Q_u(s_k) -\widetilde \expec[Q_u(s_k)]|
	> \eta c s_ku^{\tau-2}, \SS_u\in[0,M/u]\Big)\\
	&\qquad \leq \sum_{k\colon s_k\in  [Tu^{-(\tau-2)},\vep]}
	\widetilde \expec\Big[\widetilde \prob\Big(|Q_u(s_k) -\widetilde \expec[Q_u(s_k)\mid \JJ(u)]|
	> \eta' c s_ku^{\tau-2}\mid \JJ(u)\Big)\indic{\SS_u\in [0,M/u]}\Big]\nn\\
	&\qquad\qquad +
	\widetilde \prob\Big(\exists s_k\in [Tu^{-(\tau-2)},\vep]\colon |\widetilde \expec[Q_u(s_k)\mid \JJ(u)]-\widetilde \expec[Q_u(s_k)]|>\eta' c s_k u^{\tau-2}\Big).\nn
	}
%\todo[inline]{Remco: Replace first term using union bound and sum over $k$. Discuss 19 (once again)}
We will bound both contributions separately, and start by setting the stage. We compute that
	\eqan{
	\lbeq{diff-q-cond-exp}
	Q_u(s) -\widetilde \expec[Q_u(s)\mid \JJ(u)]
	&= -\sum_{i=2}^{\infty} c_i u\big[\indic{T_i\in (u-s,u]}-\widetilde \prob(T_i\in (u-s,u]\mid \JJ(u))]\\
	&=-\sum_{i=2}^{\infty} c_i u\big[\indic{T_i\in (u-s,u]}-p_{i,u}(s)\big],\nn
	}
where we abbreviate
	\eqn{
	\lbeq{p-iu-t}
	p_{i,u}(s)=\widetilde \prob(T_i\in (u-s,u]\mid \JJ(u))=\widetilde \prob(T_i\in (u-s,u]\mid i\in\JJ(u)).
	}
It turns out that both contributions in \refeq{cond-expec-split} can be expressed in terms of
$p_{i,u}(s)$, and we continue our analysis by studying this quantity in more detail.

\paragraph{{\bf Proof for $t\in[u-\vep,u-Tu^{-(\tau-2)}]$: Analysis of $p_{i,u}(s)$.}} We next analyse the
conditional probability $p_{i,u}(s)$. We compute (recall \refeq{def-ci}, \refeq{II-proc-tilt} and \refeq{index-def})
	\eqn{
	p_{i,u}(s)=\widetilde \prob(T_i\in (u-s,u]\mid i\in\JJ(u))
	=\frac{\widetilde \prob(T_i\in (u-s,u])}{\widetilde \prob(T_i\leq u)}
	=\frac{\widetilde \prob(T_i\leq u)-\widetilde \prob(T_i\leq u-s)}{\widetilde \prob(T_i\leq u)}.
	}
Using the distribution of $T_i$ formulated in Lemma
\ref{lem-ind-tilt}, we obtain, for any $s\in [0,u]$,
	\eqn{
	\lbeq{IIit-probs}
	\widetilde \prob(T_i\leq u-s) =
    	\frac{\e^{\theta c_iu}(1-\e^{-c_i (u-s)})}{\e^{\theta c_iu}(1-\e^{-c_iu})
    	+\e^{-c_iu}},
	}
so that
	\eqan{
	\lbeq{piut-eq}
	p_{i,u}(s)=
	\frac{\e^{\theta c_iu}(1-\e^{-c_i u})-\e^{\theta c_iu}(1-\e^{-c_i (u-s)})}
	{\e^{\theta c_iu}(1-\e^{-c_i u})}
	=\frac{\e^{-c_i u}(\e^{c_i s}-1)}
	{1-\e^{-c_i u}}=\frac{\e^{c_i s}-1}
	{\e^{c_i u}-1}.
	}
We start by bounding $p_{i,u}(s)$, for $s\in[0,\vep]$, by
	\eqn{
	\lbeq{piut-bd}
	p_{i,u}(s)\leq
	O(s/u),
	\qquad
	\text{and}
	\qquad
	p_{i,u}(s)\leq
	O(c_i s) \e^{-c_iu} (c_i u\wedge 1)^{-1}.
	}
Moreover, %for $s\in[0,\vep]$,
%	\eqn{
%	\lbeq{cond-prob-Is}
%	up_{i,u}(s)-s=s\frac{c_i u-\e^{c_i u}+1}
%	{\e^{c_i u}-1} +u\frac{\e^{c_i s}-1-c_is}
%	{\e^{c_i u}-1}=s\frac{c_i u-\e^{c_i u}+1}
%	{\e^{c_i u}-1} + O(us^2)\frac{c_i^2}
%	{\e^{c_i u}-1},
%	}
%\todo{S-after meeting: I checked, we do not use \eqref{cond-prob-Is} and \eqref{piut-t-bd} ...}
%\todo[inline]{S-add I think this yields $+ u O(s^2)\frac{c_i^2}{\e^{c_i u}-1}$ - note that we do not use it in what follows. Remco: You are right, I had a factor $u$ too many... This is corrected.}
%where the $O(s^2)$ term is non-negative. Using that $\e^{c_i u}\geq 1+c_iu$ and $\e^{c_i u}\geq 1+c_iu +(c_iu)^2/2$, we conclude that,
for $u$ sufficiently large,
	\eqn{
	\lbeq{piut-t-bd}
	|up_{i,u}(s)-s|\leq s(c_iu\wedge 1).
	}

\paragraph{{\bf Proof for $t\in[u-\vep,u-Tu^{-(\tau-2)}]$: Completion first term \refeq{cond-expec-split}.}}
For the first term in \refeq{cond-expec-split}, we use Markov's inequality in the form $\prob(|X-\expec[X]|>a)\leq a^{-4} \expec[(X-\expec[X])^4]$ to obtain
	\eqn{
	\lbeq{diff-q-cond-exp-repeat}
	\widetilde \prob\big(|Q_u(s) -\widetilde \expec[Q_u(s)\mid \JJ(u)]|
	> \eta' c su^{\tau-2}\mid \JJ(u)\big)
	\leq (\eta' c su^{\tau-2})^{-4}
	\widetilde \expec\big[(Q_u(s) -\widetilde \expec[Q_u(s)\mid \JJ(u)])^4\mid \JJ(u)\big],
	}
and recall from \refeq{diff-q-cond-exp} that
	\eqn{
	Q_u(s) -\widetilde \expec[Q_u(s)\mid \JJ(u)]
	= -\sum_{i=2}^{\infty} c_i u\big[\indic{T_i\in (u-s,u]}-\widetilde \prob(T_i\in (u-s,u]\mid \JJ(u))] = -\sum_{i=2}^{\infty} c_i u\big[\indic{T_i\in (u-s,u]}-p_{i,u}(s)\big].
	}
The summands are conditionally independent given $\JJ(u)$ and identically 0 when $\II_i(u)=0$, so that
	\eqan{
	\widetilde \expec\big[(Q_u(s) -\widetilde \expec[Q_u(s)\mid \JJ(u)])^4\mid \JJ(u)\big]
	&\leq \sum_{i\geq 2} c_i^4 u^4 p_{i,u}(s)\II_i(u)\\
	&\qquad +\sum_{i,j\geq 2\colon i\neq j} c_i^2 c_j^2 u^4
	p_{i,u}(s)(1-p_{i,u}(s))\II_i(u) p_{j,u}(s)(1-p_{j,u}(s))\II_j(u).\nn
	}
By the second bound in \refeq{piut-bd} and Corollary \ref{cor-exp-mom-cs},
%\todo[inline]{S-add proof not existent in other file (but in tex-version in gray). Remco: Discuss 20, also applies to LD paper.}
the first term is at most
	\eqn{
	O(1) s u^4 \sum_{i\geq 2} c_i^5 \e^{-c_iu} (c_i u\wedge 1)^{-1}
	\leq O(1) s u^4 \sum_{i\geq 2} c_i^5 \e^{-c_iu} [1+(c_i u)^{-1}]=O(su^{\tau-2}).
	}
By \refeq{first-bd-sec-mom} in Lemma \ref{lem-LD-par}, we may assume
%\todo[inline]{S-add I would drop this formulation and add $+o(1)$ in the two following equations; afterwards not necessary anymore; then we can also drop the comment {\it on the event ...} following the third equation from here onwards. Remco: Note sure, since we are both working with conditional as well as with regular expectations...}
that $\sum_{i=2}^{\infty} c_i^2 \II_i(u)\leq K u^{\tau-3}$, since the complement has a probability that is $o(u^{-(\tau-1)/2})$. Then, in a similar way, using the first bound in \refeq{piut-bd}, the second term is at most
	\eqan{
	\Big(\sum_{i\geq 2} c_i^2 u^2
	p_{i,u}(s)\II_i(u)\Big)^2
	&\leq O(1)\Big(s\sum_{i\geq 2} c_i^2 u\II_i(u)\Big)^2=O\big((su^{\tau-2})^2\big).
	}
As a result,
	\eqn{
	\widetilde \expec\big[(Q_u(s) -\widetilde \expec[Q_u(s)\mid \JJ(u)])^4\mid \JJ(u)\big]
	\leq O(su^{\tau-2})+ O\big((su^{\tau-2})^2\big).
	}
Since $s\geq Tu^{-(\tau-2)}$, this can be simplified to
	\eqn{
	\widetilde \expec\big[(Q_u(s) -\widetilde \expec[Q_u(s)\mid \JJ(u)])^4\mid \JJ(u)\big]
	\leq O\big((su^{\tau-2})^2\big).
	}
We conclude using \eqref{diff-q-cond-exp-repeat} that, on the event that $\{\sum_{i=2}^{\infty} c_i^2 \II_i(u)\leq K u^{\tau-3}\},$
	\eqn{
	\widetilde \prob\big(|Q_u(s) -\widetilde \expec[Q_u(s)\mid \JJ(u)]|
	> \eta' c su^{\tau-2}\mid \JJ(u)\big)\leq  \frac{(cs u^{\tau-2})^2}{(\eta' c su^{\tau-2})^{4}}
	=O(\eta^{-4} (su^{\tau-2})^{-2}),
	}
so that, also using that $\widetilde \prob(\SS_u\in[0,M/u]) = O(u^{-(\tau-1)/2})$ by Proposition \ref{prop-dens-f},
%\todo[inline]{S-? I still think we should replace the following two probabilities with the first term (on the r.h.s.) from \refeq{cond-expec-split}. After all, we are in the part called {\it Completion {\bf first} term \refeq{cond-expec-split}.}}
	\eqan{
	&u^{(\tau-1)/2}
	\widetilde \expec\Big[\widetilde \prob\Big(|Q_u(s_k) -\widetilde \expec[Q_u(s_k)\mid \JJ(u)]|
	> \eta' c s_ku^{\tau-2}\mid \JJ(u)\Big)\indic{\SS_u\in [0,M/u]}\Big]\nn\\
	&\qquad\leq O(\eta^{-4} (su^{\tau-2})^{-2})u^{(\tau-1)/2}\widetilde \prob(\SS_u\in[0,M/u])=O(\eta^{-4} (su^{\tau-2})^{-2}).
	}
This bound is true for any $s\in[Tu^{-(\tau-2)},\vep]$. Taking $s=s_k=k u^{-(\tau-2)}$ and summing out over $k\geq T$ leads to
	\eqan{
	&u^{(\tau-1)/2}\sum_{k\colon s_k\in  [Tu^{-(\tau-2)},\vep]}\widetilde \expec\Big[\widetilde \prob\Big(|Q_u(s_k) -\widetilde \expec[Q_u(s_k)\mid \JJ(u)]|
	> \eta' c s_ku^{\tau-2}\mid \JJ(u)\Big)\indic{\SS_u\in [0,M/u]}\Big]\\
	&\qquad \leq O(\eta^{-4}) \sum_{k\geq T}  k^{-2}=O(\eta^{-4}/T)=o_{\sss T}(1),\nn
	}
when we take $T=T(\eta)$ sufficiently large, as required.
%\todo{Remco: Update from here!}

\paragraph{{\bf Proof for $t\in[u-\vep,u-Tu^{-(\tau-2)}]$: Completion second term \refeq{cond-expec-split}.}}
For the second term in \refeq{cond-expec-split}, we need to bound
	\eqn{
	\widetilde \prob\Big(\exists s_k\in [Tu^{-(\tau-2)}, \vep]\colon|\widetilde\expec[Q_u(s_k)\mid \JJ(u)]-\widetilde\expec[Q_u(s_k)]|>\eta' c s_ku^{\tau-2}\Big).
	}
We compute using \refeq{Qu-def-repeat}
	\eqn{
	\widetilde\expec[Q_u(s)]
	=-\sum_{i=2}^{\infty} c_i\big[u\widetilde \prob(T_i\in (u-s,u])-s\widetilde \prob(i\in \JJ(u))\big],
	}
while
	\eqn{
	\widetilde\expec[Q_u(s)\mid \JJ(u)]
	%=\sum_{i=2}^{\infty} c_i[u\widetilde \prob(T_i\in (u-t,u]\mid \JJ(u))-t\II_i(u)]
	=-\sum_{i=2}^{\infty} c_i \II_i(u) \big[u\widetilde \prob(T_i\in (u-s,u]\mid i\in\JJ(u))-s\big].
	}
As a result, using \refeq{p-iu-t},
	\eqan{
	\widetilde\expec[Q_u(s)\mid \JJ(u)]-\widetilde\expec[Q_u(s)]
	&=-\sum_{i=2}^{\infty} c_i \big[\II_i(u)-\widetilde \prob(i\in \JJ(u))\big]
	[up_{i,u}(s)-s]=:
	s X+Y(s),
	}
where with \refeq{piut-eq}
	\eqan{
	X=-\sum_{i=2}^{\infty} c_i \big[\II_i(u)-\widetilde \prob(i\in \JJ(u))\big]\frac{1+c_iu-\e^{c_iu}}{\e^{c_iu}-1},
	\qquad
	Y(s)=-u\sum_{i=2}^{\infty} c_i \big[\II_i(u)-\widetilde \prob(i\in \JJ(u))\big]
	\frac{\e^{c_is}-1-c_is}{\e^{c_iu}-1}.
	}
As a result,
	\eqan{
	\lbeq{sum-X-Y-parts}
	&\widetilde \prob\Big(\exists s_k\in [Tu^{-(\tau-2)},\vep]\colon |\widetilde \expec[Q_u(s_k)\mid \JJ(u)]-\widetilde \expec[Q_u(s_k)]|>\eta' c s_ku^{\tau-2}\Big)\\
	&\qquad\leq \widetilde \prob\big(|X|\geq \eta' c u^{\tau-2}/2\big)+
	\widetilde \prob\Big(\exists s_k\in [Tu^{-(\tau-2)},\vep]\colon |Y(s_k)|\geq \eta' c s_ku^{\tau-2}/2\Big).
	\nn
	}
For both terms, we use the Chebychev inequality.

For $X$, as $\widetilde \expec[X]=0$, this leads to
	\eqn{
	\widetilde \prob\big(|X|\geq \eta' c u^{\tau-2}/2\big)
	\leq \frac{4}{(\eta' c u^{\tau-2})^2}\Var(X).
	}
We use Lemma \ref{lem-ind-tilt} to see that
$\widetilde \prob(i\in \JJ(u))=\frac{1-\e^{-c_i u}}{1-\e^{-c_iu}
   +\e^{-c_iu(1+\theta)}}$, so that
	\eqn{
	\lbeq{Iiu-not-Iiu-bd}
	\widetilde \prob(i\in \JJ(u))\widetilde \prob(i\not\in \JJ(u))
	=\widetilde \prob(\II_i(u)=1)\widetilde \prob(\II_i(u)=0)
	=\frac{(1-\e^{-c_i u})\e^{-c_iu(1+\theta)}}
	{(1-\e^{-c_iu}+\e^{-c_iu(1+\theta)})^2}\leq O(1)c_iu\e^{-c_iu(1+\theta)},
	}
since $1-\e^{-x}+\e^{-x(1+\theta)}$ is uniformly bounded from below away from 0 for all $x\geq 0$.
We use this together with Corollary \ref{cor-exp-mom-cs} to compute that
	\eqan{
	\Var(X)&=\sum_{i=2}^{\infty} c_i^2 \widetilde \prob(\II_i(u)=1)
	\widetilde \prob(\II_i(u)=0)\Big(\frac{1+c_iu-\e^{c_iu}}{\e^{c_iu}-1}\Big)^2\\
	&\leq O(1) u\sum_{i=2}^{\infty} c_i^3 \e^{-c_iu(1+\theta)}=O(1) u^{\tau-3}.\nn
	}
Therefore,
	\eqn{
	\lbeq{first-bound-X-Y-part}
	\widetilde \prob\big(|X|\geq \eta' c u^{\tau-2}/2\big)
	\leq O(1) (\eta')^{-2} u^{-2(\tau-2)}\Var(X)=O(1) (\eta')^{-2} u^{-(\tau-1)}=o(u^{-(\tau-1)/2}),
	}
as required below.

For the term involving $Y(s)$, we start by using the union bound to obtain
	\eqan{
	\lbeq{union-Y(t)}
	\widetilde \prob\Big(\exists s_k\in [Tu^{-(\tau-2)}, \vep]\colon |Y(s_k)|\geq \eta' c s_ku^{\tau-2}/2\Big)
	\leq \vep u^{\tau-2} \max_{k\colon s_k\in [Tu^{-(\tau-2)}, \vep]}
	\widetilde \prob\big(|Y(s_k)|\geq \eta' c s_ku^{\tau-2}/2\big).
	}
Then, by the Chebychev inequality and as $\widetilde \expec[Y(s_k)]=0$,
	\eqan{
	\lbeq{union-Y(t)b}
	\widetilde \prob\big(|Y(s_k)|\geq \eta' c s_k u^{\tau-2}/2\big)
	&\leq \frac{4}{(\eta' c s_k u^{\tau-2})^2} \Var(Y(s_k)),
	}
where, using \refeq{Iiu-not-Iiu-bd},  $\e^{c_is}-1-c_is=O(s^2c_i^2)$ and $\e^{c_iu}-1\geq c_iu$,
	\eqan{
	\Var(Y(s)) &=u^2\sum_{i=2}^{\infty} c_i^2 \widetilde \prob(i\in \JJ(u))\widetilde \prob(i\not\in \JJ(u))
	\Big(\frac{\e^{c_is}-1-c_is}{\e^{c_iu}-1}\Big)^2\\
	&\leq O(s^4) \sum_{i=2}^{\infty} c_i^4\e^{-c_iu(1+\theta)}=O(s^4) O(u^{\tau-5}),\nn
	}
where we used Corollary \ref{cor-exp-mom-cs} in the last equality.
Substituting this into \refeq{union-Y(t)} and \refeq{union-Y(t)b}, we arrive at
	\eqan{
	\lbeq{second-bound-X-Y-part}
	\widetilde \prob\Big(\exists s_k\in [Tu^{-(\tau-2)}, \vep]\colon |Y(s_k)|\geq \eta' c s_ku^{\tau-2}/2\Big)
	&\leq \vep u^{\tau-2}
	\max_{k\colon s_k\in [Tu^{-(\tau-2)}, \vep]} O(s_k^2 u^{\tau-5} u^{-2(\tau-2)})(\eta')^{-2}\nn\\
	&=O(\vep^3/(\eta')^2) u^{-3}=o(u^{-(\tau-1)/2}),
	}
since $\tau\in(3,4)$. Combining \refeq{first-bound-X-Y-part} and \refeq{second-bound-X-Y-part} in \refeq{sum-X-Y-parts} completes the proof.
\qed
\bigskip

We now know that with high probability the process does not deviate much from its mean
when observed at the discrete times $k\delta_u \in[\vep,u-T\delta_u]$.
We continue to show that this actually holds with high probability
on the whole interval $[\vep,u-T\delta_u]$.
We complete the preparations for the proof of Proposition \ref{prop-no-middle-ground}
by proving that it is unlikely for the process to deviate far from the mean for all times $t\in[\vep, u-T\delta_u]$ simultaneously:

\begin{Proposition}[Probability to deviate far from mean at some time]
%\todo{S-? at some time?}
\label{prop-no-hit-small-interval}
For every $\eta>0$ and $M>0$,
	\eqn{
	\lbeq{prob-dev-far-all-times}
	\limsup_{u\to \infty} u^{(\tau-1)/2} \widetilde \prob\big(\exists \, t\in[\vep,u-T\delta_u]\colon
	|\SS_t-\widetilde\expec[\SS_t]|\geq 2\eta \widetilde\expec[\SS_t], \SS_u\in[0,M/u]\big)=o_{\sss T}(1).
	}
\end{Proposition}

%\todo[inline, color=magenta]{Remco: Update from here!}

\proof
Fix $T>0$ and recall that $\delta_u=u^{-(\tau-2)}$. Let
	\eqan{
	\lbeq{def-Eu-prob-dev-far-all-times}
	E_u&=\{|\SS_{k\delta_u}-\widetilde\expec[\SS_{k\delta_u}]|\leq \eta \widetilde\expec[\SS_{k\delta_u}]~\forall k\text{ s.t. }k\delta_u\in [\vep,u-T\delta_u]\}\\
	&\qquad
	\cap \Big\{\sum_{i=2}^{\infty} c_i [1-\II_i(u/2)]\big(\e^{\lambda c_i}-1-\lambda c_i\big)
	\leq K\lambda^2u^{\tau-4}\Big\}\cap \Big\{\sum_{i=2}^{\infty} c_i^2 \II_i(u)\leq K u^{\tau-3}\Big\},\nn
	}
where we take $\lambda=\delta u$ with $\delta>0$ sufficiently small and $K \geq 1$ as in Lemma \ref{lem-LD-par}. %\todo{check}
We first give a bound on
$\widetilde \prob(E_u^c \cap \{\SS_u\in[0,M/u]\})$. We apply \refeq{first-bd-sec-mom} in Lemma \ref{lem-LD-par}
to obtain that
	\eqn{
	\widetilde \prob(\sum_{i=2}^{\infty} c_i^2 \II_i(u)
	\geq K u^{\tau-3})=O(u^{-(\tau-1)})=o(u^{-(\tau-1)/2}),
	}
which is contained in the error term in \refeq{prob-dev-far-all-times}. Further,
by \refeq{sec-bd-sec-mom} in Lemma \ref{lem-LD-par}
	\eqn{
	\widetilde \prob(\sum_{i=2}^{\infty} c_i [1-\II_i(u/2)]\big(\e^{\lambda c_i}-1-\lambda c_i\big)
	\geq K\lambda^2u^{\tau-4})=O(u^{-(\tau-1)})=o(u^{-(\tau-1)/2}).
	}
Combined with Proposition \ref{prop-integer-times}, this ensures that
	\eqn{
	\limsup_{u\to\infty} u^{(\tau-1)/2} \widetilde \prob(E_u^c\cap \{\SS_u\in[0,M/u]\} ) =o_{\sss T}(1).
	}
As a result, we are left to control the fluctuations of the process on any interval $I_k=[k\delta_u,(k+1)\delta_u]$. We use Boole's inequality to bound
	\eqan{
	&\widetilde \prob(E_u, \exists \, t\in[\vep,u-T\delta_u]\colon
	|\SS_t-\widetilde\expec[\SS_t]|\geq 2\eta \widetilde\expec[\SS_t])\\
	&\quad \leq \sum_{k\colon k\delta_u\in [\vep,u-T\delta_u]}\!\!\!
	\widetilde \prob\big(E_u, \exists \, t\in I_k\colon
	|\SS_t-\widetilde\expec[\SS_t]|\geq 2\eta \widetilde\expec[\SS_t]\big).\nn
	}
Let $t_k=k\delta_u$, so that $I_k=[t_{k},t_{k+1}]$.
We split the analysis into four cases, depending on whether $t_k\leq u/2$
or not, and on whether $\SS_t-\widetilde\expec[\SS_t]\geq 2\eta \widetilde\expec[\SS_t]$
or $\SS_t-\widetilde\expec[\SS_t]\leq -2\eta \widetilde\expec[\SS_t]$, which we refer to
as `large upper' and `large lower' deviations, respectively.

\paragraph{Part 1: The case $t_k\leq u/2$ and a large upper deviation.}
We start by bounding the probability that there exists a $t\in I_k=[t_{k},t_{k+1}]$, $\epsilon \le t_k \le u/2$ such that
$\SS_t-\widetilde\expec[\SS_t] \geq 2\eta \widetilde\expec[\SS_t]$.
Using that $\widetilde\expec[\SS_{t}]=\widetilde\expec[\SS_{t_k}](1+o(1))$ for any $t\in I_k$ by Corollary \ref{cor-asympt-mean-S}(c), we bound
	\eqn{
	\widetilde \prob(E_u, \exists \, t\in I_k\colon
	\SS_t-\widetilde\expec[\SS_t]\geq 2\eta \widetilde\expec[\SS_t])
	\leq \widetilde \prob(\exists \, t\leq \delta_u\colon
	\SS_{t+t_k}-\SS_{t_k}\geq \frac{\eta}{2}\widetilde\expec[\SS_{t_k}]).
	}
By \eqref{SS-def-ref},
	\eqn{
	\SS_{t+t_k}-\SS_{t_k}=
	\betam t+\sum_{i=2}^{\infty} c_i [\II_i(t+t_k)-\II_i(t_k)-c_i t],
	}
which can be stochastically dominated by the process $\betam t + \RR_t$ with
$\RR_t \equiv \sum_{i=2}^{\infty} c_i [N_i(t)-c_i t],$
where $(N_i(t))_{t\geq 0}$ is a Poisson process with rate $c_it$. As a result,
	\eqn{
	\widetilde \prob(E_u, \exists \, t\in I_k\colon
	\SS_t-\widetilde\expec[\SS_t]\geq 2\eta \widetilde\expec[\SS_t])
	\leq \widetilde \prob(\exists \, t\leq \delta_u\colon
	\betam t + \RR_t\geq \frac{\eta}{2}\widetilde\expec[\SS_{t_k}]).
	}
%\todotemp{I changed what followed about L\'evy processes. I could not find a proper reference and indeed would guess, that the claimed formulation does not hold for $\RR_t$ but rather for $-\RR_t$ as we only have positive jumps?!}
%Since $(-\RR_t)_{t\geq 0}$ is a L\'evy process with no positive jumps, we obtain from \cite[Section VII.1]{Bert96} that $\widetile\expec[ \e^{-\lambda \RR_t} ] = \e^{t \phi(\lambda)}$ for all $\lambda \geq 0$ and where $\phi(\lambda) = -\Psi(-i\lambda) = \log(\widetilde\expec[ \e^{-\lambda \RR_1}])$ (cf. \cite[Section I.1]{Bert96}). For $\lambda\ge 0$, we now define the exponential martingale (use that a L\'evy-process $\RR_t$ has independent increments and that the law of $\RR_t-\RR_s$ is the same as that of $\RR_s$)
%	\eqn{
%	\MM_t(\lambda)=\e^{-\lambda\RR_t-t\phi(\lambda)},
%	\qquad
%	\text{where}
%	\qquad
%	\phi(\lambda)=\log{\widetilde\expec[\e^{-\lambda\RR_1}]}
%	=\sum_{i=2}^{\infty} c_i[\e^{-c_i\lambda}-1+c_i\lambda].
%	}
Since $(\RR_t)_{t\geq 0}$ is a finite-variance L\'evy process, it is well-concentrated. In more detail,
for $\lambda\in \R$, we define the exponential martingale
	\eqn{
	\lbeq{exp-mart}
	\MM_t(\lambda)=\e^{\lambda\RR_t-t\phi(\lambda)},
	\qquad
	\text{where}
	\qquad
	\phi(\lambda)=\log{\widetilde\expec[\e^{\lambda\RR_1}]}
	=\sum_{i=2}^{\infty} c_i[\e^{c_i\lambda}-1-c_i\lambda].
	}
Then, for every $\lambda\geq 0$, using that $\phi(\lambda)\geq 0$ and by
Doob's inequality,
	\eqan{
	\widetilde \prob(\exists \, s\leq t\colon
	\betam s + \RR_s\geq x)
	&\leq \widetilde \prob(\exists \, s\leq t\colon
	\MM_s(\lambda)\geq \e^{x\lambda-t\phi(\lambda)-t |\betam| \lambda}) \\
  &\leq \e^{-2[x\lambda-t\phi(\lambda)-t |\betam| \lambda]}\widetilde\expec[\MM_t(\lambda)^2]
	=\e^{-[2x\lambda-t\phi(2\lambda)-2t |\betam| \lambda}].\nn
	}
We apply this inequality to $x=\frac{\eta}{2}\widetilde\expec[\SS_{t_k}]$, $t=\delta_u$
and $\lambda=1$, and Corollary \ref{cor-asympt-mean-S}(a)
implies that $\widetilde\expec[\SS_{t_k}]\geq c t_k u^{\tau-3}$ for $t_k=k\delta_u\in [\vep,u/2]$.
Therefore (using $t_k=k\delta_u$)
	\eqn{
	\widetilde \prob(E_u, \exists \, t\in I_k\colon
	\SS_t-\widetilde\expec[\SS_t]\geq 2\eta \widetilde\expec[\SS_t])
	\leq (1+o(1))\e^{-ck\delta_u u^{\tau-3}},
	}
which is small even when summed out over $k$ as above.
%\todo[inline]{check (S: on first glance the above does not seem summable, but it is, as one has to consider the range of $k$ to make it work. So I added ``as above''.)}

\paragraph{Part 2: The case $t_k\leq u/2$ and a large lower deviation.}
We continue with bounding the probability that there exists a $t\in I_k=[t_{k},t_{k+1}]$, $\epsilon \le t_k \le u/2$ such that
$\SS_t-\widetilde\expec[\SS_t]\leq -2\eta \widetilde\expec[\SS_t]$, which is slightly more involved.
Again using that $\widetilde\expec[\SS_{t}]=\widetilde\expec[\SS_{t_k}](1+o(1))$ for any $t\in I_k$ by Corollary \ref{cor-asympt-mean-S}(c), we bound
	\eqn{
	\widetilde \prob(E_u, \exists \, t\in I_k\colon
	\SS_t-\widetilde\expec[\SS_t]\leq -2\eta \widetilde\expec[\SS_t])
	\leq \widetilde \prob(\exists \, t \leq \delta_u\colon
	\SS_{t+t_k}-\SS_{t_k}\leq -\frac{\eta}{2}\widetilde\expec[\SS_{t_k}]).
	}
Further,
	\eqan{
	\lbeq{split-Stk-LD}
	\SS_{t+t_k}-\SS_{t_k}&=
	\betam t+\sum_{i=2}^{\infty} c_i [\II_i(t+t_k)-\II_i(t_k)-c_i t]\\
	&=
	\betam t+\sum_{i=2}^{\infty} c_i [1-\II_i(t_k)]\big[\II_i(t+t_k)-\II_i(t_k)-c_i t\big]
	-t\sum_{i=2}^{\infty} c_i^2 \II_i(t_k)\nn\\
	&=\betam t+\sum_{i=2}^{\infty} c_i [1-\II_i(t_k)] [N_i(t)-c_i t]
	-\sum_{i=2}^{\infty}
	c_i [1-\II_i(t_k)] [N_i(t)-\indic{N_i(t)\geq 1}]-t\sum_{i=2}^{\infty} c_i^2 \II_i(t_k)\nn\\
	&\geq \betam t + \RR_t'-\DD_t-\delta_u\EE_u,\nn
	}
where we set
	\eqan{
	\RR_t'&=\sum_{i=2}^{\infty} c_i [1-\II_i(t_k)] [N_i(t)-c_i t]
	}
and
	\eqn{
	\DD_t=\sum_{i=2}^{\infty} c_i [1-\II_i(t_k)] [N_i(t)-\indic{N_i(t)\geq 1}],
	\quad
	\EE_u=\sum_{i=2}^{\infty} c_i^2 \II_i(u).
	}
Thus, conditionally on $(\II_i(t_k))_{i\geq 1}$, the process $(\RR_t')_{t\geq 0}$
is a L\'evy process similar to the L\'evy process investigated in Part 1 above,
$\DD_t$ is the contribution due to $i$ for which $N_i(t)\geq 2$,
while $\EE_u$ yields an upper  bound for the decrease in the drift of our process.
We deal with the three terms one by one, starting with $(\RR_t')_{t\geq 0}$.
As in the previous part,
	\eqn{
	\widetilde \prob(\exists \, s\leq \delta_u\colon
	\betam s + \RR_s'\leq -\frac{\eta}{4}\widetilde\expec[\SS_{t_k}])
	}
is small enough even when summed out over $k$ such that $t_k\in [\vep,u/2]$. This again follows by
Doob's inequality and the bound that for any $\lambda\geq 0$, and with $\FF_{t_k}$ the $\sigma$-algebra generated by $(\SS_t)_{t\in[0,t_k]}$,
	\eqan{
	\lbeq{part-2-ref-deltauEu}
	\widetilde \prob(\exists \, s\leq t\colon
	\betam s + \RR_s'\leq -x \mid \FF_{t_k})&\leq \widetilde \prob(\exists \, s\leq t\colon
	\MM_s'(-\lambda)\geq \e^{x\lambda-t\phi'(-\lambda)-t |\betam| \lambda}\mid \FF_{t_k})\\
	&\leq \e^{-[2x\lambda-t\phi'(-\lambda)-t |\betam| \lambda]}
	\widetilde\expec[\MM_t'(-\lambda)^2 \mid \FF_{t_k}]=\e^{-[2x\lambda-t\phi'(-2\lambda)-2t |\betam| \lambda]},\nn
	}
where
	\eqn{
	\lbeq{def-MM-prime}
	\MM_t'(\lambda)= \e^{\lambda \RR_t'-t\phi'(\lambda)},
	\qquad
	\text{with}
	\qquad
	\phi'(\lambda)=\log{\widetilde \expec[\e^{\lambda \RR'_1}\mid \FF_{t_k}]}.
	}
We compute that
	\eqn{
	\lbeq{comp-phi-prime}
	\phi'(\lambda)
	=\sum_{i=2}^{\infty} \log{\widetilde \expec[\e^{\lambda c_i[1-\II_i(t_k)][N_i(1)-c_i]}\mid \FF_{t_k}]}
	=\sum_{i=2}^{\infty} c_i[1-\II_i(t_k)]\big(\e^{\lambda c_i}-1-\lambda c_i\big).
	}	
Now follow the same steps as in Part 1, using that $0 \leq \phi'(-2) \leq \mbox{const.}$
By Lemma \ref{lem-LD-par}, the term $\delta_u\EE_u$ is, with probability at least $1-Cu^{-(\tau-1)}$
bounded by $\delta_u K u^{\tau-3}=K/u$, which is $o(\widetilde \expec[\SS_{t_k}])$ as $\widetilde\expec[\SS_{t_k}]\geq c t_k u^{\tau-3}$ for $t_k\in [\vep,u/2]$ by Corollary \ref{cor-asympt-mean-S}(a).
We continue to bound $\DD_t$ by
bounding
	\eqn{
	\widetilde \prob(\exists \, t\leq \delta_u\colon
	\DD_t\geq \eta\widetilde\expec[\SS_{t_k}]/4)
	=\widetilde \prob(\DD_{\delta_u}\geq \eta\widetilde\expec[\SS_{t_k}]/4),
	}
since the process $t\mapsto \DD_t$ is non-decreasing.
By the Markov inequality,
	\eqn{
	\lbeq{part-2-ref-Dt}
	\widetilde \prob(\DD_{t}\geq x)\leq x^{-1} \widetilde\expec[\DD_t]
	\leq x^{-1} \sum_{i=2}^{\infty} c_i \widetilde\expec[N_i(t)-\indic{N_i(t)\geq 1}]
	\leq cx^{-1} \sum_{i=2}^{\infty} c_i (c_i t)^2\leq c t^2/x.
	}
Applying this to $x=\eta\widetilde\expec[\SS_{t_k}]/4$ with
$\widetilde\expec[\SS_{t_k}] \geq c t_k u^{\tau-3}$
and $t=\delta_u=u^{-(\tau-2)}$ yields
	\eqn{
	\widetilde \prob(\exists \, t\leq \delta_u\colon \DD_{t}\geq \eta\widetilde\expec[\SS_{t_k}]/4)\leq
	c u^{-2(\tau-2)} (\eta t_k)^{-1} u^{-(\tau-3)}
	=c k^{-1} u^{-(2\tau-5)}.
	}
%\todo{check}
When summing this out over $k$ such that $t_k=k\delta_u\in [\vep,u/2]$
we obtain a bound $c(\log{u}) u^{-(2\tau-5)}=o(u^{-(\tau-1)/2}),$
since $(2\tau-5)>(\tau-1)/2$ precisely when $\tau>3$. This proves
that
	\eqn{
	\sum_{k\colon k\delta_u\in [\vep,u/2]}
	\widetilde \prob(\DD_{\delta_u}\geq \eta\widetilde\expec[\SS_{t_k}]/4)
	=o(u^{-(\tau-1)/2}),
	}
as required. Collecting terms completes Part 2.

\paragraph{Part 3: The case $t_k\geq u/2$ and a large upper deviation.}
This proof is more subtle. We fix $k$ such that $t_k\in [u/2, u-T\delta_u]$ and
condition on $\FF_{t_k}$, which is the $\sigma$-field generated by $(\SS_t)_{t\leq t_k}$
to write (recall \eqref{def-Eu-prob-dev-far-all-times})
	\eqan{
	\lbeq{part-3-cond}
	&\widetilde \prob\big(E_u, \SS_u\in[0,M/u], \exists \, t\in I_k\colon
	\SS_t-\widetilde\expec[\SS_t]\geq 2\eta \widetilde\expec[\SS_t]\big)\\
	&\qquad\leq
	\widetilde \expec\Big[\indic{|\SS_{t_k}-\widetilde\expec[\SS_{t_k}]|
	\leq \eta \widetilde\expec[\SS_{t_k}],
	\sum_{i=2}^{\infty} c_i [1-\II_i(u/2)]\big(\e^{\lambda c_i}-1-\lambda c_i\big)
	\leq K\lambda^2u^{\tau-4}}\widetilde \prob\big(\exists \, t\in I_k\colon
	\SS_t-\widetilde\expec[\SS_t]\geq 2\eta \widetilde\expec[\SS_t]\mid \FF_{t_k}\big)\Big].\nn
	}
First observe that on $\{ |S_{t_k}-\widetilde\expec[\SS_{t_k}]| \leq \eta \widetilde\expec[\SS_{t_k}]\}$, we have
  \eqn{
  \widetilde\prob\big(\exists \, t\in I_k\colon
	\SS_t-\widetilde\expec[\SS_t]\geq 2\eta \widetilde\expec[\SS_t]\mid \FF_{t_k}\big)
	\le \widetilde \prob\big(\exists \, t\leq \delta_u\colon
	\SS_{t+t_k}-\SS_{t_k}\geq \frac{\eta}{2}\widetilde\expec[\SS_{t_k}] \mid \FF_{t_k}\big)
  }
by using that $\widetilde\expec[\SS_{t}]=\widetilde\expec[\SS_{t_k}](1+o_T(1))$ for any $t\in I_k$
by Corollary \ref{cor-asympt-mean-S}(d).
Similar to \refeq{split-Stk-LD}, we bound from above
	\eqan{
%	\lbeq{split-Stk-LD2}
	\SS_{t+t_k}-\SS_{t_k}
	&=
	\betam t+\sum_{i=2}^{\infty} c_i [1-\II_i(t_k)]\big[\II_i(t+t_k)-\II_i(t_k)-c_i t\big]
	-t\sum_{i=2}^{\infty} c_i^2 \II_i(t_k)\\
	&\leq \betam t+\sum_{i=2}^{\infty} c_i [1-\II_i(t_k)] [N_i(t)-c_i t] =\betam t + \RR_t',\nn
	}
where we note that $\RR_t'$ is as in Part 2.
Conditionally on $\FF_{t_k}$, the process $(\RR_t')_{t\geq 0}$ is a L\'evy process, and we use
	\eqan{
	\widetilde \prob(\exists \, s\leq t\colon
	\betam s + \RR_s'\geq x\mid \FF_{t_k})&\leq \widetilde \prob(\exists \, s\leq t\colon
	\MM_s'(\lambda)\geq \e^{x\lambda-t\phi''(\lambda)-t |\betam| \lambda}\mid \FF_{t_k})\\
	&\leq \e^{-2[x\lambda-t\phi''(\lambda)-t|\betam|\lambda]}
	\widetilde\expec[\MM_t'(\lambda)^2\mid \FF_{t_k}]
	=\e^{-2x\lambda+t\phi''(2\lambda)+2t|\betam|\lambda},\nn
	}
where we recall equations \eqref{def-MM-prime} and \eqref{comp-phi-prime}.
Since $\e^{\lambda c_i}-1-\lambda c_i\geq 0$ for every $\lambda\in {\mathbb R}$, and
since $1-\II_i(t_k)\leq 1-\II_i(u/2)$ for every $t_k\geq u/2$, a.s.\
	\eqn{
	\phi'(\lambda)
	\leq
	\sum_{i=2}^{\infty} c_i[1-\II_i(u/2)]\big(\e^{\lambda c_i}-1-\lambda c_i\big).
	}
On the event $\{\sum_{i=2}^{\infty} c_i [1-\II_i(u/2)]\big(\e^{\lambda c_i}-1-\lambda c_i\big)
	\leq K\lambda^2u^{\tau-4}\}$ (recall \eqref{part-3-cond}), we have that
$\phi'(\lambda)\leq K\lambda^2u^{\tau-4}$, so that we can further bound, choosing $\lambda=\delta u$ and $t=\delta_u=u^{-(\tau-2)}$,
	\eqn{
	\widetilde \prob(\exists \, s\leq t\colon
	\betam s + \RR_s'\geq x\mid \FF_{t_k})
	\leq \e^{-2x\lambda+tK4\lambda^2u^{\tau-4}+2t|\betam|\lambda}
	\leq \e^{-2x\delta u+K4\delta^2 +2 u^{-(\tau-3)} |\betam| \delta}.
	}
We take $x=\widetilde\expec[\SS_{t_k}]$ and note that Corollary \ref{cor-asympt-mean-S}(b) and Lemma \ref{lem-expec-S}(d) yield that $\widetilde\expec[\SS_{t_k}]\geq c(u-t_k) u^{\tau-3}$ for
$t_k\in [u/2,u-T\delta_u]$. Then,
	\eqn{
	\widetilde \prob(\exists \, s\leq t\colon
	\betam s + \RR_s'\geq x\mid \FF_{t_k})
	\leq c \e^{-c\delta (u-t_k) u^{\tau-2}}.
	}
Summing over $k$ with $t_k=k \delta_u\in [u/2, u-T\delta_u]$ and $\delta_u=u^{-(\tau-2)}$, using Proposition \ref{prop-dens-f} and $\widetilde\expec[\SS_{t_k}]\leq c(u-t_k) u^{\tau-3}$ by Corollary \ref{cor-asympt-mean-S}(b) and Lemma \ref{lem-expec-S}(d)
yields as an upper bound (recall also \eqref{part-3-cond} and the definition of $E_u$ from \eqref{def-Eu-prob-dev-far-all-times})
%\todo[inline]{S-? Why can't we just bound the probability in the second line to come by $1$? Remco: Then  the powers of $u$ do not add up... It fact, it is quite delicate and works {\em precisely}. Nothing can be thrown away.}
	\eqan{
	&\sum_{k\colon t_k\in [u/2, u-T\delta_u]}
	\widetilde \prob\big(E_u, \SS_u\in[0,M/u], \exists \, t\in I_k\colon
	\SS_t-\widetilde\expec[\SS_t]\geq 2\eta \widetilde\expec[\SS_t]\big)\\
	&\qquad \leq c \sum_{k\colon t_k\in [u/2, u-T\delta_u]}
	\widetilde \prob(|\SS_{t_k}-\widetilde\expec[\SS_{t_k}]|\leq \eta \widetilde\expec[\SS_{t_k}])
	\e^{-c\delta (u-k\delta_u) u^{\tau-2}}\nn\\
	&\qquad\leq c \sum_{k\colon t_k\in [u/2, u-T\delta_u]}
	\eta \widetilde\expec[\SS_{t_k}]
	\Big(\sup_w \widetilde f_{\SS_{t_k}}(w)\Big)\e^{-c\delta (u-k\delta_u) u^{\tau-2}}\nn\\
	&\qquad\leq c u^{-(\tau-3)/2} \sum_{k\colon t_k\in [u/2, u-T\delta_u]} C(u-k\delta_u) u^{\tau-3}
	\e^{-c\delta (u-k\delta_u) u^{\tau-2}}
	\nn\\
	&\qquad\leq
	c u^{-(\tau-1)/2} \sum_{k\colon t_k\in [u/2, u-T\delta_u]} C(u-k\delta_u) u^{\tau-2}
	\e^{-c\delta (u-k\delta_u) u^{\tau-2}}\nn\\
	&\qquad\leq
	cu^{-(\tau-1)/2}\e^{-c\delta T}=o_{\sss T}(1)u^{-(\tau-1)/2},\nn
	}
as required.

\paragraph{Part 4: The case $t_k\geq u/2$ and a large lower deviation.}
We again start from \refeq{split-Stk-LD}, and note that the bounds on
$\DD_t$ and $\delta_u\EE_u$ proved in Part 2 still apply, now using that by Corollary \ref{cor-asympt-mean-S}(b) and Lemma \ref{lem-expec-S}(d) $\widetilde\expec[\SS_{t_k}]\geq c(u-t_k) u^{\tau-3}$ for $t_k\in [u/2, u-T\delta_u]$ with $\delta_u=u^{-(\tau-2)}$ below \eqref{part-2-ref-Dt}. We further use this estimate to replace the statement that $\delta_u \EE_u$ is $o(\widetilde \expec[\SS_{t_k}])$ from below \eqref{part-2-ref-deltauEu} by $\delta_u \EE_u \leq o_{\sss T}(1) \widetilde\expec[\SS_{t_k}]$. The exponential
martingale bound for $\RR_t'$ performed in Part 3 can easily be adapted to
deal with a large lower deviation as well. We omit further details.
\qed

\medskip

\noindent
{\it Proof of Proposition \ref{prop-no-middle-ground}.}
The proof follows by combining Propositions \ref{prop-integer-times} and \ref{prop-no-hit-small-interval}.
Indeed, choose $\eta=1/4$ and observe that $\widetilde\expec[\SS_t] > 0$ on $[\vep,u-T\delta_u]$, using that $\widetilde\expec[\SS_t]\geq c t u^{\tau-3}$ for $t\in [\vep,u/2]$ and $\widetilde \expec[\SS_t]\ge c(u-t)u^{\tau-3}$ for all $t \in [u/2,u-T \delta_u]$ by Corollary \ref{cor-asympt-mean-S}(a),(b) and Lemma \ref{lem-expec-S}(d).
\qed

\section{Conditional expectations given $u\SS_u=v$}
\label{sec-cond-Suv}

A major difficulty in the proof of Proposition
\ref{prop-weak-cond-conv} is the fact that, while the summands
in the definition of $Q_u(t)$ in \refeq{Qu-def} are
independent, this property is lost due to the fact that
we \emph{condition on $\SS_u$}. The following lemma
allows us to deal with such expectations:

\begin{Lemma}[Conditional expectations given a continuous random variable]
\label{lem-cond-expec-SSu}
Let $\FSS((\SS_s)_{s\geq 0})$ be a functional of the process $(\SS_s)_{s\geq 0}$
such that $\FSS((\SS_s)_{s\geq 0})\geq 0$ $\widetilde\prob$-a.s.,
and $0<\widetilde\expec[\FSS((\SS_s)_{s\geq 0})]<\infty$.
Then, for every $w\in \R$,
    \eqn{
    \lbeq{FT-cond-expec-1}
    \widetilde\expec\big[\FSS((\SS_s)_{s\geq 0})\mid \SS_u=w\big]=
    \frac{1}{\widetilde f_{\SS_u}(w)} \int_{-\infty}^{+\infty}
    \e^{-\im k w} \widetilde\expec\big[\FSS((\SS_s)_{s\geq 0})\e^{\im k\SS_u}\big]
    \frac{dk}{2\pi},
    }
where $\im$ denotes the imaginary unit.
\end{Lemma}

For $\FSS((\SS_s)_{s\geq 0})=1$, \refeq{FT-cond-expec-1} is just the usual Fourier inversion
theorem applied to the (continuous) random variable $\SS_u$.
The expectation $\widetilde\expec\big[\FSS((\SS_s)_{s\geq 0})\e^{\im k\SS_u}]$
\emph{factorizes} when $\FSS((\SS_s)_{s\geq 0})$ is of product form
in the underlying random variables $(\II_i(s))_{s\geq 0}$.
In our applications, $\widetilde\expec\big[\FSS((\SS_s)_{s\geq 0})\mid \SS_u=w\big]$
will be close to constant in $w$. Then, in order to compute its asymptotics,
it suffices to check that the computation in the proof of Proposition \ref{prop-dens-f}
is hardly
affected by the presence of $\FSS((\SS_s)_{s\geq 0})$.

\proof Define the measure $\widetilde\prob^{\sss \FSS}$ by
    \eqn{
    \lbeq{prob-F-def}
    \widetilde\prob^{\sss \FSS}(E)=\frac{\widetilde\expec\big[\FSS((\SS_s)_{s\geq 0})\indicwo{E}\big]}{\widetilde\expec\big[\FSS((\SS_s)_{s\geq 0})\big]}.
    }
Under the measure $\widetilde\prob^{\sss \FSS}$, the random variable
$\SS_u$ is again continuous, since $0<\widetilde\expec[\FSS((\SS_s)_{s\geq 0})]<\infty$.
Let $\widetilde f_{\SS_u}^{\sss \FSS}$ denote the density of $\SS_u$ under the measure
$\widetilde\prob^{\sss \FSS}$. Then, we obtain, by the Fourier inversion theorem
applied to $\widetilde\prob^{\sss \FSS}$, that
    \eqn{
    \lbeq{fF-FI-res}
    \widetilde f_{\SS_u}^{\sss \FSS}(w)=\int_{-\infty}^{+\infty}
    \e^{-\im k w} \widetilde\expec^{\sss \FSS}\big[\e^{\im k\SS_u}]\frac{dk}{2\pi}.
    }
Now, by \refeq{prob-F-def},
    \eqn{
    \widetilde f_{\SS_u}^{\sss \FSS}(w)=\frac{\widetilde\expec\big[\FSS((\SS_s)_{s\geq 0})\mid \SS_u=w\big]}
    {\widetilde\expec\big[\FSS((\SS_s)_{s\geq 0})\big]} \widetilde f_{\SS_u}(w),
    }
while
    \eqn{
    \widetilde\expec^{\sss \FSS}\big[\e^{\im k\SS_u}]=\frac{\widetilde\expec\big[\FSS((\SS_s)_{s\geq 0})\e^{\im k\SS_u}]}
    {\widetilde\expec\big[\FSS((\SS_s)_{s\geq 0})\big]}.
    }
Therefore, substituting both sides in \refeq{fF-FI-res} and multiplying through by
$\widetilde\expec\big[\FSS((\SS_s)_{s\geq 0})\big]$ proves the claim.
\qed
\medskip

Let $\widetilde\prob_v$ denote $\widetilde \prob$
conditionally on $u\SS_u=v$, so that Lemma \ref{lem-cond-expec-SSu} implies that
    \eqn{
    \lbeq{FT-cond-expec-2}
    \widetilde\expec_v\big[\FSS((\SS_s)_{s\geq 0})\big]=
    \frac{1}{\widetilde f_{\SS_u}(v/u)} \int_{-\infty}^{+\infty}
    \e^{-\im k v/u} \widetilde \expec\big[\FSS((\SS_s)_{s\geq 0})\e^{\im k\SS_u}]\frac{dk}{2\pi}.
    }
In many cases, it shall prove to
be convenient to rewrite the above using
    \eqn{
    \widetilde \expec\big[\FSS((\SS_s)_{s\geq 0})\e^{\im k\SS_u}\big]
    =\widetilde \expec\Big[\e^{\im k\SS_u}\widetilde \expec\big[\FSS((\SS_s)_{s\geq 0})\mid \JJ(u)\big]\Big],
    }
since the random variables $(T_i)_{i\in \JJ(u)}$
are, conditionally on $\JJ(u)$, independent with
    \eqn{
    \lbeq{cond-distr-T}
    \widetilde \prob(T_i\leq u-t \mid T_i\leq u)
    =\frac{1-\e^{-c_i(u-t)}}{1-\e^{-c_iu}}.
    }
In the following lemma, we investigate the effect on $\prob(i\in \JJ(u))$ of conditioning on
$\SS_u=w$:

\begin{Lemma}[The set $\JJ(u)$ conditionally on $\SS_u=w$]
\label{lem-difference-probs}
%\ \\
%(i)
There exists a constant $d>0$ such that for any $i$ and $w=o(u^{(\tau-3)/2})$,
    \eqn{
    \Big|\widetilde \prob(j\in \JJ(u) \mid \SS_u=w) -\widetilde \prob(j\in \JJ(u))\Big| \le d c_i \widetilde \prob(j\in \JJ(u)) \widetilde \prob(j\not\in \JJ(u)) u^{-(\tau-3)/2}.
    }
%(ii) Similarly, there exists a constant $d>0$ such that for any $i\neq j$ and $w=o(u^{(\tau-3)/2})$,
%    \eqan{
%    & \Big|\widetilde \prob(i\in \JJ(u),j\in \JJ(u) \mid \SS_u=w) - \widetilde \prob(i\in \JJ(u),j\in \JJ(u))  \Big| \\
%    &\le d u^{-(\tau-3)}c_ic_j\widetilde\prob(i\in \JJ(u),j\in \JJ(u)) \prob(i\not\in \JJ(u),j\not\in \JJ(u)). \nn
%    }
\end{Lemma}

%\todotemp{Remco: I do not trust (ii)!! We also do not give its proof. Finally, we do not seem to use it???}
\proof
By Lemma \ref{lem-cond-expec-SSu} (for the second term use $\FSS \equiv 1$)
    \eqan{
    & \left| \widetilde \prob(j\in \JJ(u) \mid \SS_u=w) -\widetilde \prob(j\in \JJ(u)) \right| \\
    &= \frac{1}{\widetilde f_{\SS_u}(w)} \left| \int_{-\infty}^{+\infty}
    \e^{-\im k w} \widetilde \expec\!\left[ \left( \indic{j\in \JJ(u)} -\widetilde \prob(j\in \JJ(u)) \right) \e^{\im k\SS_u} \right] \frac{dk}{2\pi} \right| \nn\\
    &= \frac{1}{u^{(\tau-3)/2} \widetilde f_{\SS_u}(w)} \left| \int_{-\infty}^{+\infty}
    \e^{-\im k u^{-(\tau-3)/2} w} \widetilde \expec\!\left[ \left( \indic{j\in \JJ(u)} -\widetilde \prob(j\in \JJ(u)) \right) \e^{\im k u^{-(\tau-3)/2} \SS_u} \right] \frac{dk}{2\pi} \right|. \nn
    }
Recall Lemma \ref{lem-ind-tilt}. Under the measure $\widetilde \prob$, the distribution of the indicator processes
$(\II_j(t))_{t\geq 0}$ is that of independent indicator processes.
Define $\SS_u^{\sss (j)}=\SS_u-c_j(\II_j(u)-c_ju)$. By \eqref{SS-def-ref} and \eqref{index-def}, the random variables $\II_j(u)$ and $\SS_u^{\sss (j)}$ are independent under $\widetilde \prob$. This yields
    \eqan{
    \lbeq{diff-prob-1}
    & \left| \widetilde \prob(j\in \JJ(u) \mid \SS_u=w) -\widetilde \prob(j\in \JJ(u)) \right| \\
    &\quad\leq \frac{1}{u^{(\tau-3)/2} \widetilde f_{\SS_u}(w)} \int_{-\infty}^{+\infty}
    \left| \widetilde \expec\!\left[ \left( \indic{j\in \JJ(u)} -\widetilde \prob(j\in \JJ(u)) \right) \e^{\im k u^{-(\tau-3)/2} c_j(\II_j(u)-c_ju)} \right] \right| \left| \widetilde \expec\!\left[ \e^{\im k u^{-(\tau-3)/2} \SS_u^{(j)}} \right] \right| \frac{dk}{2\pi} \nn\\
    &\quad= \frac{1}{u^{(\tau-3)/2} \widetilde f_{\SS_u}(w)} \int_{-\infty}^{+\infty}
    \widetilde \prob(j\not\in \JJ(u))\widetilde \prob(j\in \JJ(u)) \left| \e^{\im k u^{-(\tau-3)/2} c_j} - 1 \right|
    \left| \widetilde \expec\!\left[ \e^{\im k u^{-(\tau-3)/2} \SS_u^{\sss (j)}} \right] \right| \frac{dk}{2\pi}. \nn
    }
Next we claim that there exist constants $C_1, C_2$ such that for all $j \geq 2$
    \eqn{
    \lbeq{diff-prob-2}
    \left| \widetilde \expec\!\left[ \e^{\im k u^{-(\tau-3)/2} \SS_u^{\sss (j)}} \right] \right|
    \leq C_1 \e^{-C_2 |k|^{\tau-2}}.
    }
Indeed, for $\SS_u^{\sss (j)}$ replaced by $\SS_u$ the result was derived in the proof of Proposition \ref{prop-dens-f} in \cite{AidHofKliLee13b}. To prove the same for $\SS_u^{\sss (j)}$ with $j \geq 2$ arbitrary, and following the approach in \cite{AidHofKliLee13b},
%we only have to adapt the estimates in \eqref{density-1} and \eqref{density-2}.
% Proceeding as in \eqref{density-1},
 we obtain for $\frac{k}{2\pi} u^{-(\tau-3)/2-1} \leq 1/8$ the bound
    \eqan{
    \log\!\left( \left| \widetilde \expec\!\left[ \e^{\im k u^{-(\tau-3)/2} \SS_u^{(j)}} \right] \right| \right)
    &\leq -c u^{4-\tau} k^2 \sum_{i\geq 2: c_i<1/u, i\neq j} c_i^3\leq -c_0 |k|^{\tau-2} + c u^{4-\tau} k^2 c_j^3 \indic{c_j<1/u}\\
    &\leq -c_0 |k|^{\tau-2} + c u^{-(\tau-1)} k^2\leq -c_0 |k|^{\tau-2} + c, \nn
    }
%while in \eqref{density-2} and
while for $y_k=8\frac{k}{2\pi}u^{-(\tau-3)/2}>u$,
    \eqan{
    \log\!\left( \left| \widetilde \expec\!\left[ \e^{\im k u^{-(\tau-3)/2} \SS_u^{(j)}} \right] \right| \right)
    &\leq -c_0 |k|^{\tau-2} + cu^{4-\tau}k^2 c_j^3 \indic{c_j<1/y_k}\leq -c_0 |k|^{\tau-2} + cu^{4-\tau}k^2 y_k^{-3}\\
    &\leq -c_0 |k|^{\tau-2} + cu^{4-\tau}u^{(\tau-3)} u^{-1}= -c_0 |k|^{\tau-2} + c. \nn
    }
Substituting  \eqref{diff-prob-2} in \eqref{diff-prob-1} yields
    \eqan{
    & \left| \widetilde \prob(j\in \JJ(u) \mid \SS_u=w) -\widetilde \prob(j\in \JJ(u)) \right| \\
    &\leq \frac{1}{u^{(\tau-3)/2} \widetilde f_{\SS_u}(w)} \int_{-\infty}^{+\infty}
    \widetilde \prob(i\not\in \JJ(u)) \widetilde \prob(i\in \JJ(u)) \left| \e^{\im k u^{-(\tau-3)/2} c_j} - 1 \right| C_1 \e^{-C_2 |k|^{\tau-2}} \frac{dk}{2\pi}. \nn
    }
We further have
    \eqn{
    \left| \e^{\im k u^{-(\tau-3)/2} c_j} - 1 \right|
    = \left( 2(1-\cos(k u^{-(\tau-3)/2} c_j)) \right)^{1/2}
    \leq \sqrt{2} |k| u^{-(\tau-3)/2} c_j,
    }
which yields
    \eqan{
    & \left| \widetilde \prob(j\in \JJ(u) \mid \SS_u=w) -\widetilde \prob(j\in \JJ(u)) \right| \\
    &\leq C_3 \frac{1}{u^{(\tau-3)/2} \widetilde f_{\SS_u}(w)} \int_{-\infty}^{+\infty}
    \widetilde \prob(i\not\in \JJ(u)) \widetilde \prob(i\in \JJ(u)) k u^{-(\tau-3)/2} c_j \e^{-C_2 |k|^{\tau-2}} \frac{dk}{2\pi} \nn\\
    &= C_3 \widetilde \prob(i\not\in \JJ(u))\widetilde \prob(i\in \JJ(u)) u^{-(\tau-3)/2} c_j \frac{1}{u^{(\tau-3)/2} \widetilde f_{\SS_u}(w)} \int_{-\infty}^{+\infty} k \e^{-C_2 |k|^{\tau-2}} \frac{dk}{2\pi}. \nn
    }
For $w=o(u^{(\tau-3)/2})$ and by Proposition \ref{prop-dens-f}, $u^{(\tau-3)/2} \widetilde f_{\SS_u}(w)=B(1+o(1))$ uniformly in $w$ and the claim in (i) follows. %The claim in (ii) follows similarly.
\qed
\medskip

\begin{Corollary}
\label{cor-difference-probs}
There exists a constant $C>0$ such that for any $i$ and $w=o(u^{(\tau-3)/2})$,
    \eqn{
    \widetilde \prob(i\in \JJ(u) \mid \SS_u=w)
    \leq C (1\wedge c_iu).
    }
\end{Corollary}

\proof
The bound by 1 is obvious. The bound by $Cc_iu$ follows once we
recall \eqref{Ti-def-tilt} and observe that for $c_j\leq 1/u$,
$\widetilde \prob(T_j \leq u) = \widetilde \prob(j\in \JJ(u))
\leq C(\tau) c_j u$. Now use Lemma \ref{lem-difference-probs}(i).
\qed
\medskip

%%%%%%%%%%%%%%%%%%%%%%%%%%%%%%%

\section{The near-end ground: Proof of Proposition \ref{prop-weak-cond-conv}}
\label{section-pf-A-B}
In this section, we prove Proposition \ref{prop-weak-cond-conv}.
The proof is divided into several key parts. In Section \ref{sec-conv-A}, we show convergence
of the mean process $A_u$ in Proposition \ref{prop-weak-cond-conv}(a). In Section \ref{sec-conv-B},
we prove the convergence of $B_u$ in Proposition \ref{prop-weak-cond-conv}(b).

%%%%%%%%%%%%%%%%%%%%%%%%%%%%%%%%%%%%%%%%%%%%%%%%%%%%%%%

\subsection{Convergence of the mean process $A_u$}
\label{sec-conv-A}

%%%%%%%%%%%%%%%%%%%%%%%%%%%%%%%%%%%%%%%%%%%%%%%%%%%%%%%

Recall the definition of $A_u$ from \refeq{AB-defs}. By \refeq{cond-distr-T},
    \eqn{
    A_u(tu^{-(\tau-2)})
    =-\sum_{j\in \JJ(u)} c_j \left[ u\frac{\e^{-c_j(u-tu^{-(\tau-2)})}-\e^{-c_ju}}{1-\e^{-c_ju}}-tu^{-(\tau-2)} \right].
    }
We use that $|\e^x-1-x| \leq \e^D x^2/2$ for $0 \leq x \leq D$ with
$x=c_jtu^{-(\tau-2)}$, where for $0 \leq t \leq T$, $c_jtu^{-(\tau-2)} \leq tu^{-(\tau-2)} \leq \mbox{const.}$,
to obtain
    \eqn{
    \lbeq{A-form-1}
    A_u(tu^{-(\tau-2)})
    = -\sum_{j\in \JJ(u)} c_jtu^{-(\tau-2)} \left[ c_ju\frac{\e^{-c_ju}}{1-\e^{-c_ju}}-1 \right]
    +E_u(t)
    }
with an error term $E_u(t)$ bounded by
    \eqan{
    \lbeq{error-Taylor}
    |E_u(t)|
    &\leq C \sum_{j\in \JJ(u)} \left( c_jtu^{-(\tau-2)} \right)^2 c_ju \frac{\e^{-c_ju}}{1-\e^{-c_ju}} \\
    &\leq C T^2 u^{-2(\tau-2)}
    \Big[u\sum_{j\in \JJ(u)\colon c_j>1/u} c_j^3 +
    \sum_{j\in \JJ(u)\colon c_j\leq 1/u} c_j^2\Big],\nn
    }
uniformly in $t \leq T$. Since $\sum_{j \geq 2} c_j^3 < \infty$ and $u^{-2(\tau-2)+1}=u^{5-2\tau}=o(1)$,
the first term vanishes. Further, by Corollary \ref{cor-difference-probs} with $w=v/u$,
    \eqn{
    u^{-2(\tau-2)}\widetilde \expec_v\Big[\sum_{j\in \JJ(u)\colon c_j\leq 1/u} \shift c_j^2\Big]
    =u^{-2(\tau-2)}\shift \sum_{j\in \JJ(u)\colon c_j\leq 1/u} \shift c_j^2\widetilde \prob_v(j\in \JJ(u))
    \leq  u^{5-2\tau}
    \shift
	\sum_{j\in \JJ(u)\colon c_j\leq 1/u}
	\shift
	c_j^3=o(1),
    }
so that also the second term is $o_{\sss \widetilde \prob_v}(1)$.

In the above proof, we see that it is useful to split
a sum over $j \in \JJ(u)$ into $j\in \JJ(u)$ such that
$c_j>1/u$ and $j\in \JJ(u)$ such that $c_j \leq 1/u$. Then we use upper bounds similar to the ones in Corollary~\ref{cor-difference-probs} to bound the
arising sums. We will follow this strategy often below.

We further rewrite \eqref{A-form-1} into
    \eqn{
    %\lbeq{A-and-q}
    A_u(tu^{-(\tau-2)})
    = t \sum_{j\in \JJ(u)} q_j(u)+E_u(t)
    \qquad
	\text{with}
	\qquad
    \lbeq{def-qj}
    q_j(u) \equiv u^{-(\tau-2)}c_j\e^{-c_ju}\frac{\e^{c_ju}-1-c_ju}{1-\e^{-c_ju}}.
    }
Note that $0 \leq q_j(u) \leq 1$ for $u$ big. Below, we will frequently rely
on the bounds
    \eqn{
    \lbeq{qj-approx}
    q_j(u)
    \leq C(\tau) u^{-(\tau-2)}c_j (1\wedge c_j u)
    }
and, using \eqref{Ti-def-tilt} for $t=u$,
    \eqn{
	\lbeq{lem-approx}
    \widetilde \prob(T_j \leq u)
    \leq C(\tau)(1\wedge c_j u),
    \qquad
    1-\widetilde \prob(T_j \leq u)
    \leq \e^{-c_ju(1+\theta)}.
    }

By \refeq{error-Taylor}, to prove the claim of Proposition \ref{prop-weak-cond-conv}(a),
it is enough to show that
    \eqn{
    \lbeq{conv-lin-term}
    \kappa_u
    \equiv \sum_{j\in \JJ(u)} q_j(u)
    = \sum_{j\geq 2} \II_j(u) q_j(u)
    \convpv \kappa.
    }
For this, we compute the Laplace transform of $\kappa_u$ under the measure
$\widetilde \prob_v$ using Lemma \ref{lem-cond-expec-SSu} and a change of variable.
For $a\geq 0$,
    \eqn{
    \lbeq{laplace-A}
    \widetilde \expec_v [\e^{-a \kappa_u}]
    =\frac{1}{u^{(\tau-3)/2} \widetilde f_{\SS_u}(v/u)} \int_{-\infty}^{+\infty}
    \e^{-\im k vu^{-(\tau-1)/2}} \widetilde\expec\big[\e^{-a \kappa_u}\e^{\im k u^{-(\tau-3)/2}\SS_u}\big]\frac{dk}{2\pi}.
    }
By Proposition \ref{prop-dens-f}, for each $v>0$, $u^{(\tau-3)/2} \widetilde f_{\SS_u}(v/u)\rightarrow B$.
We aim to use \emph{dominated convergence} on the integral appearing in \refeq{laplace-A},
for which we have to prove (a) pointwise convergence for each $k\in \R$; and (b) a uniform bound that
is integrable. We start by proving pointwise convergence:

%......................................................

\begin{Lemma}[Pointwise convergence]
\label{lem-one-dim-A}
For $a \geq 0$ arbitrary, $v=o(u^{(\tau-1)/2})$, and with $\kappa_u$ as in \eqref{conv-lin-term},
    \eqn{
    \e^{-\im k vu^{-(\tau-1)/2}}\widetilde\expec\big[\e^{-a \kappa_u}\e^{\im k u^{-(\tau-3)/2} \SS_u}\big]
    = \e^{-a \kappa} \e^{-I_{\sss V}(1) k^2/2}+ o(1).
    }
\end{Lemma}

%......................................................

\proof Trivially, $\e^{-\im k vu^{-(\tau-1)/2}}\to 1$
pointwise when $v=o(u^{(\tau-1)/2})$.
To compute $\widetilde\expec\big[\e^{-a \kappa_u}\e^{\im k u^{-(\tau-3)/2}\SS_u}\big]$,
recall the definition of $\SS_u$ from \eqref{SS-def-ref} and recall that the indicator
processes $\II_j(t)=\indic{T_j\leq t}$ are independent under the measure
$\widetilde \prob$ (cf. Lemma \ref{lem-ind-tilt}), to see that
    \eqan{
    \widetilde\expec\big[\e^{-a \kappa_u}\e^{\im k u^{-(\tau-3)/2} \SS_u}\big]
    &=\e^{\im k u^{-(\tau-3)/2}+\betam \im k u^{-(\tau-5)/2}} \\
    &\quad \times\prod_{j\geq 2} \e^{-\im k u^{-(\tau-5)/2}c_j^2}\Big(1+\big(\e^{-a q_j(u)+\im k u^{-(\tau-3)/2} c_j}-1\big) \widetilde \prob(T_j \leq u) \Big). \nn
    }
The remainder of the proof proceeds in three steps.
\paragraph{Step 1: Asymptotic factorization.} We start by proving that
    \eqn{
    \lbeq{limit-laplace-A}
    \widetilde\expec\big[\e^{-a \kappa_u}\e^{\im k u^{-(\tau-3)/2} \SS_u}\big]
    = \e^{-a \widetilde\expec[\kappa_u]} \widetilde\expec\big[\e^{\im k u^{-(\tau-3)/2} \SS_u}\big] + o(1).
    }
To this end, we first use
    \eqan{
    \lbeq{prod-diff-error}
    \Big| \prod_{j \geq 2} a_j - \prod_{j \geq 2} b_j \Big|
    & \leq \sum_{j \geq 2} \prod_{j_1<j} |a_{j_1}| |a_j-b_j| |\prod_{j_2>j} |b_{j_2}| %\\
    %&
	\leq \sum_{j \geq 2} |a_j-b_j| \qquad \mbox{ if } \sup_j (|a_j| \vee |b_j|) \leq 1,%\nn
    }
to get (recall that $q_j(u)\geq 0$)
    \eqan{
    \lbeq{prod-diff-laplace-A}
    \Big| \prod_{j \geq 2} a_j - \prod_{j \geq 2} b_j \Big|
    &= \Big| \widetilde\expec\big[\e^{-a \kappa_u}\e^{\im k u^{-(\tau-3)/2} \SS_u}\big] - \e^{-a \widetilde\expec[\kappa_u]} \widetilde\expec\big[\e^{\im k u^{-(\tau-3)/2} \SS_u}\big] \Big| \\
    &= \Big| \prod_{j\geq 2} \widetilde\expec\big[ \e^{-a \II_j(u) q_j(u) + \im k u^{-(\tau-3)/2} c_j \II_j(u)} \big] - \prod_{j\geq 2} \e^{-a q_j(u)\widetilde \prob(T_j \leq u)} \widetilde\expec\big[ \e^{\im k u^{-(\tau-3)/2} c_j \II_j(u)} \big] \Big| \nn\\
    &\leq \sum_{j \geq 2} \Big| \widetilde\expec\big[ \e^{-a \II_j(u) q_j(u) + \im k u^{-(\tau-3)/2} c_j \II_j(u)} \big] - \e^{-a q_j(u)\widetilde \prob(T_j \leq u)} \widetilde\expec\big[ \e^{\im k u^{-(\tau-3)/2} c_j \II_j(u)} \big] \Big| \nn\\
    &\equiv \sum_{j \geq 2} \Delta_j(-a q_j(u)), \nn
    }
where we abbreviate $q \equiv -a q_j(u) \leq 0$ such that
    \eqn{
    \lbeq{Delta-j-q-def}
    \Delta_j(q)=\left|\widetilde\expec\big[ \e^{\II_j(u)q+ \im k u^{-(\tau-3)/2} c_j \II_j(u)} \big] -
    \e^{q\widetilde \prob(T_j \leq u)} \widetilde\expec\big[ \e^{\im k u^{-(\tau-3)/2} c_j \II_j(u)}\big]
    \right|.
    }
To bound $\Delta_j(q)$, we write
$\e^{\im k u^{-(\tau-3)/2} c_j}=1+(\e^{\im k u^{-(\tau-3)/2} c_j}-1)$ and
use the triangle inequality to bound each summand by
    \eqan{
    \lbeq{a-triangle}
    \Delta_j(q)
    &= \left| \left( 1-\widetilde \prob(T_j \leq u)+\e^{q} \e^{\im k u^{-(\tau-3)/2} c_j} \widetilde \prob(T_j \leq u) \right) \right. \\
    &\quad \left. - \e^{q\widetilde \prob(T_j \leq u)} \left( 1-\widetilde \prob(T_j \leq u)+\e^{\im k u^{-(\tau-3)/2} c_j} \widetilde \prob(T_j \leq u) \right) \right| \nn\\
    &\leq \left| 1-\widetilde \prob(T_j \leq u) +\e^{q} \widetilde \prob(T_j \leq u) - \e^{q\widetilde \prob(T_j \leq u)}\right| \nn\\
    &\quad + \left| \e^{q} - \e^{q\widetilde \prob(T_j \leq u)} \right| \left| \e^{\im k u^{-(\tau-3)/2} c_j} - 1 \right| \widetilde \prob(T_j \leq u). \nn
    }
We can bound
    \eqn{
    \left| \e^{q} - \e^{q\widetilde \prob(T_j \leq u)} \right| \leq |q|\e^{(q\vee 0)}
    \qquad
	\text{and}
	\qquad
    \left| \e^{\im k u^{-(\tau-3)/2} c_j} - 1 \right|
    \leq |k| u^{-(\tau-3)/2} c_j,
    }
which gives a bound $|q|\e^{(q\vee 0)}|k| u^{-(\tau-3)/2} c_j \widetilde \prob(T_j \leq u)$
on the last line of \refeq{a-triangle}.

To bound the first line of \refeq{a-triangle},
we use the error bounds $|\e^{-x}-1+x| \leq |x|^2$ for all $x\geq 0$
to all the exponential functions in it, to obtain
    \eqn{
    \left| 1-\widetilde \prob(T_j \leq u) +\e^q \widetilde \prob(T_j \leq u) - \e^{q\widetilde \prob(T_j \leq u)}\right|\leq C q^2\widetilde \prob(T_j \leq u).
    }
Together, this leads us to
    \eqan{
    \lbeq{diff-laplace-A}
    \Delta_j(-a q_j(u)) = \Delta_j(q)
    &\leq C |q| \e^{(q\vee 0)} \Big(|q|+|k|u^{-(\tau-3)/2} c_j\Big)\leq C(a) q_j(u) \Big(q_j(u)+|k|u^{-(\tau-3)/2} c_j\Big) \equiv \Xi_j.
    }
To prove \eqref{limit-laplace-A}, by \eqref{prod-diff-laplace-A} and \eqref{diff-laplace-A}
it is enough to show that $\sum_{j \geq 2} \Xi_j= o(1).$
Consider the sum over $c_j> 1/u$ first. By \refeq{qj-approx},
    \eqan{
    \sum_{j\geq 2\colon c_j > 1/u} \shift \Xi_j
    &\leq C \sum_{j\geq 2\colon c_j > 1/u} u^{-(\tau-2)} c_j\Big(u^{-(\tau-2)}c_j+c_ju^{-(\tau-3)/2}\Big)\\
    &\leq C \sum_{j\geq 2\colon c_j > 1/u} u^{-(\tau-2)}u^{-(\tau-3)/2} c_j^2
    \leq C \sum_{j\geq 2\colon c_j > 1/u} u^{-3(\tau-3)/2}c_j^3= o(1),\nn
    }
where we have used that $\sum_j c_j^3<\infty$ and $\tau>3$ in the last equality. For $c_j\leq 1/u$ and
by \refeq{qj-approx}, we similarly get
    \eqn{
    \lbeq{laplace-A-1-2nd}
    \sum_{j\geq 2\colon c_j \leq 1/u} \shift \Xi_j
    \leq C\shift \sum_{j\geq 2\colon c_j \leq 1/u} \shift u^{-(\tau-3)} c_j^2\Big(u^{-(\tau-3)}c_j^2+c_ju^{-(\tau-3)/2}\Big)
    \leq C \shift\sum_{j\geq 2\colon c_j \leq 1/u} \shift u^{-3(\tau-3)/2} c_j^3=o(1).
    }
This completes the proof that $\sum_{j \geq 2} \Xi_j= o(1)$ and thus of the claim
in \eqref{limit-laplace-A}.
\qed

\paragraph{Step 2: The limit of $\widetilde \expec[\kappa_u]$.}
We proceed by showing that
$\lim_{u \rightarrow \infty} \widetilde \expec[\kappa_u] = \kappa$ with $\kappa>0$ as in \eqref{kappa-def}.
By definition of $\kappa_u$ in \eqref{conv-lin-term}, $q_j(u)$ in \eqref{def-qj} and $\widetilde \prob(T_j \leq u)$ in \eqref{Ti-def-tilt},
    \eqan{
    \lim_{u \rightarrow \infty} \widetilde \expec[\kappa_u]
    &= \lim_{u \rightarrow \infty} \sum_{j\geq 2} q_j(u)\widetilde \prob(T_j \leq u) \lbeq{kappa-mean-conv} \\
    &= \lim_{u \rightarrow \infty} u^{-(\tau-1)} \sum_{j \geq 2} c_ju\e^{-c_ju}\frac{\e^{c_ju}-1-c_ju}{1-\e^{-c_ju}} \frac{\e^{\theta c_ju}(1-\e^{-c_ju})}{\e^{\theta c_ju}(1-\e^{-c_ju})
    +\e^{-c_ju}} \nn\\
    &= \lim_{\Delta \rightarrow 0^+} \Delta \sum_{j \geq 2} x_j^{-\alpha}\e^{-x_j^{-\alpha}} \left[ \e^{x_j^{-\alpha}}-1-x_j^{-\alpha} \right] \frac{\e^{\theta x_j^{-\alpha}}}{\e^{\theta x_j^{-\alpha}}(1-\e^{-x_j^{-\alpha}})
    +\e^{-x_j^{-\alpha}}} \nn\\
    &= \int_0^{\infty} x^{-\alpha} \frac{\e^{\theta x^{-\alpha}}\e^{-x^{-\alpha}}}
    {\e^{\theta x^{-\alpha}}(1-\e^{-x^{-\alpha}})+\e^{-x^{-\alpha}}}
    \big[\e^{x^{-\alpha}}-1-x^{-\alpha} \big]dx,\nn
    }
with $\Delta=u^{-(\tau-1)}$ and $x_j=j \Delta, \ j \geq 2$. Here we used that
the integrand in the last line of \refeq{kappa-mean-conv} is continuous and
integrable over $(0,\infty)$. Set $-x^{-\alpha}=z$ to get the representation
\eqref{kappa-def} for $\kappa$.
\qed

\paragraph{Step 3: Completion of the proof.}
%The convergence of $\widetilde\expec\big[\e^{\im k u^{-(\tau-3)/2} \SS_u}\big]
%\rightarrow \e^{-k^2 I_{\sss V}(1)/2}$
%is already proved in \refeq{conv-FT-Su}.
By Proposition \ref{prop-laplace}, we know that
    \eqn{
    \lbeq{conv-FT-Su}
    \widetilde \expec[\e^{\ii k  u^{-(\tau-3)/2}\SS_u}]\rightarrow \e^{-k^2 I_{\sss V}(1)/2}.
    }
Therefore, Steps 1-2
and \refeq{conv-FT-Su} complete the proof of pointwise convergence in
Lemma \ref{lem-one-dim-A}.
\qed
\medskip

To show that the dominated convergence theorem can be applied, it remains
to show that the integrand in \eqref{laplace-A} has an integrable dominating function:
%......................................................

\begin{Lemma}[Domination by an integrable function]
\label{lem-dom-cvg-A}
    \eqn{
    \lbeq{dom-cvg-A-1st-int}
    \int_{-\infty}^\infty \sup_{u \geq u_0} \left| \widetilde\expec\big[\e^{-a \kappa_u}\e^{\im k u^{-(\tau-3)/2}\SS_u}\big] \right| dk < \infty.
    }
\end{Lemma}

%......................................................

\proof By definition of $\SS_u$ from \eqref{SS-def-ref} and the independence in Lemma \ref{lem-ind-tilt},
    \eqan{
    \lbeq{laplace-A-11}
    \left| \widetilde \expec\Big[ \e^{-a \kappa_u} \e^{\im k u^{-(\tau-3)/2} \SS_u} \Big] \right|^2
    &\leq \prod_{j \geq 2} \left| \widetilde \prob(T_j \leq u) \e^{-a q_j(u)} \e^{\im k u^{-(\tau-3)/2} c_j} + (1-\widetilde \prob(T_j \leq u)) \right|^2 \\
    &= \prod_{j \geq 2} \left[ 1 + \widetilde \prob(T_j \leq u)^2 (\e^{-2a q_j(u)}+1) \right. \nn\\
    &\quad \left. + 2 \widetilde \prob(T_j \leq u) \cos(k u^{-(\tau-3)/2} c_j) \e^{-a q_j(u)} (1-\widetilde \prob(T_j \leq u)) -2\widetilde \prob(T_j \leq u) \right]. \nn
    }
We can rewrite each factor as
    \eqan{
    \lbeq{laplace-A-12}
    &1 -2\widetilde \prob(T_j \leq u) \left\{ \widetilde \prob(T_j \leq u)(1-\e^{-2a q_j(u)})/2+
    (1-\cos(k u^{-(\tau-3)/2} c_j) \e^{-a q_j(u)}) (1-\widetilde \prob(T_j \leq u)) \right\}\\
    &\quad \leq 1 -2\widetilde \prob(T_j \leq u)
    (1-\cos(k u^{-(\tau-3)/2} c_j) \e^{-a q_j(u)}) (1-\widetilde \prob(T_j \leq u)),
    \nn
    }
since $q_j(u)\geq 0$. We then use $\log(1+x) \leq x$ for $x \geq -1$ to obtain
    \eqn{
    \lbeq{laplace-A-2}
    \log\!\left( \left| \widetilde \expec\Big[\e^{-a \kappa_u} \e^{\im k u^{-(\tau-3)/2} \SS_u} \Big] \right|^2 \right) \leq \sum_{j \geq 2} 2 \widetilde \prob(T_j \leq u) \big(\cos(k u^{-(\tau-3)/2} c_j) \e^{-a q_j(u)} -1\big) (1-\widetilde \prob(T_j \leq u)).
    }
The latter equals
    \eqn{
    \lbeq{laplace-A-3}
    \sum_{j \geq 2} 2 \widetilde \prob(T_j \leq u) (\cos(k u^{-(\tau-3)/2} c_j) -1) (1-\widetilde \prob(T_j \leq u)) + e_j(u),
    }
with an overall error term (using that $\sup_j q_j(u)$ is arbitrarily small for $u$ big enough)
    \eqn{
    \sum_{j \geq 2} |e_j(u)|
    \leq C(a) \sum_{j \geq 2} \widetilde \prob(T_j \leq u) q_j(u) (1-\widetilde \prob(T_j \leq u)).
    }
Applying \refeq{lem-approx}, we get
    \eqan{
    \sum_{j \geq 2} |e_j(u)|
    &\leq C(a) \Big\{ \sum_{j \geq 2: c_j>1/u} u^{-(\tau-2)}c_j \e^{-c_ju(1+\theta)} + \sum_{j \geq 2: c_j\leq 1/u} c_ju u^{-(\tau-2)}c_j^2 u \Big\}\leq C(a),
    }
where we have used the bounds
    \eqn{
    \lbeq{lem-help-(a)-(b)}
    \sum_{i: c_i>1/u} u^{-(\tau-2)}c_j \e^{-c_ju(1+\theta)}\leq 1,
    \quad
    u^{-(\tau-4)} \sum_{i: c_i\leq 1/u} c_i^3 \leq C(\tau),
    }
whose proof is straightforward.

Together with \eqref{laplace-A-2} and \eqref{laplace-A-3}, we obtain
    \eqn{
    \log\!\left( \left| \widetilde \expec\Big[ \e^{-a \kappa_u} \e^{\im k u^{-(\tau-3)/2} \SS_u} \Big] \right| \right) \leq \sum_{j \geq 2} \widetilde \prob(T_j \leq u) (\cos(k u^{-(\tau-3)/2} c_j) -1) (1-\widetilde \prob(T_j \leq u)) + C(a).
    }
%The sum over $j$ coincides with the right-hand side of \eqref{density-calc}, which was used
%to determine the asymptotics of the density of $\SS_u$ near zero.
As all summands are nonpositive we obtain together with \eqref{Ti-def-tilt}
    \eqn{
    \log\!\left( \left| \widetilde \expec\Big[ \e^{-a \kappa_u} \e^{\im k u^{-(\tau-3)/2} \SS_u} \Big] \right| \right) \leq C \sum_{j \geq 2: c_j \leq 1/u} c_j u (\cos(k u^{-(\tau-3)/2} c_j) -1) + C(a).
    }
Following the proof of \cite[Proposition 2.5, (6.7)-(6.10)]{AidHofKliLee13b}, we obtain
%the obtained bound in \refeq{rt-xi-bd} yields
    \eqn{
    \lbeq{laplace-A-13}
    \log\!\left( \left| \widetilde \expec\Big[\e^{-a \kappa_u} \e^{\im k u^{-(\tau-3)/2} \SS_u} \Big] \right| \right) \leq C_1-C_2 |k|^{\tau-2}
    }
and integrability of $| \widetilde \expec[ \e^{-a \kappa_u} \e^{\im k u^{-(\tau-3)/2} \SS_u}]|$ against $k$ uniformly in $u$ follows.
\qed

\paragraph{Completion of the proof of  Proposition \ref{prop-weak-cond-conv}(a).}
By the dominated convergence theorem, Lemmas~\ref{lem-one-dim-A} and \ref{lem-dom-cvg-A}
complete the proof of Proposition \ref{prop-weak-cond-conv}(a).
\qed

%%%%%%%%%%%%%%%%%%%%%%%%%%%%%%%%%%%%%%%%%%%%%%%%%%%%%%%

\subsection{Convergence of the process $B_u$}
\label{sec-conv-B}

%%%%%%%%%%%%%%%%%%%%%%%%%%%%%%%%%%%%%%%%%%%%%%%%%%%%%%%

In this section, we investigate the convergence of the
$B_u$ process and prove Proposition \ref{prop-weak-cond-conv}(b).
Since the limit is a \emph{random process}, this part is more involved than the
previous section. We first note that
    \eqn{
    \lbeq{Bu-one-dim}
    B_u(tu^{-(\tau-2)})=\sum_{i\in \JJ(u)} c_iu[\indic{T_i\in (u-tu^{-(\tau-2)},u]}-\widetilde \prob(T_i>u-tu^{-(\tau-2)} \mid T_i\leq u)],
    }
and the processes $(\indic{T_i\in (u-tu^{-(\tau-2)},u]})_{t\geq 0}$ are,
conditionally on $\JJ(u)$, independent. Thus, $(B_u(tu^{-(\tau-2)}))_{t\geq 0}$
is, conditionally on $\JJ(u)$, a sum of (conditionally) independent
processes having zero mean. We make crucial use of this observation, as well
as the technique in Lemma \ref{lem-cond-expec-SSu}, to compute expectations
of various functionals of the process $(B_u(tu^{-(\tau-2)}))_{t\geq 0}$.

In order to prove the stated convergence in distribution, we follow
the usual path of first proving weak convergence of the one-dimensional marginals,
followed by the weak convergence of all finite-dimensional distributions,
and complete the proof by showing tightness. We now discuss each of these steps
in more detail.

%%%%%%%%%%%%%%%%%%%%%%%%%%%%%%%%%%%%%%%%%%%%%%%%%%%%%%%

\subsubsection{Convergence of the one-dimensional marginal of $B_u$}

%%%%%%%%%%%%%%%%%%%%%%%%%%%%%%%%%%%%%%%%%%%%%%%%%%%%%%%

We start by computing the one-dimensional marginal of $B_u(tu^{-(\tau-2)})$
(recall \eqref{Bu-one-dim}) and show that
it is consistent with the claimed L\'evy process limit.
We achieve this by computing the Laplace transform
    \eqn{
    \lbeq{psi-a-def}
    \psi_{u,v}(a)=\widetilde\expec_v\big[\e^{-a B_u(tu^{-(\tau-2)})}\big],
    }
and proving that it converges to the Laplace transform
of the claimed L\'evy process limit at time $t$.
The main result in this section is the following proposition:

\begin{Proposition}[One-time marginal of $B_u(tu^{-(\tau-2)})$]
\label{prop-one-dim}
There exists a measure $\Pi$ such that, for every $v, a>0$ fixed and as
$u\rightarrow \infty$,
    \eqn{
    \psi_{v,u}(a)
    \rightarrow \e^{t\int_0^{\infty} (\e^{-a z}-1+a z)\Pi(dz)},
    }
which is the Laplace transform of a L\'evy process $(L_s)_{s\geq 0}$ with non-negative jumps
and characteristic measure $\Pi$
   \eqn{
   \lbeq{charact-meas-def-again}
   \Pi(dz) \equiv \e^z \frac{\e^{-\theta z}}{\e^{-\theta z}(1-\e^z)+\e^z} (\tau-1) (-z)^{-(\tau-1)} dz.
   }
Therefore, the one-dimensional marginals of the process $(B_u(su^{-(\tau-2)}))_{t\geq 0}$ converge
to those of $(L_s)_{s\geq 0}$.
\end{Proposition}

The remainder of this section is devoted to the proof
of Proposition \ref{prop-one-dim}.
As for $A_u$, we use Lemma~\ref{lem-cond-expec-SSu} and a change of
variables to rewrite
    \eqan{
    \lbeq{laplace-B-complete}
    \psi_{v,u}(a)&\equiv \widetilde\expec_v\Big[\e^{-a B_u(tu^{-(\tau-2)})}\Big]
    =\frac{1}{u^{(\tau-3)/2} \widetilde f_{\SS_u}(v/u)} \int_{-\infty}^{+\infty}
    \e^{-\im k v u^{-(\tau-1)/2}} \widetilde \expec\big[\e^{-a B_u(tu^{-(\tau-2)})}\e^{\im k u^{-(\tau-3)/2} \SS_u}]\frac{dk}{2\pi}\nn\\
    &= \frac{1}{u^{(\tau-3)/2} \widetilde f_{\SS_u}(v/u)} \int_{-\infty}^{+\infty}
    \e^{-\im k v u^{-(\tau-1)/2}} \widetilde \expec\Big[\psi_{\sss \JJ}(a) \e^{\im k u^{-(\tau-3)/2} \SS_u} \Big] \frac{dk}{2\pi},
    }
where
    \eqan{
    \lbeq{laplace-B}
    \psi_{\sss \JJ}(a)
    &\equiv \widetilde\expec\Big[\e^{-a B_u(tu^{-(\tau-2)})}\mid \JJ(u)\Big]\\
    &=\prod_{j\in \JJ(u)} \e^{a c_ju \widetilde \prob(T_j>u-tu^{-(\tau-2)} \mid T_j\leq u)}
    \Big(1+(\e^{-a c_ju}-1) \widetilde \prob(T_j>u-tu^{-(\tau-2)} \mid T_j\leq u)\Big)\nn\\
    &=\prod_{j\in \JJ(u)} \e^{a c_ju p_{j,t}^u}
    \Big(1+(\e^{-a c_ju}-1)p_{j,t}^u\Big),\nn
    }
and where we abbreviate
    \eqn{
    \lbeq{def-p-jtu}
    p_{j,t}^u
    = \widetilde \prob(T_j>u-tu^{-(\tau-2)} \mid T_j\leq u) = \frac{\e^{c_jtu^{-(\tau-2)}}-1}{1-\e^{-c_ju}} \e^{-c_ju},
    }
by \eqref{cond-distr-T}. %Naturally, $\psi_{u,v}(a)=\widetilde\expec_v[\psi_{\sss \JJ}(a)]$.
We again wish to use dominated convergence on the integral in \refeq{laplace-B-complete}.

We proceed along the lines of the proof of the convergence of the mean process $A_u$.
Basically, in the proof below, we replace $-a \kappa_u$ in \eqref{laplace-A} (recall the definition of $\kappa_u$ and $q_j(u)$ from \eqref{conv-lin-term} and \eqref{def-qj}) by $\sum_{j\in \JJ(u)} r_{j,t}^u$, where
we define
    \eqn{
    \lbeq{def-rj}
    r_{j,t}^u
    \equiv (\e^{-a c_ju}-1+a c_ju) p_{j,t}^u
    = (\e^{-a c_ju}-1+a c_ju) \frac{\e^{c_jtu^{-(\tau-2)}}-1}{1-\e^{-c_ju}} \e^{-c_ju}.
    }
In what follows, we frequently make use of the bounds
    \eqn{
    \lbeq{cond-distr-T-B}
    p_{j,t}^u
    \leq C t u^{-(\tau-1)} (c_j u \e^{-c_ju}\wedge 1)
    \leq C t u^{-(\tau-1)},
    }
and
    \eqan{
    \lbeq{approx-rj}
    r_{j,t}^u
    \leq C(a,T) c_j u^{-(\tau-2)}(1\wedge c_j u).
    }
We again start by proving pointwise convergence:

\begin{Lemma}[Pointwise convergence revisited]
\label{lem-one-dim-B}
For $a \geq 0$ arbitrary, $v=o(u^{(\tau-1)/2})$,
    \eqn{
    \lbeq{limit-laplace-B}
    \e^{-\im k vu^{-(\tau-1)/2}}\widetilde\expec\big[\psi_{\sss \JJ}(a)\e^{\im k u^{-(\tau-3)/2} \SS_u}\big]
    = \e^{t\int_{-\infty}^0 (\e^{a z}-1-a z) \Pi(dz)}\e^{-I_{\sss V}(1) k^2/2}+ o(1).
    }
\end{Lemma}

\proof The first factor on the left-hand side of \refeq{limit-laplace-B}
converges to 1. We identify the limit of the expectation in the following steps that mimic
the pointwise convergence proof in Lemma \ref{lem-one-dim-A}. It will be convenient to
split the asymptotic factorization in Step 1 of that proof into
two parts, denoted by Steps~1(a) and 1(b). We start by showing that
we can simplify $\psi_{\sss \JJ}(a)$:

\paragraph{Step 1(a): Simplification of $\psi_{\sss \JJ}(a)$.}
As a first step towards the identification of the pointwise limit,
we show that we can simplify the expectation in \refeq{limit-laplace-B}
as follows:
    \eqn{
    \lbeq{laplace-B-0}
    \widetilde \expec\!\left[\left|\psi_{\sss \JJ}(a)
    - \e^{\sum_{j\in \JJ(u)} r_{j,t}^u} \right|\right] = o(1).
    }

To prove \refeq{laplace-B-0}, we denote the difference in \refeq{laplace-B-0} by
    \eqn{
    \lbeq{Eu-t-def}
    E_u(t)
    = \prod_{j\in \JJ(u)} \e^{a c_ju p_{j,t}^u} \left| \prod_{j\in \JJ(u)}
    \Big(1+(\e^{-a c_ju}-1) p_{j,t}^u \Big) - \prod_{j\in \JJ(u)} \e^{(\e^{-a c_ju}-1) p_{j,t}^u} \right|,
    }
so that
	\eqn{
	\lbeq{laplace-B-0-bd-Eu-t}
	\Big|\widetilde \expec\!\left[\left|\psi_{\sss \JJ}(a)
    - \e^{\sum_{j\in \JJ(u)} r_{j,t}^u} \right|\right]\Big|
	\leq \widetilde \expec[E_u(t)].
	}
	
Using the first line of \eqref{prod-diff-error} and applying the error bound $|\e^x-(1+x)| \leq |x|^2$ for $|x| \leq 1$ to the differences $|a_j-b_j|$, the error of the approximation can be bounded by
    \eqn{
    E_u(t)
    \leq C \prod_{j\in \JJ(u)} \e^{a c_ju p_{j,t}^u} \sum_{j\in \JJ(u)} \prod_{i\in \JJ(u), i<j} \Big(1+(\e^{-a c_iu}-1) p_{i,t}^u \Big) \left[ (\e^{-a c_ju}-1) p_{j,t}^u \right]^2 \prod_{i\in \JJ(u), i>j} \e^{(\e^{-a c_iu}-1) p_{i,t}^u}.
    }
Next use that $1-x \leq \e^{-x}$ for $x \geq 0$ to obtain as a further bound to the above
    \eqn{
    C \sum_{j\in \JJ(u)} \e^{a c_ju p_{j,t}^u} \prod_{i\in \JJ(u) \backslash \{j\}} \e^{(\e^{-a c_iu}-1+a c_iu) p_{i,t}^u} \left[ (\e^{-a c_ju}-1) p_{j,t}^u \right]^2.
    }
For $t \leq T$ with $T>0$ fixed, we further have by \refeq{cond-distr-T-B} that $\e^{a c_ju p_{j,t}^u} \leq C(a,T)$. Together with $\e^{-x}-1+x\geq 0$ for $x \geq 0$, we obtain
    \eqn{
    E_u(t)
    \leq C(a,T) \prod_{i\in \JJ(u)} \e^{(\e^{-a c_iu}-1+a c_iu) p_{i,t}^u} \sum_{j\in \JJ(u)} \left[ (\e^{-a c_ju}-1) p_{j,t}^u \right]^2.
    }
The bound $\e^{-x}-1+x\leq x^2/2, \,\forall x \geq 0$ yields
    \eqn{
    \lbeq{laplace-B-4}
    E_u(t)
    \leq C(a,T) \e^{ \frac{a^2}{2} \sum_{i\in \JJ(u)} (c_iu)^2 p_{i,t}^u} \sum_{j\in \JJ(u)} \left[ (\e^{-a c_ju}-1) p_{j,t}^u \right]^2.
    }

We first bound the sum in \eqref{laplace-B-4}. With $1-\e^{-x} \leq x$ for $x \geq 0$ and
by \refeq{cond-distr-T-B} we obtain
    \eqn{
    \sum_{j\in \JJ(u)} \left[ (\e^{-a c_ju}-1) p_{j,t}^u \right]^2
    \leq \sum_{j\in \JJ(u)} (a c_ju)^2 \left( p_{j,t}^u \right)^2
    \leq C(a,T) u^{-2(\tau-3)-1} \Big\{ C + \sum_{j\in \JJ(u)\colon c_j\leq 1/u} c_j^2 u^{-1} \Big\}.
    }
This yields as an upper bound for \eqref{laplace-B-0} (recall \eqref{laplace-B-0-bd-Eu-t}),
    \eqan{
    \widetilde \expec[E_u(t)]
    &\leq C(a,T) u^{-2(\tau-3)-1} \Big\{ \widetilde \expec\!\left[ \e^{ \frac{a^2}{2} \sum_{i\in \JJ(u)} (c_iu)^2 p_{i,t}^u} \right]\\
    &\quad + \sum_{j \colon c_j\leq 1/u} c_j^2 u^{-1} \widetilde \expec\!\Big[\indic{j \in \JJ(u)} \e^{ \frac{a^2}{2} (c_ju)^2 p_{j,t}^u } \Big] \widetilde \expec\!\Big[ \e^{ \frac{a^2}{2} \sum_{i\in \JJ(u) \backslash \{j\}} (c_iu)^2 p_{i,t}^u} \Big] \Big\} \nn\\
    &\leq C(a,T) u^{-2(\tau-3)-1} \{ 1 + \sum_{j \colon c_j\leq 1/u} c_j^3 \} \widetilde \expec\!\left[ \e^{ \frac{a^2}{2} \sum_{i\in \JJ(u)} (c_iu)^2 p_{i,t}^u} \right], \nn
    }
where we have used \refeq{lem-approx} in the last line.

The claim \eqref{laplace-B-0} follows once we show that $\widetilde \expec\!\Big[ \e^{ \frac{a^2}{2} \sum_{i\in \JJ(u)} (c_iu)^2 p_{i,t}^u}\Big]$ is bounded. To prove this, consider first the sum over $c_i>1/u$ only. By
\refeq{cond-distr-T-B} and \refeq{lem-help-(a)-(b)},
    \eqn{
    \sum_{i\in \JJ(u): c_i>1/u} (c_iu)^2 p_{i,t}^u
    \leq C \sum_{i\in \JJ(u): c_i>1/u} (c_iu)^2 c_i tu^{-(\tau-2)} \e^{-c_iu}
    \leq C \sum_{i: c_i>1/u} tu^{-(\tau-1)}
    \leq C(T).
    }
Using \refeq{cond-distr-T-B} once more, it remains to show the boundedness of
    \eqan{
    \widetilde \expec\!\left[ \e^{ \frac{a^2}{2} \sum_{i\in \JJ(u): c_i\leq 1/u} (c_iu)^2 p_{i,t}^u} \right]
    &\leq \widetilde \expec\!\left[ \e^{ C(a) \sum_{i\in \JJ(u): c_i\leq 1/u} (c_iu)^2 t u^{-(\tau-1)} } \right] \\
    &= \prod_{c_i\leq 1/u} \left( \widetilde \prob(T_i\leq u) \e^{ C(a) (c_iu)^2 t u^{-(\tau-1)} } + (1-\widetilde \prob(T_i\leq u)) \right), \nn
    }
which is equivalent to bounding
    \eqn{
    \lbeq{laplace-B-6}
    \sum_{c_i\leq 1/u} \log\!\left( \widetilde \prob(T_i\leq u) \e^{ C(a) (c_iu)^2 t u^{-(\tau-1)} } + (1-\widetilde \prob(T_i\leq u)) \right)
    \leq \sum_{c_i\leq 1/u} \widetilde \prob(T_i\leq u) \left( \e^{ C(a) (c_iu)^2 t u^{-(\tau-1)} } - 1 \right)
    }
appropriately. Here we used that $\log(1+x) \leq x$ for $x \geq 0$. Next bound $\widetilde \prob(T_i\leq u) \leq C c_iu$ in the above to obtain that for $c_i\leq 1/u$ we have $C(a) (c_iu)^2 t u^{-(\tau-1)} \leq C(a,T) u^{-(\tau-1)} \leq \log(2)$ for $u$ big enough. Hence we can use that $\e^x-1 \leq 2x$ for $0 \leq x \leq \log(2)$ and thus get as a further upper bound to \eqref{laplace-B-6}
    \eqn{
    C(a,T) \sum_{c_i\leq 1/u} c_iu (c_iu)^2 u^{-(\tau-1)}
    \leq C(a,T).
    }
The last inequality follows from \refeq{lem-help-(a)-(b)}. This completes the proof of \refeq{laplace-B-0}.
\qed
\medskip

\paragraph{Step 1(b): Asymptotic factorization.} We next show that
    \eqn{
    \lbeq{laplace-B-7}
    \widetilde\expec\big[\e^{\sum_{j\in \JJ(u)} r_{j,t}^u}\e^{\im k u^{-(\tau-3)/2} \SS_u}\big]
    = \e^{\widetilde\expec[\sum_{j\in \JJ(u)} r_{j,t}^u]} \widetilde\expec\big[\e^{\im k u^{-(\tau-3)/2} \SS_u}\big] + o(1).
    }
To prove \refeq{laplace-B-7}, we note that, by the definition of $r_{j,t}^u$ in \eqref{def-rj},
    \eqan{
    \bar{E}_u(t)
    &\equiv \left| \widetilde\expec\big[\e^{\sum_{j\in \JJ(u)} r_{j,t}^u}\e^{\im k u^{-(\tau-3)/2} \SS_u}\big]
    - \e^{\widetilde\expec[\sum_{j\in \JJ(u)} r_{j,t}^u]} \widetilde\expec\big[\e^{\im k u^{-(\tau-3)/2} \SS_u}\big] \right| \\
    &= \left| \prod_{j\geq 2} \widetilde\expec\big[ \e^{\II_j(u) r_{j,t}^u + \im k u^{-(\tau-3)/2} c_j \II_j(u)} \big] - \prod_{j\geq 2} \e^{r_{j,t}^u \widetilde \prob(T_j \leq u)} \widetilde\expec\big[ \e^{\im k u^{-(\tau-3)/2} c_j \II_j(u)} \big] \right|. \nn
    }
As in the calculations of the Laplace transform of $A_u$ in \eqref{prod-diff-laplace-A}, we now apply \eqref{prod-diff-error}. Note that here we cannot apply the second bound of \eqref{prod-diff-error} as $\sup_j (|a_j| \vee |b_j|)$ is not bounded by $1$ (recall that $r_{j,t}^u \geq 0$). Instead, we get
    \eqan{
    \bar{E}_u(t)
    \lbeq{laplace-B-9}
    &\leq \sum_{j \geq 2} \prod_{2\leq j_1\leq j-1} \left| \widetilde\expec\big[ \e^{\II_{j_1}(u) r_{j_1,t}^u + \im k u^{-(\tau-3)/2} c_{j_1} \II_{j_1}(u)} \big] \right| \\
    &\quad \times \left| \widetilde\expec\big[ \e^{\II_j(u) r_{j,t}^u + \im k u^{-(\tau-3)/2} c_j \II_j(u)} \big] - \e^{r_{j,t}^u\widetilde \prob(T_j \leq u)} \widetilde\expec\big[ \e^{\im k u^{-(\tau-3)/2} c_j \II_j(u)} \big] \right| \nn\\
    &\quad \times \prod_{j_2\geq j+1} \left| \e^{r_{j_2,t}^u\widetilde \prob(T_{j_2} \leq u)} \widetilde\expec\big[ \e^{\im k u^{-(\tau-3)/2} c_{j_2} \II_{j_2}(u)} \big] \right|. \nn
    }
We proceed to prove that the first and the third product are bounded by constants. Indeed, we can bound the third product using \refeq{lem-approx} by
    \eqn{
    \prod_{j_2\geq j+1} \e^{r_{j_2,t}^u\widetilde \prob(T_{j_2} \leq u)}
    \leq \e^{C\sum_{j \geq 2} r_{j,t}^u (c_j u\wedge 1)},
    }
where, by \refeq{approx-rj} and \refeq{lem-help-(a)-(b)},
    \eqn{
    \lbeq{laplace-B-8}
    \sum_{j \geq 2} r_{j,t}^u (c_j u\wedge 1)\leq C.
    }
For the first product in \eqref{laplace-B-9}, we obtain as an upper bound
    \eqan{
    \prod_{j\geq 2} \widetilde\expec\big[ \e^{\II_j(u) r_{j,t}^u} \big]
    &= \prod_{j\geq 2} \left( \widetilde \prob(T_j \leq u) \e^{r_{j,t}^u} + (1-\widetilde \prob(T_j \leq u)) \right)
    \lbeq{laplace-B-11}\\
    &= \prod_{j\geq 2} \left(1+\widetilde \prob(T_j \leq u)(\e^{r_{j,t}^u}-1)\right)
    \leq \e^{\sum_{j \geq 2} \widetilde \prob(T_j \leq u) (\e^{r_{j,t}^u}-1)}.
    \nn
    }
As $r_{j,t}^u$ is uniformly bounded for $u$ big enough the above is again
bounded by \eqref{laplace-B-8}.

Hence, it suffices to bound the middle part of \eqref{laplace-B-9}, that is, it remains to show that
    \eqn{
    \lbeq{laplace-B-10}
    \sum_{j \geq 2} \left| \widetilde\expec\big[ \e^{\II_j(u) r_{j,t}^u + \im k u^{-(\tau-3)/2} c_j \II_j(u)} \big] - \e^{r_{j,t}^u\widetilde \prob(T_j \leq u)} \widetilde\expec\big[ \e^{\im k u^{-(\tau-3)/2} c_j \II_j(u)} \big] \right|
    =\Delta_j(r_{j,t}^u)= o(1),
    }
where we recall the definition of $\Delta_j(q)$ in \refeq{Delta-j-q-def}.
By \refeq{diff-laplace-A},
    \eqn{
    \Delta_j(r_{j,t}^u)\leq C r_{j,t}^u \e^{r_{j,t}^u}(r_{j,t}^u+|k| u^{-(\tau-3)/2} c_j)
    }
for $u=u(k)$ big enough. The bound on $r_{j,t}^u$ in \refeq{approx-rj}
is equal to $C(a,T)$ times the bounds on $q_j(u)$ in \refeq{qj-approx}. The remaining calculations for $A_u$
in \eqref{diff-laplace-A}-\eqref{laplace-A-1-2nd} therefore directly carry over, so that \eqref{laplace-B-10} follows.
\qed
\medskip

\paragraph{Step 2: The limit of $\expec\![\sum_{j\in \JJ(u)} r_{j,t}^u ]$.}
In this step, we identify the limit of $\expec\![\sum_{j\in \JJ(u)} r_{j,t}^u ]$.
For this, we use that by definition of $r_{j,t}^u$ in \eqref{def-rj},
that of $p_{j,t}^u$ in \eqref{def-p-jtu}, and \eqref{Ti-def-tilt} with $t=u$,
    \eqan{
%    \lbeq{limit-B-guess}
    \widetilde \expec\Big[ \sum_{j\in \JJ(u)} r_{j,t}^u \Big]
    &=\sum_{j \geq 2} (\e^{-a c_ju}-1+a c_ju) \frac{\e^{c_jtu^{-(\tau-2)}}-1}{1-\e^{-c_ju}} \e^{-c_ju} \frac{\e^{\theta c_ju}(1-\e^{-c_ju})}{\e^{\theta c_ju}(1-\e^{-c_ju})
    +\e^{-c_ju}} \\
    &= tu^{-(\tau-1)} \sum_{j \geq 2} (\e^{-a c_ju}-1+a c_ju) \frac{\e^{c_ju tu^{-(\tau-1)}}-1}{t u^{-(\tau-1)}} \e^{-c_ju} \frac{\e^{\theta c_ju}}{\e^{\theta c_ju}(1-\e^{-c_ju})
    +\e^{-c_ju}} \nn\\
    &\stackrel{u \rightarrow \infty}{\rightarrow} t\int_0^\infty (\e^{-a x^{-\alpha}}-1+a x^{-\alpha}) x^{-\alpha} \e^{-x^{-\alpha}} \frac{\e^{\theta x^{-\alpha}}}{\e^{\theta x^{-\alpha}}(1-\e^{-x^{-\alpha}})+\e^{-x^{-\alpha}}} dx. \nn
    }
The convergence of the sum to the integral follows as in \refeq{kappa-mean-conv}.
Next set $-x^{-\alpha}=z$ to get
    \eqn{
    \lbeq{limit-B-guess-2}
    \lim_{u \rightarrow \infty} \frac{1}{t} \widetilde \expec\Big[ \sum_{j\in \JJ(u)} r_{j,t}^u\Big]
    = \int_{-\infty}^0 (\e^{a z}-1-a z) \Pi(dz),
    }
with $\Pi(dz)$ as in \eqref{charact-meas-def-again}, respectively, \eqref{charact-meas-def}.
For $\Pi$ to be the L\'evy measure of a real-valued L\'evy process with no positive jumps
as in \cite[Section V.1]{Bert96}, by the L\'evy-Khintchine formula in
\cite[Section 0.2 and Theorem 1 in Section I.1]{Bert96}, we have to check
that $\Pi$ is a measure on $(-\infty,0)$ that satisfies
$\int \Pi(dz) (1 \wedge z^2) < \infty$. Indeed, close to $0$, $z^2 \Pi(dz)$
behaves like $(\tau-1) z^{-(\tau-3)} dz$, which is integrable at $0$ and for
$z \rightarrow \infty$, $\Pi(dz)$ behaves like $\e^{-z} (\tau-1) z^{-(\tau-1)} dz$,
whose integral is finite for all $n \in \mathbb{N}$.
\qed

\paragraph{Step 3: Completion of the proof.}
The convergence of $\widetilde\expec\big[\e^{\im k u^{-(\tau-3)/2} \SS_u}\big]\rightarrow \e^{-k^2 I_{\sss V}(1)/2}$
is already proved in \refeq{conv-FT-Su}. Therefore, Steps 1(a)-1(b) and 2,
together with \refeq{conv-FT-Su}, complete the proof of pointwise convergence in
Lemma \ref{lem-one-dim-B}.
\qed
\medskip

To show that the dominated convergence theorem can be applied, it again remains
to show that the integrand has an integrable dominating function:

\begin{Lemma}[Domination by an integrable function]
\label{lem-dom-cvg-B}
    \eqn{
    \lbeq{dom-cvg-B-1st-int}
    \int_{-\infty}^\infty \sup_{u \geq u_0} \left| \widetilde\expec\big[\psi_{\sss \JJ}(a)\e^{\im k u^{-(\tau-3)/2}\SS_u}\big] \right| dk < \infty.
    }
\end{Lemma}

\proof This follows in a similar way as in the proof of
Lemma \ref{lem-dom-cvg-A}. We compute
    \eqan{
    \Big|\widetilde\expec\big[\psi_{\sss \JJ}(a)\e^{\im k u^{-(\tau-3)/2}\SS_u}\big]\Big|^2
    &=\Big|\widetilde\expec\big[\e^{\im k u^{-(\tau-3)/2}\SS_u} \prod_{j\in \JJ(u)} \e^{a c_ju p_{j,t}^u}
    \big(1+(\e^{-a c_ju}-1)p_{j,t}^u\big)\big]\Big|^2\\
    &=\prod_{j\geq 2}\Big|1-\widetilde \prob(T_j \leq u)
    +\e^{\im k u^{-(\tau-3)/2}c_j}\widetilde \prob(T_j \leq u)\e^{a c_ju p_{j,t}^u}
    \big(1+(\e^{-a c_ju}-1)p_{j,t}^u\big)\Big|^2.\nn
    }
This is identical to the bound appearing in \eqref{laplace-A-11},
apart from the fact that the term $e^{-aq_j(u)}$ in \eqref{laplace-A-11}
is replaced with $b_{j,t}(u)=\e^{a c_ju p_{j,t}^u}\big(1+(\e^{-a c_ju}-1)p_{j,t}^u\big)$
in the above. Proceeding as in
\eqref{laplace-A-11} to \eqref{laplace-A-3}, we finally obtain
    \eqan{
    \lbeq{B-sum-end}
    &\log\left| \widetilde \expec\Big[\psi_{\sss \JJ}(a)\e^{\im k u^{-(\tau-3)/2} \SS_u} \Big]\right|^2\\
    &\leq \sum_{j \geq 2} 2 \widetilde \prob(T_j \leq u)
    \Big[\widetilde \prob(T_j \leq u)((b_{j,t}(u))^2-1)/2+
    \big(\cos(k u^{-(\tau-3)/2} c_j)b_{j,t}(u)-1\big)(1-\widetilde \prob(T_j \leq u)),\nn
    }
where the additional first term in comparison to \eqref{laplace-A-2} arises because
$b_{j,t}(u)\leq 1$ no longer holds. Indeed, since $\e^{xp}(1+(\e^{-x}-1)p)\geq 1$ for
$x\geq 0$ and $p\in [0,1]$, we have that $b_{j,t}(u)\geq 1$. Further,
    \eqn{
    \lbeq{bjt-ub}
    b_{j,t}(u)=\e^{a c_ju p_{j,t}^u}\big(1+(\e^{-a c_ju}-1)p_{j,t}^u\big)
    \leq \e^{(\e^{-a c_ju}-1+a c_ju) p_{j,t}^u}=\e^{r_{j,t}(u)}.
    }

The first part of the sum in \refeq{B-sum-end} can, by \refeq{bjt-ub}
and since $\e^x-1 \leq 2x$ for $0 \leq x \leq \log(2)$, be bounded by
    \eqn{
    \sum_{j \geq 2} \widetilde \prob(T_j \leq u)^2(\e^{2 r_{j,t}^u}-1) \\
    \leq C(a,T)\sum_{j \geq 2}(1\wedge c_ju)^2 r_{j,t}^u.
    }
Now we can apply \refeq{approx-rj} and \refeq{lem-help-(a)-(b)} to get as a further bound
    \eqn{
    C(a,T) \Big\{ \sum_{j \geq 2: c_j>1/u} u^{-(\tau-1)} + \sum_{j \geq 2: c_j\leq 1/u} c_j^3 u^{-(\tau-4)} \Big\}
    \leq C.
    }
For the second part of the sum \eqref{B-sum-end}, we proceed as in \eqref{laplace-A-2}-\eqref{laplace-A-13} to split it as
    \eqan{
    &-\sum_{j \geq 2} 2 \widetilde \prob(T_j \leq u)
    \widetilde \prob(T_j \leq u)\big(1-\cos(k u^{-(\tau-3)/2} c_j)\big)(1-\widetilde \prob(T_j \leq u))\\
    &\quad +\sum_{j \geq 2} \widetilde \prob(T_j \leq u) \cos(k u^{-(\tau-3)/2} c_j)(b_{j,t}(u)-1)(1-\widetilde \prob(T_j \leq u)).\nn
    }
By a second order Taylor expansion and the fact that $r_{j,t}(u)$
is bounded, there exists a constant $C$ such that $b_{j,t}(u)-1\leq C r_{j,t}(u).$
Now we can proceed as in \eqref{laplace-A-2}-\eqref{laplace-A-13}, where we again take advantage of
being able to dominate the bounds on $r_{j,t}^u$ in \refeq{approx-rj} by
the bounds on $q_j(u)$ in \refeq{qj-approx}. Integrability of
$| \widetilde \expec[\psi_{\sss \JJ}(a)
\e^{\im k u^{-(\tau-3)/2} \SS_u} ] |$ against $k$ follows.
\qed

\paragraph{Proof of Proposition \ref{prop-one-dim}.}
The claim follows from Lemmas \ref{lem-one-dim-B}, \ref{lem-dom-cvg-B}
and the dominated convergence theorem.
\qed

%%%%%%%%%%%%%%%%%%%%%%%%%%%%%%%%%%%%%%%%%%%%%%%%%%%%%%%

\subsubsection{Convergence of the finite-dimensional distributions of $B_u$}
%%%%%%%%%%%%%%%%%%%%%%%%%%%%%%%%%%%%%%%%%%%%%%%%%%%%%%%

In this section, the convergence of the
one-dimensional marginals of the process $(B_u(tu^{-(\tau-2)}))_{t\geq 0}$
gets extended to convergence of its finite-dimensional
distributions. In the same way as above, it can be shown that, for
$0<t_1\cdots<t_n$, the increments $(B_u(t_iu^{-(\tau-2)})-B_u(t_{i-1}u^{-(\tau-2)}))_{i=1}^n$
(where, by convention, $t_0=0$) converge in distribution, under $\widetilde \prob_v$,
to \emph{independent} L\'evy random variables with the correct distribution.

In what follows, we only outline some minor changes in the proof. Instead of \eqref{laplace-B}, we fix $n \in \mathbb{N}$, $\vec{a} \in (\mathbb{R}^+)^n$ and $0 = t_0 < t_1 < \cdots < t_n \leq T$ and consider
    \eqan{
    \psi_{\sss {\mathcal J}}(\vec{a})
    &\equiv \widetilde\expec\Big[\e^{-\sum_{k=1}^n a_k \left( B_u(t_k u^{-(\tau-2)})-B_u(t_{k-1} u^{-(\tau-2)}) \right)}\mid {\mathcal J}(u)\Big]\\
    &=\prod_{j\in {\mathcal J}(u)} \e^{c_ju \sum_{k=1}^n a_k p_{j,t_{k-1},t_k}^u}
    \Big(1+\sum_{k=1}^n (\e^{-a_k c_ju}-1) p_{j,t_{k-1},t_k}^u \Big) \nn
    }
with (the two-point analogue to \refeq{cond-distr-T-B})
    \eqn{
    \lbeq{p-two-times}
    p_{j,s,t}^u
    \equiv \widetilde \prob(T_j\in (u-tu^{-(\tau-2)},u-su^{-(\tau-2)}] \mid T_j\leq u)
    = \e^{-c_ju} \frac{\e^{c_jtu^{-(\tau-2)}}-\e^{c_jsu^{-(\tau-2)}}}{1-\e^{-c_ju}}
    }
for $0 \leq s \leq t \leq T$, using \eqref{cond-distr-T}.
Then, clearly, \refeq{cond-distr-T-B} is replaced with
    \eqn{
    \lbeq{p_two_times_bounds}
    p_{j,s,t}^u
    \leq C (t-s) u^{-(\tau-1)} (\e^{-c_ju} c_j u \wedge 1).
    }
We follow Steps 1(a)-(b) to Step 3 in the proof of convergence of the one-time marginal.

Similarly to Step 1(a), one can show that
    \eqn{
    \widetilde \expec\!\Big[ \Big| \psi_{\mathcal J}(\vec{a})
    - \prod_{j\in {\mathcal J}(u)} \e^{\sum_{k=1}^n (\e^{-a_k c_ju}-1+a_k c_ju) p_{j,t_{k-1},t_k}^u} \Big| \Big] = o(1).
    }
We then continue to reason as from \eqref{laplace-B-complete}
onwards, where $r_{j,t}^u$ in \eqref{def-rj} gets replaced by
    \eqn{
    r_{j,\vec{t}}^u \equiv \sum_{k=1}^n (\e^{-a_k c_ju}-1+a_k c_ju) p_{j,t_{k-1},t_k}^u.
    }
The remaining calculations are analogous to the one-dimensional case.
The asymptotic factorization in Step 1(b) is replaced with
    \eqn{
    \widetilde\expec\big[\e^{\sum_{j\in {\mathcal J}(u)} r_{j,\vec{t}}^u}\e^{\im k u^{-(\tau-3)/2} \SS_u}\big]
    = \e^{\widetilde\expec[\sum_{j\in {\mathcal J}(u)} r_{j,\vec{t}}^u]} \widetilde\expec\big[\e^{\im k u^{-(\tau-3)/2} \SS_u}\big] + o(1)
    }
and we calculate the limit of $\widetilde \expec[ \sum_{j\in {\mathcal J}(u)} r_{j,\vec{t}}^u]$
in a similar way as in Step 2 in the previous subsection as
    \eqan{
    \widetilde \expec\Big[ \sum_{j\in {\mathcal J}(u)} r_{j,\vec{t}}^u\Big]
    &= \sum_{j \geq 2} \sum_{k=1}^n (\e^{-a_k c_ju}-1+a_k c_ju) \frac{\e^{c_jt_ku^{-(\tau-2)}}-\e^{c_jt_{k-1}u^{-(\tau-2)}}}{1-\e^{-c_ju}} \e^{-c_ju} \widetilde\prob( j \in \mathcal{J}(u) ) \\
    &\stackrel{u \rightarrow \infty}{\rightarrow} \int_0^\infty \sum_{k=1}^n (\e^{-a_k x^{-\alpha}}-1+a_k x^{-\alpha}) x^{-\alpha} (t_k-t_{k-1}) \e^{-x^{-\alpha}} \frac{\e^{\theta x^{-\alpha}}}{\e^{\theta x^{-\alpha}}(1-\e^{-x^{-\alpha}})+\e^{-x^{-\alpha}}} dx \nn\\
    &= \int_{-\infty}^0 \sum_{k=1}^n (\e^{a_k z}-1-a_k z) (t_k-t_{k-1}) \Pi(dz). \nn
    }
Finally, we note that
    \eqn{
    \e^{\widetilde\expec[\sum_{j\in {\mathcal J}(u)} r_{j,\vec{t}}^u]}
    = \exp\!\left[ \int_{-\infty}^0 \sum_{k=1}^n (\e^{a_k z}-1-a_k z) (t_k-t_{k-1}) \Pi(dz) \right]
    = \expec\!\left[ \e^{\sum_{k=1}^n a_k (-(L_{t_k}-L_{t_{k-1}})) } \right],
    }
where we have used that, by definition, L\'evy processes have independent
stationary increments. This completes the convergence of the finite-dimensional
distributions of  $(B_u(tu^{-(\tau-2)}))_{t\geq 0}$.
\qed

%%%%%%%%%%%%%%%%%%%%%%%%%%%%%%%%%%%%%%%%%%%%%%%%%%%%%%%

\subsubsection{Tightness of $B_u$}

%%%%%%%%%%%%%%%%%%%%%%%%%%%%%%%%%%%%%%%%%%%%%%%%%%%%%%%

We next turn to \emph{tightness} of the process $(B_u(tu^{-(\tau-2)}))_{t\geq 0}$.
For this, we use the following tightness criterion:

%......................................................

\begin{Proposition}[Tightness criterion \protect{\cite[Theorem 15.6 and the comment following it]{Bill99}}]
The sequence $\{X_n\}$ is tight in $D([0,T],\R^d)$ if the limiting process
$X$ has a.s.\ no discontinuity at $t=T$ and there exist constants $C>0$, $r>0$ and $a>1$
such that for $0\le t_1<t_2<t_3\leq T$ and for all $n$,
    \eqn{
    \lbeq{eqTightness1}
        \expec\Big[|X_n(t_2)-X_n(t_1)|^{r}\,|X_n(t_3)-X_n(t_2)|^{r}\Big]
        \le C |t_3-t_1|^a.
    }
\end{Proposition}

%......................................................

Let
    \eqn{
    V^{\sss(u)}(t) = B_u(tu^{-(\tau-2)})=\sum_{i\in \JJ(u)} c_iu[\indic{T_i\in (u-tu^{-(\tau-2)},u]}-\widetilde \prob(T_i>u-tu^{-(\tau-2)} \mid T_i\leq u)].
    }
We show tightness of $V^{\sss(u)}(t)$ given $u\SS_u=v$. In what follows, we therefore bound
    \eqan{
    \lbeq{tightness-1}
    &\widetilde \expec_v[ ( V^{\sss(u)}(t_2)-V^{\sss(u)}(t_1) )^2 ( V^{\sss(u)}(t_3)-V^{\sss(u)}(t_2) )^2 ] \\
    &= \widetilde \expec\Big[ \widetilde \expec\big[ ( V^{\sss(u)}(t_2)-V^{\sss(u)}(t_1) )^2 ( V^{\sss(u)}(t_3)-V^{\sss(u)}(t_2) )^2 \mid \JJ(u)\big] \mid u\SS_u=v\Big]. \nn
    }
First observe that with
    \eqn{
    \II_i^u(s,t) \equiv \indic{T_i\in (u-t u^{-(\tau-2)},u-s u^{-(\tau-2)}]}
    }
we have $p_{i,s,t}^u = \widetilde \expec[ \II_i^u(s,t) \mid T_i \leq u ]$ (recall \eqref{p-two-times}) and
    \eqan{
    &\widetilde \expec[ ( V^{\sss(u)}(t_2)-V^{\sss(u)}(t_1) )^2 ( V^{\sss(u)}(t_3)-V^{\sss(u)}(t_2) )^2 \mid \JJ(u) ] \\
    &= \widetilde \expec\!\Big[ \prod_{n\in \{1,2\}} \Big( \sum_{i\in \JJ(u)} c_iu \left[ \II_i^u(t_n,t_{n+1}) - p_{i,t_n,t_{n+1}}^u \right] \Big)^2 \mid \JJ(u) \Big]. \nn
    }
By the conditional independence of the processes conditional on $\JJ(u)$ (recall comment preceding \eqref{cond-distr-T}), and as we subtract their respective expectations, we obtain
    \eqan{
    &\widetilde \expec[ ( V^{\sss(u)}(t_2)-V^{\sss(u)}(t_1) )^2 ( V^{\sss(u)}(t_3)-V^{\sss(u)}(t_2) )^2 \mid \JJ(u) ] \\
    &= \widetilde \expec\Big[ \sum_{i\in \JJ(u)} (c_iu)^4 \prod_{n\in \{1,2\}} \left( \II_i^u(t_n,t_{n+1}) - p_{i,t_n,t_{n+1}}^u \right)^2 \mid \JJ(u) \Big] \nn\\
    &\quad + \widetilde \expec\Big[ \sum_{i\in \JJ(u)} \sum_{j\in \JJ(u) \backslash \{i\}} (c_iu)^2 (c_ju)^2 \left( \II_i^u(t_1,t_2) - p_{i,t_1,t_2}^u \right)^2 \left( \II_j^u(t_2,t_3) - p_{j,t_2,t_3}^u \right)^2 \mid \JJ(u) \Big] \nn\\
    &\quad + 2 \widetilde \expec\Big[ \sum_{i\in \JJ(u)} \sum_{j\in \JJ(u) \backslash \{i\}} (c_iu)^2 (c_ju)^2 \prod_{n\in \{1,2\}} \left( \left( \II_i^u(t_n,t_{n+1}) - p_{i,t_n,t_{n+1}}^u \right) \left( \II_j^u(t_n,t_{n+1}) - p_{j,t_n,t_{n+1}}^u \right) \right) \mid \JJ(u) \Big]. \nn
    }
We can bound this from above by
    \eqn{
    C \Big\{ \sum_{i\in \JJ(u)} (c_iu)^4 p_{i,t_1,t_2}^u p_{i,t_2,t_3}^u + \prod_{n\in \{1,2\}} \Big( \sum_{i\in \JJ(u)} (c_iu)^2 p_{i,t_n,t_{n+1}}^u \Big) \Big\}.
    }
By \refeq{p_two_times_bounds},
    \eqan{
    \lbeq{tightness-2}
    &\widetilde \expec[ ( V^{\sss(u)}(t_2)-V^{\sss(u)}(t_1) )^2 ( V^{\sss(u)}(t_3)-V^{\sss(u)}(t_2) )^2 \mid \JJ(u) ] \\
    &\leq C (t_2-t_1) (t_3-t_2) \Big\{\sum_{i\in \JJ(u)} (c_iu)^4 u^{-2(\tau-1)} (\e^{-c_iu} c_i u \wedge 1)^2+\Big(\sum_{i\in \JJ(u)} (c_iu)^2 u^{-(\tau-1)} (\e^{-c_iu} c_i u \wedge 1)\Big)^2 \Big\}. \nn
    }
For the first sum, note that $(c_iu)^4 u^{-2(\tau-1)}
=c_i^4 u^{-2(\tau-3)}$, so that its sum is order $o(1)$ as $\sum_i c_i^3<\infty$
and $\tau>3$. For the second sum in \eqref{tightness-2}, we note that the sum over $i$ such that
$c_i>1/u$ is clearly bounded, since it is bounded by
    \eqn{
    \sum_{i\colon c_i>1/u} (c_iu) u^{-(\tau-1)} \e^{-c_iu},
    }
which converges to a constant as $u\rightarrow \infty$ since it is a Riemann approximation to a finite integral.
For the contributions due to $c_i\leq 1/u$, we bound
its expectation as
	\eqan{
    \lbeq{tightness-3}
	&\widetilde \expec_v\!\Big[ \Big(\sum_{i\in \JJ(u)\colon c_i\leq 1/u} (c_iu)^2 u^{-(\tau-1)} \Big)^2 \Big] \\
    &\leq \sum_{i \neq j, c_i\leq 1/u,c_j\leq 1/u} (c_iu)^2 (c_ju)^2 u^{-2(\tau-1)} c_iu c_ju + \sum_{i\in \JJ(u)\colon c_i\leq 1/u} (c_iu)^4 u^{-2(\tau-1)} c_iu \nn\\
	&\leq \Big( \sum_{i\colon c_i\leq 1/u} (c_iu)^3 u^{-(\tau-1)} \Big)^2 + \sum_{i\colon c_i\leq 1/u} c_i^3 \leq C, \nn
	}
by \refeq{lem-help-(a)-(b)}. Hence, we get with \eqref{tightness-1} and \eqref{tightness-2}-\eqref{tightness-3},
    \eqn{
    \widetilde \expec_v[ ( V^{\sss(u)}(t_2)-V^{\sss(u)}(t_1) )^2( V^{\sss(u)}(t_3)-V^{\sss(u)}(t_2) )^2 ]
    \leq C (t_2-t_1) (t_3-t_2)
    \leq C (t_3-t_1)^2,
    }
as required.
\qed

%%%%%%%%%%%%%%%%%%%%%%%%%%%%%%%%%%%%%%%%%%%%%%%%%%%%%%%

\subsubsection{Completion of the proof of Proposition \ref{prop-weak-cond-conv}(b)}

%%%%%%%%%%%%%%%%%%%%%%%%%%%%%%%%%%%%%%%%%%%%%%%%%%%%%%%

The convergence of the finite-dimensional distributions together with tightness yields $(B_u(tu^{-(\tau-2)}))_{t\geq 0}\convd (L_t)_{t\geq 0}$ by \cite[Theorem 5.1]{Bill99}.

%%%%%%%%%%%%%%%%%%%%%%%%%%%%%%%%%%%%%%%%%%%%%%%%%%%%%%%
%%%%%%%%%%%%%%%%%%%%%%%%%%%%%%%%%%%%%%%%%%%%%%%%%%%%%%%

%%%%%%%%%%%%%%%%%%%%%%%%%%%%%%%

%\input{section8_31_05_2012.tex}

%%%%%%%%%%%%%%%%%%%%%%

\paragraph{Acknowledgements.}
The work of RvdH, JvL and SK  was supported
in part by the Netherlands Organisation for Scientific Research (NWO).
The work of JvL was supported by the European Research Council (ERC).
We thank Elie A\"id\'ekon for numerous discussions.

%\listoftodos

\bibliographystyle{plain}
\bibliography{bib}
%\todo[inline]{S-? check [1 our other article: Preprint \ch{2014}] [7 BHL-Novel scaling limits], [30 Turova]. Remco: all done, except for [30] (which I need to look up).}
%\medskip

%New:

%\begin{itemize}
%    \item
%    Fixed proof of Lemma \ref{lem-error-for-tail}.
%    \item
%    Added Convergence of the finite-dimensional distributions of $B_u$.
%\end{itemize}

%Things to do:
%
%\begin{itemize}
%    \item
%    Add Sandra's references to Bruijn, Durrett, Bertoin. Billingsley-tightness is Theorem 15.6.
%    \item
%    Mention somewhere that $\SS_u$ is $\JJ(u)$-measurable etc.
%    \item
%    Replace $\gamma=(\tau-1)/2$ in appropriate places. (already in use for instance in \eqref{near-end})
%    \item
%    Sandra: Clean up Lemma \ref{lem-expec-S} (e.g. move $1+\beta t$ to order terms)
%\end{itemize}

\end{document}